%% file: olv24.tex
\documentclass[11pt,a4paper]{article}
\usepackage[body={154mm,229mm},top=34mm]{geometry} \usepackage{amscd}
\usepackage{float} \usepackage{bm} \usepackage{mathrsfs} \usepackage{epsfig}
\usepackage{amsmath} \usepackage{amssymb} \usepackage{amsthm}
\usepackage{epsfig} \usepackage{verbatim} \usepackage{url} \usepackage{xcolor}
\usepackage{hyperref} \usepackage{tikz}
\usepackage{pspicture}

 \def\fq{\mathfrak{q}}

 \def\var{\text{var}} 
\def\vol{\text{vol}} \def\Vol{\text{Vol}} \def\sgn{\mbox{sgn }}

\newcommand{\ninseps}[3]{ \begin{figure}[H] \begin{center}
\scalebox{#3}{\includegraphics{#1}} \end{center}

\vspace{-0.65cm} \caption{\hspace{0.25cm}#2\label{f:#1}} \end{figure} }

 \newtheorem{lemma}{Lemma}

\newtheorem{definition}{Definition}[section] \newtheorem{remark}{Remark}

\newtheorem{assumption}{Assumption}
\newtheorem{proposition}{Proposition}[section]

 \usepackage{makeidx}
\usepackage{mathabx} \usepackage{mathrsfs} 

\DeclareRobustCommand{\vect}[1]{\bm{#1}} \newcommand{\Ymin}{Y^{\min}} \newcommand{\Zmin}{Z^{\min}}
\newcommand{\Mmin}{M^{\min}}

\newcommand{\contrset}{{\mathcal A}_{I,\bm S}}

\begin{document} \title{Optimal Liquidation with Conditions on Minimum Price\footnote{This work was supported by TUBITAK (The Scientific and Technological Research Council of Turkey)
through project number 118F163.}}
\author{Mervan Aksu\footnote{Middle East Technical University, Institute of
Applied Mathematics, Ankara, Turkey, and Mardin Artuklu University, Faculty of
Economics and Administrative Sciences, Mardin, Turkey,
mervanaksu@artuklu.edu.tr}, Alexandre Popier\footnote{Laboratoire Manceau de
Math\'ematiques, Le Mans Universit\'e, France, apopier@univ-lemans.fr}, Ali
Devin Sezer\footnote{Middle East Technical University, Institute of Applied
Mathematics, Ankara, Turkey, devinsezer@gmail.com}}

\maketitle \begin{abstract} The classical optimal trading problem is the
	closure of an initial position in a financial asset over a fixed time
	interval; the trader tries to maximize an expected utility under the
	constraint that the position be fully closed by terminal time.  Given
	that the asset price is stochastic, the liquidation constraint may be
	too restrictive; the trader may want to relax the full liquidation
	constraint or slow down/stop trading depending on price behavior.  We
	consider two additional parameters that serve these purposes within the
	Almgren-Chriss liquidation framework: a binary valued process $I$ that
	prescribes when trading takes place and a measurable set ${\bm S}$ that
	prescribes when full liquidation is required.  We give four examples
	for ${\bm S}$ and $I$ which are all based on a lower bound specified
	for the price process.  The terminal cost of the stochastic optimal
	control problem is $\infty$ over ${\bm S}$; this represents the
	liquidation constraint.  The permanent price impact defines the
	negative part of the terminal cost over the complement of ${\bm S}.$
	The $I$ parameter enters the stochastic optimal control problem as a
	multiplier of the running cost.  Except for quadratic liquidation costs
	the problem turns out to be non-convex.  A terminal cost that can take
	negative values implies 1) the backward stochastic differential
	equation (BSDE) associated with the value function of the control
	problem can explode to $-\infty$ backward in time and 2) the existence
	results on minimal supersolutions of BSDE with singular terminal values
	and monotone drivers are not directly applicable.  To tackle these we
	introduce an assumption that balances the market volume process and the
	permanent price impact in the model over the trading horizon.  In the
	quadratic case, assuming only that the noise driving the asset price is
	a martingale, we show that the minimal supersolution of the BSDE gives
	both the value function and the optimal control of the stochastic
	optimal control problem.  A key point in our arguments is that the
	negative part of the terminal value of the BSDE (arising from permanent
	price impact) is bounded.  For the non-quadratic case, we assume a
	Brownian motion driven stochastic volatility model and focus on choices
	of $I$ and ${\bm S}$ that are either Markovian or can be broken into
	Markovian pieces. These assumptions allow us to represent the value
	functions as solutions of PDE or PDE systems. The PDE arguments are
	based on the smoothness of the value functions and do not require
	convexity.

We quantify the financial performance of the resulting liquidation algorithms by the percentage difference between the initial stock price and the average price at which the position is (partially) closed in the time interval $[0,T].$  We note that this difference can be divided into three pieces:  one corresponding to permanent price impact ($A_1$), one corresponding to random fluctuations in the price ($A_2$) and one corresponding to transaction/bid-ask spread costs ($A_3$).  $A_1$ turns out to be a linear function of the portion of the portfolio that is closed;  therefore, its distribution is fully determined by the distribution of that portion. We provide a numerical study of the distribution of the closed portion	and the conditional distributions of $A_2$ and $A_3$ given the closed portion under the assumption that the price process is Brownian for  $I=1$ and a ${\bm S}$ corresponding to a lowerbound on terminal price.

\end{abstract}



\tableofcontents

\section{Introduction}

There is a range of order types available to an investor to close a position on
an asset; the book \cite{gueant2016financial} presents the following:
implementation shortfall (IS), target close and volume weighted average price.
Given a trading horizon $[0,T]$, all of these order types are constrained to
close an initial position at a terminal time $T$.  Mathematically this is
expressed as the constraint $Q_T = 0$, where $Q_t$ denotes the position of the
investor at time $t$ (we assume $Q_0> 0$, i.e., an initial long position, 
for a more brief presentation; everything
below applies to a short position $Q_0 < 0$).
Given that the price process is stochastic, this constraint
can be too restrictive. For example, in IS orders the goal is to close an 
initial position near the initial price $S_0$;
it may happen that the price drops substantially during the trading interval and
the investor holding the position may no longer wish to be strict about closing the position.
The present work studies algorithms that offer this type of flexibility in execution.
We focus on IS type orders targeting the initial price $S_0$ because this is the most basic and most commonly
studied order in the current literature but
similar ideas can be considered for other types of 
orders and pursued in future research.
The work \cite{krus:popi:15} already presents an optimal liquidation algorithm in which the full liquidation
constraint is relaxed; the algorithms presented in the current work extend those in \cite{krus:popi:15}
in two directions: 1) the model in the present work involves a permanent market impact component and 2) 
we allow the model to impose constraints on when trading takes place 
(for example, the model can be set up so that trading is allowed
only when the current price is above a given threshold). 
The resulting stochastic optimal control problem leads to a
backward stochastic differential equation of the form
\begin{equation*}
	dY_t = \vol_t |Y_t|^p + dM_t, \quad Y_T = \xi, \quad p > 1,
\end{equation*}
with a singular terminal condition (i.e., $\xi$ can take the value $\infty$).
When compared to BSDE with singular terminal values arising
from optimal liquidation problems studied in \cite{krus:popi:15} 
and other works on optimal liquidation, 
this BSDE has three new features:
its generator/driver is not monotone in $Y$ (see below for the
precise definition), $\vol$ is allowed to be $0$, and its terminal value can take negative values.
For several versions of the problem with a Markovian formulation we also treat the partial differential equation (PDE)
versions of this equation.
As we discuss in detail below, the new features have nontrivial implications for the analysis of the BSDE and the PDE.

The IS order is commonly formulated as a stochastic optimal control problem
optimizing the expected utility of the cash position that the order generates \cite[Chapter 4]{gueant2016financial}.
Section \ref{ss:acframework} 
presents the stochastic optimal control formulation of the modified IS orders that the present work focuses on.
A great deal of the literature on optimal liquidation, including
the model presented in section \ref{ss:acframework}, is based on a model proposed by Almgren and Chriss in \cite{almg:chri:01}.
This model assumes that the price consists of three pieces: a random fluctuations term $\bar{S}$,
a transaction cost term and a permanent market impact term.
The permanent price impact is a term $\kappa(Q'_t) dt$ added
to $d\bar{S}_t$ where $ Q'_t$ is the time derivative of $Q$ at time $t$;
together with the random fluctuation term they make up the midprice process $S$.
The transaction cost term corresponds to trading commissions 
and the bid-ask spread; 
this is modeled using a so-called execution cost function
$L$. The common choice for $L$ is a power function $L(\rho) = \eta |\rho|^{\hat p}$, $\hat{p} > 1$ (we reserve
the letter $p$ for the H\"{o}lder conjugate of $\hat{p}$, which arises in the solution of the problem).
In most liquidation models the
permanent price impact is assumed to be linear in $Q'_t$ (see \cite{gath:10} for more comments on this assumption).
This and the constraint $Q_T = 0$ of the standard IS order lead to an interesting situation 
for this type of order:
the permanent price impact ends up having no role in the
stochastic optimal control formulation of the standard IS order (see \eqref{e:xtclassical} and \eqref{e:oc0}). 
These assumptions
have another important mathematical implication for the standard IS order: the
position variable can be factored out of the value function and out of the
backward stochastic differential and partial differential equations
that the value function satisfies (see Proposition \ref{prop:constrained_control_pb})- we will refer to this property of the value function
as homogeneity. 
When we relax the constraint $Q_T = 0$ the permanent price impact enters directly into the stochastic
optimal control problem as a part of the terminal cost. 
In the modifications of the IS order 
we treat in the present work we would like to keep the
homogeneity property of the value function as this greatly influences the analysis of the problem. 
This turns out to be possible if the permanent price
price impact function is chosen compatible with the execution cost function $L$ as follows:
\[
	\kappa(Q_t',Q_t) = k Q_t' | Q_t|^{\hat{p}-2}, k > 0,
\]
$\kappa$ remains linear in $Q_t'$; for $\hat{p} > 2$ ($\hat{p} < 2$) decays (increases) with the
size of the remaining position.

As we have already noted, 
the permanent price impact has no role in the
continuous time stochastic optimal control formulation of the standard IS order in the Almgren-Chriss framework.
We see a reflection of this fact in the  optimal liquidation literature that is based on this framework, as follows:
some works introduce a permanent price impact parameter in the model
(see e.g., \cite{fors:kenn:12,schi:scho:09} and \cite[Chapter 3]{gueant2016financial}),
but it ends up having no role in the stochastic optimal control
problem, while others
drop the permanent price impact parameter altogether
from the model
assumptions, see, e.g.,  \cite{anki:jean:krus:13,krus:popi:15}. 
The seminal work \cite{almg:chri:01} includes a permanent price impact parameter in a discrete time framework, which ends up having
a role in the optimal controls but this role disappears as the discrete time step size converges to $0$.
A key feature of the problem and the analysis presented in the current work is how the permanent price
impact enters and changes the analysis when the full liquidation constraint is relaxed.
We discuss this in detail below.

The output of the standard IS order is the cash position generated by the trading algorithm.
When full liquidation is no longer required, i.e., when we don't have the constraint $Q_T =0$, the output
of the trading process at time $T$ will be $(X_T, Q_T)$ where $X_T$ is the cash generated by the trading process
and $Q_T$ is the position remaining in the asset being traded.
There are several choices available to formulate an expected utility maximization problem based on this output.
One option is to use a general utility function whose input is the pair $(X_T,Q_T)$;
a simpler option is to first assign a monetary value $m(Q_T)$ to the position $Q_T$,
use a utility function whose only input is a monetary value and
apply it to $X_T + m(Q_T)$.
In the present work we use the latter approach.
For $m(Q_T)$ the present work focuses on $m(Q_T) = Q_T S_T$;
this is the market value of the position at terminal time
$T$ ignoring trading costs. The technical advantage of this choice is that it preserves the homogeneity of
the value function; a $|q|^{\widehat p}$ term can also be added that 
preserves homogeneity.
For the utility function we focus on the identity function, i.e., we consider the problem
of maximizing expected terminal position; an additional risk term can be included in the final
stochastic optimal control problem with minor modifications to the analysis. We further comment on these
points in the conclusion (Section \ref{s:conc}).

In subsection \ref{ss:mod}
we discuss two ways the stochastic optimal control problem modeling the IS 
order 
can be modified to delay/stop liquidation depending on price behavior 
1) by relaxing the full liquidation constraint if the price is too low (which was first proposed in \cite{krus:popi:15})
2) stopping/pausing trade if the price is too low (proposed in the present work).
These two modifications are parameterized in the stochastic optimal control problem 
by a measurable set ${\bm S}$ and a process $I$ taking values in $\{0,1\}$.
The set ${\bm S}$ prescribes when full liquidation is required 
and enters the stochastic optimal control problem as a part of the terminal cost;
the process $I$ prescribes when trading takes place and enters the stochastic optimal control problem by multiplying the volume process
(see \eqref{e:oc1m2}).
We give four examples for ${\bm S}$ and $I$ in
subsection \ref{ss:mod} which are all based on a lower bound specified for the price process: 
$(I^{(1)},{\bm S}^{(1)})$ puts no constraints on trading, the position is constrained to be fully closed if the closing price is above a given threshold;
$(I^{(2)},{\bm S}^{(2)})$ allows trading until the price goes below a given threshold, 
the position is constrained to be fully closed if the price remains above the given threshold across the whole trading interval;
$(I^{(3)},{\bm S}^{(3)})$ allows trading only when the price is above a given threshold, the position is constrained to be fully closed if the closing price is above the given threshold;
$(I^{(4)},{\bm S}^{(4)})$ is the same as previous except that for trading to restart the price process first has to upcross
a higher threshold.
The resulting stochastic optimal control problem for the modified IS order is given in \eqref{e:oc1m2}. This problem
has the same structure as the one studied in \cite{krus:popi:15} except for the following differences: 
1) inclusion of the permanent market impact in the model implies that the terminal cost can take negative values
2) the presence of the $I_t$ term in the running cost (for $k=0$ and $I_t =1$ the problem in fact reduces to the
one studied in \cite{krus:popi:15}). 
The solution method in \cite{krus:popi:15} is the derivation and the analysis of a BSDE associated with the value function of 
the control problem. This is also one of the solutions approaches we will pursue in the present work.
We next discuss how the above new features of the problem impact 
the associated BSDE and its analysis.

\bigskip

The set ${\bm S}$ (specifying conditions for full liquidation) 
defines the singular component $\infty \cdot {\bm 1}_{{\bm S}}$ of the terminal condition of the BSDE; 
allowing a nontrivial permanent price impact term 
introduces an additional negative term $-\frac{k}{\hat{p}} {\bm 1}_{{\bm S}^c}$ in the terminal condition of the BSDE.
The driver  $y\mapsto f_t(y), y \in {\mathbb R}$ of a BSDE is said to be monotone
if there exists $\chi \in {\mathbb R}$ such that 
$(f(t,y,\omega)-f(t,y',\omega))(y-y') \leq \chi (y-y')^2$
for any $t \geq 0$ and $y,y' \in {\mathbb R}$ almost surely.
The monotonicity of the driver is a key property in establishing the existence of solutions to BSDE.
The work \cite{krus:popi:15} focuses on non-negative terminal costs; a non-negative terminal cost corresponds to  a nonnegative
terminal condition for the associated BSDE. This and the dynamics and cost structure of the control
problem lead in \cite{krus:popi:15} to a BSDE
with the monotone driver $(t,y) \mapsto -\lambda_t y |y|^{p-1}$ (for $y \ge 0$, $y |y|^{p-1} = |y|^p = y^p$).
In the case of the stochastic optimal control problem \eqref{e:oc1m2}, the terminal condition is allowed to take negative values
and this forces us to work with the non-monotone convex driver 
$(t,y) \mapsto I_t \Vol_t |y|^{p}.$
Finally, as already noted,
the process $I$ (specifying when trading is allowed)
enters the BSDE by multiplying the driver of the BSDE.
The resulting BSDE is given in \eqref{e:convex_BSDE} and \eqref{e:xip}.

Note that although the terminal condition \eqref{e:xip} can take negative values, the negative component is bounded
above by a constant. For this reason we focus on terminal conditions with bounded negative parts.
The challenges/new points that arise in the analysis
of \eqref{e:convex_BSDE}, \eqref{e:xip} as a result of these new features are as follows:
the currently available literature doesn't contain
existence results on the minimal supersolution of a BSDE such as \eqref{e:convex_BSDE} with a singular terminal condition \eqref{e:xip} 
and with a driver involving the function $y \mapsto |y|^p.$
Secondly, a terminal condition that can take negative values and
the superlinearity of the driver imply that the solution of the BSDE can explode to $-\infty$ backwards in time. 

\medskip
We handle both of these issues
in Section \ref{s:bsde}
by deriving an apriori lower bound process $z$ on any supersolution
of the BSDE with a terminal condition whose negative part is 
bounded by a constant $K$ ;
$z$ is obtained by solving \eqref{e:convex_BSDE}
with the terminal condition $\xi=-K$ (for which \eqref{e:convex_BSDE} reduces to an
ordinary differential equation).
The interval over which $z$ is defined depends
on model parameters. We introduce an assumption on the permanent price impact
parameter and market volume that guarantees the existence and boundedness of the lower bound process over the interval $[0,T]$ (Assumption \ref{a:vol}).
The lower bound process is increasing in $t$; therefore, the value $z_0$ provides a lower bound on supersolutions on the whole
 interval $[0,T].$ We then deal with the non-monotonicity of the driver
by replacing the portion of $y \mapsto |y|^p$ over $(-\infty,z_0]$
with its linear approximation over the same interval, which results in a monotone driver. Hence existence and uniqueness of the solution for BSDE for integrable terminal condition is obtained in Proposition \ref{thm:solve_aux_BSDE}. 
The only way the $I$ term impacts the analysis of the BSDE is by entering Assumption \ref{as:singular} that guarantees that the supersolution
can attain the singular value $\infty$ at terminal time $T$; 
in particular Assumption \ref{as:singular} requires $I_t =1$ for $t \in [T-\epsilon,T]$ for some $\epsilon > 0.$ Under these conditions, Proposition \ref{prop:sovlv_sing_BSDE} provides a minimal supersolution for the BSDE with singular terminal condition.

The verification argument presented in subsection \ref{ss:verification}
connects the minimal supersolution of the BSDE \eqref{e:convex_BSDE},
\eqref{e:termcond} to the value function and optimal control
of the stochastic optimal control problem \eqref{e:oc1m2} (Proposition \ref{prop:constrained_control_pb}).
This argument is based on the convexity of the cost in $Q'$; 
under Assumption \ref{as:singular}, this turns out to be the case only for the quadratic case (Lemma \ref{lem:convex_cost}),
therefore this verification argument assumes $\widehat p =p=2.$

Assumption \ref{as:singular} is a stringent condition and $I^{(j)}$, $j=2,3,4$, proposed in subsection \ref{ss:mod} (see \eqref{e:I2}, \eqref{e:stopbelowL} and \eqref{e:stopbelowLb}) which depend on price behavior don't satisfy it. In subsections \ref{ss:red} and \ref{ssect:cases_3_4}, we break up the stochastic optimal control problem and the BSDE corresponding to these choices of $I$ into pieces where each piece involving a singular terminal condition satisfies Assumption \ref{as:singular} (Proposition \ref{prop:control_pb_stop_time} and Lemma \ref{lem:reduc_time_interval}). 

The main advantage of the BSDE arguments in Section 
\ref{s:bsde} is that we can work with a general filtration, without any further assumption except the standard ones (completeness and right-continuity) and the left-continuity at time $T$. The drawback is the restriction to the quadratic case $\widehat{p}=p=2.$

\bigskip

We call the problem Markovian if the price process is Markovian and the cost structure of the problem is a function of the price process. 
In Section \ref{s:pde} we explore the case $ p \neq 2$ and a PDE representation of the value function and the optimal
control when the problem is Markovian or can be broken into Markovian
pieces.
A popular choice for price dynamics in finance applications is the stochastic volatility model. To the best of our knowledge it is rarely treated in the context of optimal liquidation; in Section \ref{s:pde} we assume the price dynamics to follow this model.
For a direct PDE representation (i.e., identifying the value function as a supersolution of a related PDE), the process $I$ and the measurable set ${\bm S}$ must also be functions of the price process.  Of the four possible choices for $I$ and ${\bm S}$ given in \eqref{e:I1}-\eqref{e:stopbelowLb}, only $(I^{(1)}, {\bm S}^{(1)})$ is given as a function of the price process.  In subsection \ref{ss:I1pde} we present the PDE representation of this case. For $(I^{(j)},{\bm S}^{(j)})$, $j=2,3,4$, the decompositions/reductions given
in subsections \ref{ss:red} and \ref{ssect:cases_3_4} yield Markovian subproblems. The resulting PDE problems are treated in
subsections \ref{ss:pdexi2}-\ref{ss:pdeI3}.
The PDE representation of $(I^{(4)}, {\bm S}^{(4)})$ involves two coupled PDE (one for active trading and one when trading
is paused). To compute a solution we introduce an additional parameter $n$, which is the number of switches allowed between 
trading and no trading. We solve this problem recursively, letting $\mathcal N \rightarrow \infty$ gives the solution to the coupled PDE.
To the best of our knowledge, there exists no readily available results in the current literature
for the existence and smoothness of solutions of
PDE that arise in the analysis presented in Section \ref{s:pde}. We obtain solutions to these PDE as follows: we use the BSDE
results of Section \ref{s:bsde} to first obtain viscosity supersolutions. We then use regularization bootstrapping,
parameter smoothing and the regularity of the underlying price process to obtain the smoothness of these supersolutions. 
Once smoothness is proved classical verification
arguments from stochastic optimal control can be constructed to relate the solutions to the stochastic optimal control problems.
These arguments do not require $p =2$ as opposed to the verification arguments given in Section \ref{s:bsde} which are directly
based on the BSDE representation.

\bigskip

The output of the standard IS order is the cash position $X_T$ at time $T$ generated by the trading algorithm; 
under the assumption
that the price process is a Brownian motion, $X_T$ turns out to be normally distributed whose mean and variance have simple formulas
in terms of the model parameters.
When we relax the  IS order so that full liquidation is no longer required at terminal time, the output of the IS order
consists of the pair of real random variables $(Q_T,X_T)$ where $X_T$ is, as before, the total cash generated by the trading process and $Q_T$
is the remaining position at terminal time  in the asset being traded. For the relaxed/modified IS orders,
$X_T$ is not normally distributed even when the price process is taken to be Brownian and the joint distribution
of $(X_T,Q_T)$ doesn't have an explicit form. 
Define 
\begin{equation}\label{d:Sp}
	A = \frac{X_T - (Q_0-Q_T) S_0}{(Q_0-Q_T) S_0}
\end{equation}
$A$ is the percentage deviation from the target price $S_0$
of the average price at which
the position is (partially) closed in the time interval $[0,T].$ 
In Section \ref{s:numerical} we study the joint distribution of $(Q_T/Q_0,A)$. 
We note that $A$ can be divided into three pieces: 
one corresponding
to permanent price impact ($A_1$), one corresponding to transaction/bid-ask spread costs ($A_2$) and one corresponding to random fluctuations in the price ($A_3$). 
$A_1$ turns out to be a linear function of $1-Q_T/Q_0$; therefore, its distribution is fully determined by that of $Q_T/Q_0$. We provide a numerical
study of the distribution of $Q_T/Q_0$ and the conditional distributions of $A_2$ and $A_3$ given $Q_T/Q_0$ under the assumption
that $\bar{S}_t = \sigma W_t$ for the case $(I = 1, {\bm S} = \{W_T \geq \ell \})$.
The same section also provides numerical examples of the sample path behavior of the optimal controls of this modified IS order.
Section \ref{s:conc} comments further on the models presented in this work and on possible future research.

\section{Definitions}\label{ss:acframework}
The following model is based on
the Almgren Chriss framework for liquidation with price impact
(see, e.g., \cite[Chapter 3]{gueant2016financial}).
Everything is assumed to be defined on a probability space $(\Omega, {\mathcal F}, {\mathbb P})$ equipped with 
a filtration ${\mathbb F} = \{ {\mathcal F}_t, t \in [0,\infty)\}$, which satisfies the usual assumptions: completeness and right-continuity.
The market volume at time $t$ is denoted by $\Vol_t$, which is a positive process adapted to the filtration ${\mathbb F}.$ 
The initial position is denoted by $q_0 > 0.$
The position of the investor at time $t$ is $Q_t$, in particular $Q_0 = q_0.$
The process $Q$ is assumed to be absolutely continuous in the time variable, 
let $Q'_t$ denote its derivative at time $t$; $Q$ and $Q'$ are adapted to ${\mathbb F}.$

We will be working with two positive real numbers $p$ and $\widehat{p}$ that
are H\"{o}lder conjugates of each other $1/p + 1/\widehat{p} =1$; we use
$\hat{p}$ in the problem formulation, $p$ appears in the associated
backward stochastic and partial differential equations.

We suppose that the permanent price impact function $\kappa$ is 
\begin{equation}\label{e:kappavp}
	\kappa:{\mathbb  R}^2 \mapsto {\mathbb R},~~~ \kappa(q',q) = k \left| q\right|^{\widehat p-2} q';
\end{equation}
where $k >0$ is a real constant.
For $ \widehat{p} > 2$ ($\widehat{p} < 2$), $\kappa$ decreases (increases) with position size and for $\widehat p=2$
it is independent of position size.
For $\widehat{p}=2$, $\kappa$ reduces to 
$k q'$ which is the standard
choice for permanent price impact (see \cite{gath:10} or \cite[Chapter 3]{gueant2016financial}). 
The midprice process $S$ is
\begin{align}\label{e:midprice}
	S_t &=S_0 + \bar{S}_t +  \int_0^t \kappa(Q_s',Q_s)ds
\end{align}
where $\bar{S}$ is a martingale adapted to ${\mathbb F}$.

The actual trading price at time $t$ is
	$S_t + g_t(Q'_t/\Vol_t)$
where $g_t$ models transaction costs and the bid-ask spread ($g$ depends on $t$, $\omega$ and
$Q'_t/\Vol_t$).
The process $g$ is often specified via the so-called execution cost function 
$l_t$: $l_t(\rho) = \rho g_t(\rho)$. 
The actual trading price at time $t$, expressed in terms of $l$ is
\[
	S_t + \frac{\Vol_t}{Q_t'} L_t\left(\frac{Q_t'}{\Vol_t}\right).
\]
The cash position that $Q$ generates is
\begin{equation}
X_T = -\int_0^T \left(S_t + \frac{\Vol_t}{Q_t'} L_t\left(\frac{Q_t'}{\Vol_t}\right)\right) Q_t' dt
    = -\int_0^T S_t Q_t' dt -\int_0^T  \Vol_t L_t\left(\frac{Q_t'}{\Vol_t}\right)  dt
\label{e:XTp}
\end{equation}
($Q_t' < 0$ corresponds to selling, hence an increase in $X$, and $Q_t' > 0$ corresponds to buying, a decrease in $X$).
Note that $d (|x|^{\widehat p} )/dx = {\widehat p} x |x|^{{\widehat p} -2}$; this
and integrating the first term by parts give:
\begin{align*}
\int_0^T S_t Q_t' dt  &= S_TQ_T - S_0Q_0 -\int_0^T Q_t dS_t\\
	&= S_TQ_T - S_0Q_0 -\int_0^T Q_t d\bar{S}_t - \int_0^T k Q_t Q_t'  |Q_t|^{{\widehat p}-2}  dt\\
	&= S_TQ_T - S_0Q_0 -\int_0^T Q_t d\bar{S}_t - \frac{k}{\widehat{p}}
	( |Q_T|^{\widehat p} - |Q_0|^{\widehat p})
\end{align*}
Then
\begin{equation*}
	X_T = Q_0 S_0 -Q_T S_T +K
	\left( |Q_T|^{\widehat{p}} -|Q_0|^{\widehat p} \right) + \int_0^T Q_t d\bar{S}_t -\int_0^T \Vol_t L_t\left(\frac{Q'_t}{\Vol_t} \right) dt,
\end{equation*}
where
\[
	K = \frac{k}{\widehat{p}}.
\]
Let us first review the classical IS order, in which 
the position is required to be closed fully at terminal time,
i.e., we impose the constraint $Q_T = 0$ on the problem.
The terminal cash position under this constraint is
\begin{equation}\label{e:xtclassical}
	X_T = 
	Q_0 S_0 -K Q_0^{\widehat{p}}+ \int_0^T Q_t d\bar{S}_t -
	\int_0^T \Vol_t L\left(\frac{Q'_t}{\Vol_t} \right) dt.
\end{equation}
To identify the optimal liquidation strategy $Q^*$ one maximizes the expected utility ${\mathbb E}[U(X_T)]$, over the admissible strategies $Q$, where $U$ is the utility function of the trader. The standard choice for the utility function in optimal liquidation
literature is 
\begin{equation}\label{e:U}
	U(x) = -e^{-\gamma x},
\end{equation}
where $\gamma$ is the risk aversion
parameter of the investor \cite{gueant2016financial}. 
In the current work we will be focusing on the case $\gamma \rightarrow 0$, for which the problem becomes\footnote{Using $(-e^{-\gamma x}+1)/\gamma \approx x $ when $\gamma$ goes to zero.}
\begin{equation}\label{e:oc0}
	\sup_{Q \in {\mathcal A}_{1,\Omega} } {\mathbb E}[X_T] ,
	\end{equation}
with
\begin{align*}
	{\mathcal A}_{1,\Omega} = &\Bigg\{Q: \text{absolutely continuous, $Q'$ progressively measurable},   \\  
	&~~{\mathbb E}\left[
		\int_0^T \Vol_t L\left(\frac{Q'_t}{\Vol_t} \right) dt\right] <  +\infty \mbox{ and } Q_0 =q_0, Q_T =0 \ \text{a.s.} \Bigg\},
\end{align*}
see Definition \ref{d:efAIS} below. 
The term $\int_0^\cdot Q_t d\bar{S}_t $ 
is a martingale\footnote{In general it is only a local martingale. However we will show that for the optimal strategy, $Q$ is bounded. 
Hence there is no harm in assuming that this process is indeed a martingale.} ;
taking the expectation of $X_T$ in \eqref{e:xtclassical},
\eqref{e:oc0} reduces to
\begin{equation*}
	\inf_{Q \in {\mathcal A}_{1,\Omega} } 
	{\mathbb E}\left[ \int_0^T \Vol_t L\left(\frac{Q_t'}{\Vol_t} \right) dt\right] .
\end{equation*}
This is the standard version of the stochastic optimal control
formulation of the IS order in the Almgren Chriss framework
for $\gamma = 0.$ 
This problem and its generalization where $\gamma > 0$ is a well studied problem in the current literature, as in \cite{guea:15,schi:scho:09,schi:scho:tehr:10}. 
In the next subsection we propose several modifications to this problem relaxing the
constraint $Q_T = 0$ and introducing constraints on when trading takes place.

\subsection{Modifications}\label{ss:mod}

When full liquidation is no longer required, i.e., when we don't have the constraint $Q_T =0$, the output
of the trading process at time $T$ will be $(X_T, Q_T)$ 
where $X_T$ is the cash generated by the trading process
and $Q_T$ is the position remaining in the asset being traded. 
As discussed in the introduction, 
to formulate a utility maximization problem similar to 
\eqref{e:oc0},
we assign a monetary value $m(Q_T)$ to the position $Q_T$ and add it to $X_T$.
In the present work we mainly focus on the following simple choice for $m$: the market value of the position at terminal time
$T$ ignoring trading costs, i.e., $Q_T S_T$.
With this choice, the monetary value of the position $(X_T, Q_T)$ is
\begin{equation}\label{e:XTm}
	\tilde{X}_T = X_T + Q_T S_T = Q_0 S_0 -K|Q_0|^{\widehat{p}} + 
	K|Q_T|^{\widehat{p}}
	+ \int_0^T Q_t d\bar{S}_t -
	\int_0^T \Vol_t L_t\left(\frac{Q'_t}{\Vol_t} \right) dt.
\end{equation}

Recall that our goal is to modify 
the IS order to not liquidate depending on price behavior in two ways 1) by relaxing
the full liquidation constraint if the price is too low 2) stopping/pausing trade if the price is too low.
The following formulation allows both of these possibilities. Let $t\mapsto I_t$ be an adapted process taking values in 
$\{0,1\}.$ Let $\bm S \in {\mathcal F}_T$ be a measurable set. The trading set is defined by
$$\mathcal I = \{ t \in [0,T] ,\ I_t = 1\}.$$
\begin{definition}\label{d:efAIS}
Define $\contrset$ as the set of processes $Q$ that satisfy:
\begin{itemize}
	\item $Q$ is absolutely continuous in $t$ and $t\mapsto Q'_t$ is progressively measurable ;
	\item $Q_0 = q_0$ ;
\item $Q_T(\omega) =0$ if $\omega \in \bm S$ (liquidation constraint) ;
\item for $\lambda\otimes \mathbb P$-almost every $(t,\omega)\in [0,T] \times \Omega$, $I_t(\omega) = 0 \Rightarrow Q'_t(\omega) = 0$ (if $t\not \in \mathcal I$, $Q'_t = 0$) ;
\item the cost 
	\[
		 \int_0^T {\rm Vol}_t L\left(\frac{Q'_t}{{\rm Vol}_t} \right)dt 
		 -K
		|Q_T|^{\widehat{p}}
	\]
belongs to $L^\varrho(\Omega)$ for some $\varrho > 1$. 
\end{itemize}
\end{definition}
For ease of notation set
\begin{equation}\label{e:conventionA}
	{\mathcal A} = {\mathcal A}_{1,\emptyset},
\end{equation}
i.e., when $I=1$ and ${\bm S} = \emptyset$ we simply write ${\mathcal A}$ instead of ${\mathcal A}_{I,{\bm S}}.$

We modify \eqref{e:oc0} to 
\begin{equation*}
	\sup_{Q \in \contrset} {\mathbb E}[\tilde{X}_T].
\end{equation*}
The formula \eqref{e:XTm} implies that this control problem is equivalent to:
\begin{equation*}
	\inf_{Q \in \contrset} {\mathbb E}\left[ \int_0^T \Vol_t L\left(\frac{Q'_t}{\Vol_t} \right)dt -
	K|Q_T|^{\widehat{p}}\right]
\end{equation*}
As noted in the introduction we work with the
execution cost function \[
	L_t (\rho) = \eta_t |\rho|^{\widehat{p}},
\]
where $\eta$ is an adapted and positive valued process.
This choice of $L$ 
reduces the problem to
\begin{equation}\label{e:oc11m}
	\inf_{Q \in \contrset} {\mathbb E}\left[ \int_0^T  \dfrac{\eta_t}{\Vol_t^{\widehat{p}-1}} |Q'_t|^{\widehat{p}} dt 
	-K |Q_T|^{\widehat{p}}\right].
\end{equation}
Let us next note that $\eta_t$ can always be assumed to be $1$ by appropriately modifying the volume process, i.e.,
if we set 
\begin{equation}\label{e:Voltilde}
	\widetilde{\Vol}_t= \frac{ \Vol_t }{ \eta_t^{\frac{1}{{\widehat{p}}-1}}} = \frac{ \Vol_t }{ \eta_t^{p-1}} 
\end{equation}
\eqref{e:oc11m} can be written as
\[
	\inf_{Q \in \contrset} {\mathbb E}\left[ \int_0^T  \dfrac{1}{\widetilde\Vol_t^{\widehat{p}-1}} |Q_t'|^{\widehat{p}} dt 
	-K |Q_T|^{\widehat{p}}\right];
\]
in what follows we always assume $\eta_t =1$ and that $\Vol$ is modified 
to $\widetilde{\Vol}$
if the original $\eta$ process is not identically $1$.
A commonly used convention in the prior literature is
\begin{equation}\label{e:infconv}
	\infty \cdot 0 = 0, \infty \cdot c =  \begin{cases} &\infty, \text{ if }
		c > 0, \\
		&-\infty, \text{ if } c < 0.
	\end{cases}
\end{equation}
In addition to this we will also set $0/0 = 1$ and $c/0=  \infty$ for $c > 0.$
With these conventions \eqref{e:oc11m} can be written as
\begin{equation}\label{e:oc1m2}
	\inf_{Q \in {\mathcal A} } {\mathbb E}\left[ 
	\int_0^T \dfrac{1}{I_t \Vol_t^{\widehat{p}-1}}  |Q_t'|^{\widehat{p}} dt 
	+\left(-K{\bm 1}_{{\bm S}^c} + \infty
	\cdot {\bm 1}_{\bm S} \right)|Q_T|^{\widehat{p}}\right],
\end{equation}
where we use the convention \eqref{e:conventionA}.
Note that the $I$ process controls when trading takes place 
(by effectively switching the volume process on and off)
and the event ${\bm S}$ controls when full liquidation is required.

In the next section we obtain a representation of the value function
and the optimal control of the problem \eqref{e:oc1m2} via the
minimal supersolution of an associated BSDE. Before that let us
give several examples for the process $I$ and the event ${\bm S}.$
The midprice process $S$ consists of two components:
$\bar{S}_t$ and $\int_0^\cdot \kappa(Q'_s,Q_s)ds .$ 
The first component,
$\bar{S}$, is the random component of the change in
the midprice; a large and unpredictable 
drop in price that the investor may fear can only arise from
this component. Given this observation, a reasonable approach in 
choosing $I$ and ${\bm S}$ is
by putting a lower bound $\ell$ on this component. For this, define:
\begin{align*}
	\tau_\ell &= \inf\{ s \ge 0: \bar{S}_s < \ell \},\\
	\tau_{t,\ell} &= \inf\{ s \ge t: \bar{S}_s < \ell \}.
\end{align*}
Then some possible choices for $I$ and ${\bm S}$ are:
\begin{align}
	I_t &=I_t^{(1)} \doteq 1, {\bm S} ={\bm S}^{(1)} \doteq \{\bar{S}_T \geq \ell\}:
\label{e:I1}
	\intertext{ trading is allowed at all times, 
		full liquidation is forced only when the terminal price $\bar{S}_T$ is above $\ell$.}
			\label{e:I2}
			I_t &= I_t^{(2)} \doteq  {\bm 1}_{\{ t \le \tau_\ell \}}, {\bm S} = {\bm S}^{(2)} \doteq  \{\tau_\ell > T \}:
			\intertext{
		trading stops once $\bar{S}$ hits
		the lower bound $\ell$; full liquidation takes place if $\bar{S}$ remains above $\ell$ throughout $[0,T].$}
			\label{e:stopbelowL}
	I_t &= I_t^{(3)} \doteq {\bm 1}_{[\ell,\infty)}(\bar{S}_t) {\bm 1}_{\{t \le \tau_{T-\delta,\ell}\}}, {\bm S} = {\bm S}^{(3)} = \{\tau_{T-\delta,\ell} \ge T \}:
\end{align}
trading pauses when the
price $\bar{S}$ is below $\ell$, full liquidation takes place if the price process $\bar{S}$ remains above $\ell$
in the time interval $[T-\delta,T]$ and trading stops if $\bar{S}$ goes below $\ell$ in the same interval; $\delta >0$ is a small
fixed constant.

Let us comment on the $\delta >0$ parameter in this formulation: essentially we would like to continue with the liquidation
when the price is not too below our target price $S_0$ and close the position fully if the terminal price is also near our target price.
However, allowing trading (re)start 
arbitrarily close to $T$ and forcing a full liquidation implies high transaction costs (in fact, $\infty$ transaction costs under
the current model). This is the reason
for the $\delta > 0$ parameter: full liquidation is forced only if the price remains above $\ell$ in the time interval $[T-\delta,T].$

In the last formulation trading pauses once $\bar{S}$ hits $\ell$; if $\bar{S}$ is a continuous diffusion process,
once it hits $\ell$, it will hit $\ell$ infinitely often and the trading process will switch on and off infinitely often
as $\bar{S}$ crosses $\ell$. One can get a discrete sequence of on and off trading intervals by putting a buffer of
size $ b> 0$ above $\ell$ between trading and no trading; 
once trading pauses, it is turned back on once $\bar{S}$ goes above $b+\ell.$ The corresponding
$I$ and $S$ are expressed through the following sequence of hitting times:
\begin{align*}
\tau_{\ell,0} & = \tau_\ell, \quad \tau_{b,-1} = 0\\
\tau_{b,0} &= \inf \{t: t > \tau_{\ell,0}, \bar{S}_t \ge  b+\ell  \}, \\
\tau_{\ell,k} &= \inf \{t: t > \tau_{b,k-1}, \bar{S}_t <  \ell \},\\
\tau_{b,k} &= \inf \{t: t > \tau_{\ell,k}, \bar{S}_t \ge b+\ell\}.\\
	\bar\tau_{b,k} &= \begin{cases} \tau_{b,k}, &\text{ if } \tau_{b,k} + \delta < T,\\
                                T,         &  \text{otherwise}.
        \end{cases}
\end{align*}
Adding a buffer of size $b>0$ between no-trading and trading in \eqref{e:stopbelowL} amounts to the following
definitions:
\begin{equation}\label{e:stopbelowLb}
I_t = I_t^{(4)} \doteq \sum_{k=-1}^\infty {\bm 1}_{[ \tau_{b,k} , \tau_{\ell,k+1} ]}(t),~~~ {\bm S} = {\bm S}^{(4)} \doteq \{I_T = 1\}.
\end{equation}

\section{BSDE Analysis}\label{s:bsde}

In our arguments the concepts of monotonicity (of the driver of a BSDE)  and the minimal supersolution of a BSDE play
a key role, let us begin by giving a precise definition of these terms.
\begin{definition}
The function $(t,y,\omega) \mapsto f(t,y,\omega)$ is said to be monotone 
	if there exists $\chi \in {\mathbb R}$ such that a.s. and for any $t \geq 0$ and $y,y' \in {\mathbb R}$,
	\begin{equation*}
(f(t,y,\omega)-f(t,y',\omega))(y-y') \leq \chi (y-y')^2.
\end{equation*}
Let  $\xi$ be an ${\mathcal F}_T$-measurable real valued random variable.
A pair $(Y,M)$ is said to be a supersolution of the BSDE 
	\[
		dY_t =- f(t,Y_t) dt + dM_t, \quad Y_T = \xi,
	\]
	if \begin{enumerate}
		\item $Y$ is adapted to $\mathbb F$ ;
		\item $M$ is a martingale\footnote{In this paper, we always consider a right-continuous with left limits modification of the martingale} ;
		\item  $Y_s = Y_t  +\int_s^t f(u,Y_u) du + (M_t - M_s)$ for $ 0 \le  s < t < T$ ;
		\item the inequality 
\begin{equation}\label{e:minimalitycond}
	\liminf_{t\rightarrow T} Y_t \ge \xi
\end{equation}
holds a.s..
\end{enumerate}
It is said to be minimal if a.s. for any $t$, $Y'_t \ge Y_t$ for any other supersolution $(Y',M').$ 
\end{definition}
The work \cite{krus:popi:15} studies the following stochastic optimal control
problem:
\begin{equation}\label{e:sockp}
	Q_t = Q^C_t +  Q^J_t,\quad 
	\min_{Q \in {\mathcal A}'} {\mathbb E}\left[
		\int_0^T \left( \eta_t \left| \dfrac{dQ^C_t}{dt}  \right|^{\hat{p}}  
		+ \gamma_t |Q_t|^{\hat{p}} + C(Q^J)_t \right)dt
		 + 
		\xi |Q_T|^{\widehat{p}} \right],
\end{equation}
where: $\xi \in {\mathcal F}_T$ is a non-negative random variable that is allowed to take the value $+\infty$, $Q^C$ is the absolutely continuous part of $Q$, $Q^J$ is the jumping part of $Q$, $C(Q^J)$ a running cost associated
with $Q^J$, and ${\mathcal A}'$ an appropriate modification of ${\mathcal A}.$
On the set $\{\xi= +\infty\}$, the constraint $Q_T = 0$ is necessary to ensure a finite cost. 
Compared to \eqref{e:oc1m2} this problem has an additional term $Q^J$ in
its dynamics and two additional terms ($\gamma_t |Q_t|^{\hat{p}}$ and
$C(Q^J)_t$) in its cost structure.  To focus on the novelties associated
with the terminal cost we will assume these
terms to be $0$.
The work \cite{krus:popi:15} identifies the value function of \eqref{e:sockp} as $Y_t |Q_t|^{\hat{p}}$ where $Y$ is the minimal supersolution
of the BSDE
\begin{equation}\label{e:bsdepk}
	dY_t = (\hat p-1)\frac{Y_t |Y_t|^{p-1}}{\eta^{p-1}_t} +dM_t, \qquad Y_T = \xi \ge 0.
\end{equation}
The generator of this BSDE is monotone and establishing the existence of the minimal supersolution makes use of this property of the generator.
The reason that \cite{krus:popi:15} is able to use this monotone generator is the assumption $\xi \ge 0$.
In the present work we allow $\xi$ to take negative values therefore it is no longer possible to work with a generator involving the function $y\mapsto y |y|^{p-1}$.
The BSDE corresponding to the stochastic optimal control problem 
\eqref{e:oc1m2} turns out to be
\begin{equation} \label{e:first_convex_BSDE}
	dY_t = (\hat p-1) I_t \Vol_t |Y_t|^p dt + dM_t
\end{equation}
with terminal condition
\begin{equation}\label{e:termcond}
Y_T = \xi,  
\end{equation}
\begin{equation}\label{e:xip}
 \xi  = -K {\mathbf 1}_{{\bm S}^c} +  \infty \cdot  {\mathbf 1}_{{\bm S}}
\end{equation}
(a rigorous link between this BSDE and \eqref{e:oc1m2} will be established
in subsection \ref{ss:verification} via a verification argument).

Define
\begin{equation}\label{e:vol}
	\vol_t = (\hat p-1) I_t \Vol_t ;
\end{equation}
$\vol$ is the restricted volume process available to the trader (up to some constant). Then we can write \eqref{e:first_convex_BSDE}
also as
\begin{equation} \label{e:convex_BSDE}
	dY_t =  \vol_t |Y_t|^p dt + dM_t.
\end{equation}
Compared to \eqref{e:bsdepk} the novel features of \eqref{e:convex_BSDE}, \eqref{e:termcond}
are the following: 
its terminal condition is allowed to take negative values, its driver
is convex in $Y$ and not monotone and its generator can take the value $0$ (because of the presence
of the $I$ term). 
We introduce several assumptions
to deal with these new features when obtaining the existence of a minimal supersolution to this BSDE.
First, note that the terminal condition $\xi$ of \eqref{e:xip} is bounded 
below by $-K$; for this reason
for our purposes it suffices to focus on terminal conditions whose negative 
parts are bounded above by a constant, i.e., all of the terminal conditions
we consider satisfy
\[
	\xi^- \le K.
\]

The  generator of \eqref{e:convex_BSDE} is defined in terms of 
the function $y \mapsto |y|^{p}$
and its terminal value can take negative
values: a consequence of these facts is 
that any solution to \eqref{e:convex_BSDE} and \eqref{e:termcond} can explode to
$-\infty$ backward in time (see Lemma \ref{lem:lower_bound} below). To deal with this, we introduce
the following assumption that ensures that an explosion doesn't happen 
in $[0,T]$:
\begin{assumption}\label{a:vol}
${\rm vol}$ is non-negative and one of the next two conditions holds:
\begin{itemize}
\item ${\rm vol}$ is deterministic and satisfies
\begin{equation}\label{e:a1v1}
K^{p-1} \int_0^T \mathrm{vol}_t dt <\dfrac{1}{p-1},
\end{equation}
\item ${\rm vol}$ is bounded by a constant $\overline{{\rm vol}} > 0$ such that
\begin{equation}\label{e:a1v2}
K^{p-1}  T\overline{{\rm vol}}  < \dfrac{1}{p-1}.
\end{equation}
\end{itemize}

\end{assumption}
We will assume throughout that either \eqref{e:a1v1} or \eqref{e:a1v2} holds. 
Note that when $(\hat p-1)\Vol$ satisfies
one of these, $\vol = (\hat p-1)I~\Vol$ also does because $I \in \{0,1\}.$
This assumption balances the negative part of the terminal condition (determined
by the permanent price impact parameter $k$) with the trading volume
available to the trader.

If ${\mathbb P}({\bm S}) > 0$,
 $\xi$ can take the value $\infty$ with positive probability and
 the terminal condition is said to be singular (terminal condition is said to be singular also when
 $\xi$ doesn't belong to $L^\varrho (\Omega)$ for some $\varrho$).
We need a further assumption to deal with this possible singularity.
To ensure that there exists a solution which is finite on $[0,T)$, the generator should not be equal to zero close to time $T$. 
The corresponding assumption in  \cite{krus:popi:15} is \cite[(A.6)]{krus:popi:15}, which is: 
$\mathbb E \int_0^T \eta_s ds < +\infty$. In the present work,
to guarantee the existence of a minimal supersolution, we make the following assumption on $\vol = I~ \Vol$:
\begin{assumption}\label{as:singular}
There exists some $\varsigma > 1$ and some $\epsilon > 0$ such that 
\begin{equation} \label{eq:cond_eta_sing}
	\mathbb E \left[ \int_{T-\epsilon}^T \dfrac{1}{(\mathrm{vol}_s)^{(\widehat p-1)\varsigma }} ds \right] <+\infty. 
\end{equation}
Evoke that $\widehat p$ is the H\"older conjugate of $p$.
\end{assumption}
Assumption \ref{as:singular} can be interpreted as the availability of liquidity (through $\Vol$ and $I$)
at terminal time. In particular it means that $I_t = 1$ on $[T-\epsilon,T]$.

We begin our analysis by deriving a lower bound process $z$ which we will
use to guarantee that the minimal supersolution to \eqref{e:bsdepk}
doesn't explode to $-\infty$ in $[0,T]$.
Under Assumption \ref{a:vol} the lower bound process $z$ is defined as follows:
\begin{equation}\label{e:z1}
z_{t} = -
 \left( \dfrac{1}{ K^{1-p} - (p-1)\int_t^T \mathrm{vol}_s ds } \right)^{ \frac{1}{p-1}},
\end{equation}
if \eqref{e:a1v1} holds;
and
\begin{equation}\label{e:z2}
z_{t} = -
 \left( \dfrac{1}{ K^{1-p} - (p-1)\overline{\mathrm{vol}}(T-t) } \right)^{ \frac{1}{p-1}},
\end{equation}
if \eqref{e:a1v2} holds.
\begin{lemma} \label{lem:lower_bound}
	$z$ of \eqref{e:z1} satisfies
	\begin{equation}\label{e:ode}
\frac{dz}{dt} - \mathrm{vol}_t |z|^p = 0,
	\end{equation}
	and $z$ of  \eqref{e:z2} satisfies
	\begin{equation}\label{e:ode2}
	 \dfrac{dz}{dt} - \overline{\mathrm{vol}} \ |z|^p=0.
	\end{equation}
 Both $z$ satisfy $z_T = -K$.
	Under Assumption \ref{a:vol},
	$z$ is increasing on $[0,T]$ and satisfies
	$ -\infty < z_0 \le  z_t \le -K$ 
	for any $t \in [0,T]$.
\end{lemma}
\begin{proof}
	Assumption \ref{a:vol} implies
\begin{equation*}
K^{1-p} - (p-1) \int_t^T \vol_s ds > 0\quad \mbox{or}\quad K^{1-p} - (p-1) \overline{\vol} (T-t) > 0,
\end{equation*}
	for $t \in [0,T].$ Therefore, $z_t < 0$ on $[0,T]$. Non-negativity of $\vol$ and $\overline{\vol}$ imply that $z$ is increasing.
	One can check by differentiation that $z$ of \eqref{e:z1} satisfies \eqref{e:ode} and $z$ of \eqref{e:z2} satisfies \eqref{e:ode2}.
\end{proof}

The standard way to obtain the minimal supersolution of a BSDE 
with a singular terminal condition 
is approximation from below, i.e., we truncate the terminal condition $\xi$ to $\xi \wedge n$,
solve the resulting BSDE and let $n \nearrow \infty$. 
Therefore, the treatment of singular terminal values requires the solution of the same BSDE with bounded/integrable terminal values.
The next proposition addresses such terminal values: 
\begin{proposition} \label{thm:solve_aux_BSDE} 
Suppose that Assumption \ref{a:vol} holds. Furthermore assume:
\begin{itemize}
\item $\xi^+ \in L^\varrho(\Omega)$ for some $\varrho >1$,
\item $\xi^- \le K$ where $K$ is the constant appearing in Assumption \ref{a:vol}.
\end{itemize}
	Then BSDE \eqref{e:convex_BSDE} 
has a unique solution $(Y,M)$ such that $Y^-$ is bounded and 
$$\mathbb E \left( \sup_{t\in [0,T]} |Y^+_t|^\varrho + \left[ M \right]_T^{\varrho/2} \right) <+\infty.$$
Moreover if $\xi^+$ is bounded, $Y$ is also bounded. 
\end{proposition}
In the Brownian setting, $M$ is replaced by $\int_0^\cdot Z_s dW_s$ and 
$$ \left[ M \right]_T = \int_0^T |Z_s|^2 ds.$$
\begin{proof}
Recall that the generator $y \mapsto  -\vol_t |y|^p$ is not monotone.
However, if the negative part $Y^-$ of the solution is bounded by some constant 
$c_1 > 0$, that is, if $Y$ is bounded from below by $-c_1$, 
then we can replace the generator $(t,y)\mapsto -\vol_t |y|^p$ by a monotone continuous generator
defined as follows:
\begin{equation} \label{eq:def_mono_gene}
\widetilde f_{-c_1}(s,y) = -\vol_t |y|^p \mathbf 1_{y\geq -c_1} + \vol_t  c_1^{p-1} (py + (p-1)c_1)  \mathbf 1_{y < -c_1}.
\end{equation}
This generator is indeed monotone since
	\[(y-y')(\widetilde f_{-c_1}(t,y)-\widetilde f_{-c_1}(t,y')) \leq p\vol_t c_1^{p-1} (y-y')^2
	\]
for any $t,y,y'.$
Since $Y\geq -c_1$,
$$
Y_t = \xi +\int_t^T f(s,Y_s) ds - (M_T-M_t) = \xi +\int_t^T \widetilde f_{-c_1}(s,Y_s) ds - (M_T-M_t). 
$$
Furthermore, by Assumption  \ref{a:vol}, $\vol$ belongs to $L^1(0,T)$ almost surely. These imply that 
the uniqueness result for BSDE driven by a monotone generator
	(\cite[Proposition 5.24]{pard:rasc:14} or \cite[Lemma 5]{krus:popi:14}) apply in our current setting;
	therefore, if it exists and if its negative part is bounded, the solution of \eqref{e:convex_BSDE} is unique.

We know from Lemma \ref{lem:lower_bound} that the process $z$ is bounded from below by $z_0<0$. Consider the generator $\widetilde f_{z_0}$ and the BSDE 
\begin{equation} \label{eq:modified_convex_BSDE}
Y_t = \xi +\int_t^T \widetilde f_{z_0}(s,Y_s) ds -(M_T-M_t).
\end{equation}
Since $\widetilde f_{z_0}$ is monotone with respect to $y$ and since $\vol \in L^1 (0,T)$, BSDE \eqref{eq:modified_convex_BSDE} has a unique solution $(Y,M)$ (see again \cite[Proposition 5.24]{pard:rasc:14} or \cite[Theorems 1 and 2]{krus:popi:14}). Note that for $\xi = -K$, the solution is $(z,0)$. Since $\xi \geq -K$, the comparison principle  (\cite[Proposition 4]{krus:popi:14} or \cite[Proposition 5.33]{pard:rasc:14}) states that 
	$Y_t \geq z_t \geq z_0$
a.s. for any $t$.
In other words $(Y,M)$ is a solution of the BSDE \eqref{e:convex_BSDE}, and this achieves the proof of the proposition. 
\end{proof}

\begin{remark}[On the negative part of $Y$]\label{rem:boundedness_neg_part}
It\^o-Tanaka formula  (applied to $y \mapsto y^-$ and $Y$) 
implies that $Y^-$ is a subsolution of the BSDE
\begin{equation}\label{eq:BSDE_neg_part}
U_t  = \xi^- + \int_t^T  \mathrm{vol}_s |U_s|^p  \mathbf 1_{U_s\geq 0} ds - \int_t^T dN_s 
\end{equation}
The generator $\widetilde f : (t,y) \mapsto \mathrm{vol}_t y^p \mathbf 1_{y\geq 0}$ is not monotone. 
However, it is increasing and positive. 
From Lemma \ref{lem:lower_bound},  $(U^*,V^*) = (-z,0)$ is a bounded supersolution of this BSDE. 
Following \cite{drap:heyn:kupp:13}, we deduce the existence of a minimal bounded supersolution $(U,N)$ which is also bounded and non-negative (see \cite[Theorems 3.3 and 4.1]{drap:heyn:kupp:13}). 
Using again $\widetilde f_{z_0}$, we deduce 
that $(U,N)$ in fact is the unique solution of the BSDE \eqref{eq:BSDE_neg_part}. 
Therefore, assumptions in the previous proposition 
can be replaced by the existence of a supersolution to the BSDE \eqref{eq:BSDE_neg_part}. 
In other words, these assumptions are sufficient to obtain a supersolution, but not necessary. 
As a by-product of these calculations, we obtain a better bound: 
$0 \leq Y^-_t \leq U_t \leq -z_t$. 
almost surely for any $t \in [0,T].$
\end{remark}

\begin{remark}[On the positive part of $Y$]\label{rem:boundedness_pos_part}
Consider the BSDE 
\begin{align*} \nonumber
\Upsilon_t & =  \xi^+  + \int_t^T \mathrm{vol}_s \left[- (\Upsilon_s-U_s)^p + (U_s)^p\right]  \mathbf 1_{\Upsilon_s \geq 0} ds  -( \mathcal M_T -  \mathcal M_t) \\ 
& = \xi^+  + \int_t^T\hat f(s,\Upsilon_s) ds  -( \mathcal M_T -  \mathcal M_t) ,
\end{align*}
where $U$ is the solution of \eqref{eq:BSDE_neg_part}. 
$U$ is bounded by $K$ and $\partial_y \hat f(s,y) \leq p\mathrm{vol}_s K^{p-1}$
imply that the driver $\hat f$ is monotone:
$$(y-y')(\hat f(t,y)-\hat f(t,y'))\leq p \mathrm{vol}_t  K^{p-1} (y-y')^2$$
and the existence and uniqueness of the solution holds, if $\xi^+$ belongs to some space $L^\varrho(\Omega)$.
Define $Y = \Upsilon - U$ and $M = \mathcal M - N$:
\begin{align*}
Y_t &= \Upsilon_t - U_t = \xi^+ - \xi^- +\int_t^T (\hat f(s,\Upsilon_s) - \widetilde f(s,U_s))  ds \\
& -(\mathcal M_T - N_T - \mathcal M_t + N_t) \\
& = \xi  - \int_t^T \mathrm{vol}_s  (\Upsilon_s - U_s)^p  ds  - (M_T-M_t).
\end{align*}
Hence, $(Y,M)$ solves the BSDE \eqref{e:convex_BSDE} and $0\leq Y^+_t \leq \Upsilon_t$ holds almost surely for $t \in [0,T].$

Another estimate can be obtained using the It\^o-Tanaka formula since
\begin{align*}
(Y_t)^+ & = \xi^+ -\int_t^T \mathrm{vol}_s (Y_s)^p \mathbf 1_{Y_s\geq 0} ds - \int_t^T  \mathbf 1_{Y_s\geq 0} dM_s - \dfrac{1}{2} \int_t^T dL^Y_s. 
\end{align*}
If $\widetilde \Upsilon$ solves the BSDE with monotone generator
$$
\widetilde \Upsilon_t  = \xi^+ -\int_t^T \mathrm{vol}_s (\widetilde \Upsilon_s)^p \mathbf 1_{\widetilde \Upsilon_s\geq 0} ds - \int_t^T d \widetilde{\mathcal M}_s ,
$$
the comparison principle implies that a.s. for any $t$, $Y^+_t \leq \widetilde \Upsilon_t$.
\end{remark}

Our main result on the BSDE \eqref{e:convex_BSDE} when its terminal condition is singular is the following:
\begin{proposition} \label{prop:sovlv_sing_BSDE}
	Suppose  $\xi^- \le K $ and Assumptions \ref{a:vol} and \ref{as:singular} hold. Assume that the filtration $\mathbb F$ is left-continuous at time $T$. 
	Then there exists a minimal\footnote{among all supersolutions with bounded negative part} supersolution $(\Ymin,\Mmin)$ to the BSDE \eqref{e:convex_BSDE} with terminal
	condition $\xi$ such that $\Ymin$ has a left-limit at time $T$ and the negative part of this minimal supersolution is bounded.
\end{proposition}
Several points: there is no condition on $\xi^+$ , the only condition on $\xi^{-}$ is $\xi^- \le K$, in particular,  the terminal
condition \eqref{e:xip} arising from the stochastic optimal control problem satisfies the conditions of this proposition.
Condition \eqref{eq:cond_eta_sing} implies that $\vol$ should remain away from zero close to $T$.

\begin{proof}
The proof proceeds parallel to that of \cite[Proposition 3]{krus:popi:15}. Let us consider for any $n \geq 0$
$$\xi^n = \xi \wedge n.$$
The solution $(Y^n,M^n)$ of 
$$Y^n_t = \xi^n - \int_t^T \vol_s |Y_s^n|^p ds - (M^n_T-M^n_t),$$
has the same upper bound $U=-z$ for the negative part $(Y^n)^-$ for any $n$. 
Using the comparison principle for monotone BSDE, arguing as in \cite{krus:popi:15} gives
$$Y^n \nearrow Y.$$

From Remark \ref{rem:boundedness_pos_part}, $Y^n \leq \Upsilon^n$ with
$$
\Upsilon^n_t  = (\xi^+\wedge n) -\int_t^T \mathrm{vol}_s (\Upsilon^n_s)^p \mathbf 1_{\Upsilon^n_s\geq 0} ds - (N^n_T-N^n_t).
$$
From \cite[Lemma 1]{krus:popi:15}, we have 
\begin{align*}
\Upsilon^n_t & \leq \dfrac{1}{(T-t)^{\widehat p}} \mathbb E \left[ \int_t^T  \left( \dfrac{\widehat p-1}{\vol_s} \right)^{\widehat p-1} ds \bigg| \mathcal F_t\right]
\end{align*}
for any $t \in [T-\epsilon,T]$,
almost surely. 
Therefore, from \eqref{eq:cond_eta_sing}, $Y^n_{t}$, $T-\epsilon \leq t < T$, is finite and bounded in $L^\varsigma(\Omega)$ uniformly with respect to $n$. 
In particular, for any $\epsilon' < \epsilon$,  there exists a constant $C$ such that for any $n$, 
$$\mathbb E \left( |Y^n_{T-\epsilon'}|^\varsigma \right) \leq C.$$
Stability result for BSDE (\cite[Theorem 5.10]{pard:rasc:14}) shows that $(Y^n,M^n)$ converges to $(Y,M)$:
 $$\lim_{n\to +\infty} \mathbb E \left( \sup_{t \in[0,T-\epsilon']} |Y^n_t - Y_t|^\varsigma + \left[ M^n-M\right]^{\varsigma/2}_{T-\epsilon'} \right) = 0.$$
Thus $(Y,M)$ solves \eqref{e:convex_BSDE}: for any $0\leq t\leq r < T$
$$Y_t = Y_r -\int_t^r \vol_s |Y_s|^p ds -(M_r-M_t).$$
Moreover $Y^- \leq U$. And since the filtration is left-continuous at time $T$, we obtain that a.s.
$$\liminf_{t\to T} Y_t \geq \xi.$$
Finally, minimality can be obtained as in the proof of \cite[Proposition 4]{krus:popi:15}. If $(\widetilde Y, \widetilde M)$ is another supersolution, we add to both solutions $\widetilde Y$ and $Y^n$ the quantity $-z$ and the same arguments on $\widetilde Y -z$ and $Y^n -z$ lead to a.s. $\widetilde Y -z \geq Y^n -z$.

The only remaining problem is the existence of a limit at time $T$. 
Compared to \cite[Theorem 2.1]{popi:16}, the novelty is the negative part of $\Ymin$ or of $Y^n$, which approximates $\Ymin$. 
To deal with the negative part,
we can apply the arguments of the proof of \cite[Theorem 2.1]{popi:16} using the function
$$\Theta(y) = \int_y^\infty \dfrac{1}{1+|w|^p} dw$$ 
($\Theta$ can be defined in terms of the hypergeometric functions):
$$\Theta(Y_t) = \mathbb E \left[ \Theta(\xi) | \mathcal F_t \right] + \psi_t^+ - \psi_t^- $$
where $\psi^+$ and $\psi^-$ are two non-negative supermartingales such that $\psi^+$ converges a.s. to zero. 
To obtain this result, we crucially use that the negative part of $Y^n$ is bounded uniformly with respect to $n$ and also that the martingale $N^n = \int \mathbf 1_{Y^n \leq 0} dM^n$ is uniformly bounded in the sense that there exists a constant $C$ such that for any $n$
$$\mathbb E \left[ N^n\right]^2_T \leq C.$$ 
To obtain this last inequality, the It\^o-Tanaka formula for $(Y^n)^-$ is applied (see Remark \ref{rem:boundedness_neg_part} and Equation \eqref{eq:BSDE_neg_part}). Since the negative part of $\xi$ is bounded, a apriori estimate for BSDE (\cite[Proposition 5.7]{pard:rasc:14} or \cite[Proposition 2]{krus:popi:14}) leads to this uniform estimate on $N^n$. The other arguments can be copied from the proof of \cite[Theorem 2.1]{popi:16} with straightforward modifications. This achieves the proof. 
\end{proof}

Our next task is to relate the solution/minimal supersolution of the BSDE  \eqref{e:convex_BSDE} to the solution of the stochastic optimal control problem \eqref{e:oc1m2}.

\subsection{Solution of the quadratic stochastic optimal control problem}\label{ss:verification}
The goal of this subsection is to relate the
minimal supersolution of the BSDE \eqref{e:convex_BSDE} with terminal
condition \eqref{e:termcond} to the solution of the stochastic optimal
control problem \eqref{e:oc1m2}. This will be achieved through a verification
argument based on the convexity of the cost structure of
\eqref{e:oc1m2}. The convex structure holds only for $\widehat{p}=2$,
for this reason in this subsection we assume $\widehat{p} = p =2.$
We deal with the case $\widehat{p} \neq 2$ in the next section in a Markovian
framework.

For ease of reference let us restate our stochastic optimal control problem \eqref{e:oc1m2}:
\begin{equation}\label{e:soc1}
	\inf_{Q \in \contrset } {\mathbb E}\left[ \int_0^T  \frac{|Q'_t|^2}{\vol_t}  dt + \xi |Q_T|^2\right] = 
	\inf_{Q \in {\mathcal A} } {\mathbb E}\left[ \int_0^T  \frac{|Q'_t|^2}{\vol_t}  dt + \xi |Q_T|^2\right]; 
\end{equation}
we remind the reader that in the second formulation we are using the convention
\eqref{e:infconv}; the definitions of $\contrset$, ${\mathcal A}$ and $\vol$ 
are given Definition \eqref{d:efAIS}, \eqref{e:conventionA} and \eqref{e:vol}.

\begin{remark}\label{rem:int_cond_control}
Under Assumption \ref{a:vol}  on $\mathrm{vol}$ (with $p=2$) 
	the Cauchy-Schwarz inequality implies
\begin{align*}
\sup_{t\in [0,T]} |Q_t-Q_0|^2 & = \sup_{t\in [0,T]} \left| \int_0^t Q'_s ds \right|^2   \leq  \left( \int_0^T \mathrm{vol}_s ds \right) \left( \int_0^T\dfrac{|Q'_s|^2}{ \mathrm{vol}_s} ds \right)\\
& \leq  \dfrac{1}{K}  \left( \int_0^T\dfrac{|Q'_s|^2}{ \mathrm{vol}_s} ds \right). 
\end{align*}
Hence if $Q \in \contrset$, $\mathbb E \sup_{t\in [0,T]} |Q_t-Q_0|^{2\varrho} < +\infty$.  
In particular 
$$|Q_T|^2 \leq 2 |Q_0|^2 + 2 \dfrac{1}{K} \left( \int_0^T\dfrac{|Q'_s|^2}{ \mathrm{vol}_s} ds \right).$$
Hence for any bounded $\xi$, the integrability condition in Definition \ref{d:efAIS} is equivalent to
$$\mathbb E \left[ \left( \int_0^T\dfrac{|Q'_s|^2}{ \mathrm{vol}_s} ds \right)^\varrho \right]  < +\infty.$$
\end{remark}

Our goal is to prove the following result:
\begin{proposition} \label{prop:constrained_control_pb}
	Suppose $p = \widehat{p} = 2.$
	Suppose Assumptions \ref{a:vol} and \ref{as:singular} hold, suppose $\xi^- \le K$
	and let $(\Ymin,\Mmin)$ be the minimal supersolution of \eqref{e:convex_BSDE}, \eqref{e:termcond}.
	Then
\begin{equation}\label{e:optimalcontrol1}
Q^*_t = Q_0 \exp \left( - \int_0^t \Ymin_s \mathrm{vol}_s ds \right), t \in [0,T),
\end{equation}
	(equivalently, $(Q^*)'_t = -  \Ymin_t \mathrm{vol}_t Q^*_t$)
	is the optimal state process for the stochastic optimal control problem \eqref{e:soc1}.
	Moreover the value function of \eqref{e:soc1} at time $t$, namely
	\begin{equation*}
	\inf_{Q \in \contrset(t)} {\mathbb E}\left[ \int_t^T  \frac{|Q'_s|^2}{\mathrm{vol}_s}  ds + \xi |Q_T|^2 \bigg| \mathcal F_t\right], 
\end{equation*}
	 is given by $V(t,q) = q^2 \Ymin_t$ ($\contrset(t)$ is defined by Definition \ref{d:efAIS}, but the process $Q$ starts at time $t$ from 
	 the deterministic position $q$). 
\end{proposition}
The proof directly follows from  the next two lemmas and is given at the end of this subsection.
Let's call $J$ the expression inside the $\min$ in \eqref{e:soc1}: for $v = Q'$
\[
J(v) =  \int_0^T \frac{v_t^2 }{\vol_t} dt +\xi \left(Q_0+ \int_0^T v_t dt\right)^2.
\]
We start with the following observation:
\begin{lemma} \label{lem:convex_cost}
If 
\begin{equation} \label{eq:convex_cond}
\xi^- \int_0^T \mathrm{vol}_t dt < 1
 \end{equation}
 almost surely,
then the functional $v\mapsto J(v)$ is strictly convex. The G\^ateaux derivative of $J$ at point $v$ in direction $w$, is given by
$$\langle DJ(v), w\rangle =2 \int_0^T  \frac{v_tw_t }{ \mathrm{vol}_t } dt +2\xi \left(Q_0+ \int_0^T v_t dt\right)\left( \int_0^T w_t dt\right).$$
\end{lemma}
\begin{remark}
	Assumption \ref{a:vol} with $p=\widehat p = 2$ and $\xi^- \le K$  imply \eqref{eq:convex_cond}. 
\end{remark}
\begin{proof}
Taking $v$ and $\widetilde v$ and $\theta \in [0,1]$, we have
\begin{align*}
& J(\theta v + (1-\theta) \widetilde v) - \theta J(v) - (1-\theta)J(\widetilde v) \\
& = -\theta(1-\theta) \left[ \int_0^T  \frac{(v_t-\widetilde v_t)^2 }{\vol_t} dt + \xi \left( \int_0^T (v_t-\widetilde v_t) dt\right)^2 \right] \\
& \leq -\theta(1-\theta) \left[ \int_0^T  \frac{(v_t-\widetilde v_t)^2 }{\vol_t} dt - \xi^- \left( \int_0^T (v_t-\widetilde v_t) dt\right)^2 \right] \\
& \leq \theta(1-\theta) \left[-1 + \xi^-  \int_0^T \mathrm{vol}_t dt \right]  \int_0^T  \frac{(v_t-\widetilde v_t)^2 }{\vol_t} dt \leq 0.
\end{align*}
We use the Cauchy-Schwarz inequality for the inequality. Now for any $\epsilon > 0$ and $v$ and $w$ 
\begin{align*}
&\dfrac{1}{\epsilon} (J(v + \epsilon w) -J(v)) \\
& = 2 \int_0^T  \frac{v_tw_t }{\vol_t} dt  + \epsilon  \int_0^T  \frac{w_t^2 }{\vol_t} dt \\
& \quad +2\xi \left(Q_0+ \int_0^T v_t dt\right)\left( \int_0^T w_t dt\right) - \epsilon  \xi \left( \int_0^T w_t dt\right)^2 .
\end{align*}
Letting $\epsilon$ to zero gives the desired formula. 
\end{proof}

The last intermediate result we need is a version of Proposition \ref{prop:constrained_control_pb} where
$\xi$ is bounded. 
\begin{lemma} \label{lem:unconstrainted_control_pb}
	Suppose $\xi$ is bounded. 
	If $(Y,M)$ is the solution of  \eqref{e:convex_BSDE} with terminal condition $Y_T = \xi$, 
	then the optimal state process $Q^*$ (resp. optimal control $v^*=(Q^*)'$) of \eqref{e:soc1} is given by 
$$Q^*_t = Q_0 - \int_0^t (Y_s \mathrm{vol}_s) Q^*_s ds \quad (\mbox{resp. } -Y_s \mathrm{vol}_s Q^*_s ).$$
Moreover, $Y_0 (Q_0)^2$ is the value function of the control problem. 
\end{lemma}
\begin{proof}
	Note that for bounded $\xi$, the BSDE \eqref{e:convex_BSDE} with terminal condition $Y_T = \xi$
	has a unique solution and it equals the minimal supersolution: $Y = \Ymin$, $M=\Mmin$. 
	Since $\xi$ is bounded, from Proposition \ref{thm:solve_aux_BSDE}, $Y$ is bounded and for any $\varrho>1$
$$\mathbb E \left( \left[ M \right]_T^{\frac{\varrho}{\varrho-1}} \right) <+\infty.$$
Thus
$$Q^*_t = Q_0 \exp \left( -\int_0^t Y_s \vol_s ds \right)$$
is also bounded and from Assumption \ref{a:vol}
$$\int_0^T \dfrac{(v^*_s)^2}{\vol_s} ds = \int_0^T  \mathrm{vol}_s (Y_sQ^*_s)^2  ds$$
is bounded. Thus $Q^*$ is in $\contrset$. 

Let $Q$ be in $\contrset$ with $L^\varrho$-integrability, and define $v=Q'$, $w=v^*-v$ and  $\widehat Q = \int_0 w_s ds$. If 
\[
	N_t \doteq 2Y_t Q^*_t= -2 \dfrac{v^*_t}{\vol_t},
\]
integration by parts gives:
\[
N_t = 2Y_0 Q_0 + 2 \int_0^t Q^*_s \vol_s(Y_s)^2 ds - \int_0^t 2 Y_s Y_s \vol_s Q^*_s ds +2\int_0^t Q^*_s dM_s =2\int_0^t Q^*_s dM_s.
\]
Hence $N$ and $\int (Q^*)^2 dM$ are martingales (since $Q^*$ is bounded). Moreover 
$$\int_0^T ( \widehat Q_s Q^*_s )^2d[M]_s  \leq C \sup_{t\in[0,T]} ( \widehat Q_s)^2 [M]_T;$$
this and \ref{rem:int_cond_control} imply
$$\mathbb E \int_0^T ( \widehat Q_s Q^*_s) d[M]_s  \leq C \left[\mathbb E \sup_{t\in[0,T]} ( \widehat Q_s)^{2\varrho} \right]^{\frac{1}{\varrho}} \left[ \mathbb E \left[ M \right]_T^{\frac{\varrho}{\varrho-1}} \right]^{\frac{\varrho-1}{\varrho}} <+\infty.$$
Hence $\int \widehat Q Q^* dM$ is also a martingale.

Now integration by parts implies:
\[
\int_0^T w_t \frac{2 v^*_t }{\vol_t} dt = \int_0^T \widehat Q_t  dN_t - N_T \widehat Q_T = \int_0^T \widehat Q_t  dN_t  - 2 \xi q^*_T  \widehat Q_T.
\]
This and Lemma \ref{lem:convex_cost} give
\[
\langle DJ(v^*), w\rangle = \int_0^T \widehat Q_t  dN_t= 2  \int_0^T \widehat Q_t  Q^*_t dM_t,
\]
which is a martingale. With the convexity of $J$ we obtain
\[
\mathbb E (J(v^*) - J(v) ) \leq \mathbb E \langle DJ(v^*), v^*-v\rangle = 0.
\]
Therefore, $v^*$ is the optimal control (unique from the strict convexity of $J$).  
It\^o's formula applied to $Y_t (Q^*_t)^2$ gives
\begin{align*}
d(Y_t (Q^*_t)^2) &= (Q^*_t)^2 \vol_t (Y_t)^2 dt + (Q^*_t)^2 dM_t + 2Y_t q^*_t (-Y_t \vol_t q^*_t) dt \\
& = - \vol_t (Q^*_t Y_t)^2 dt + (Q^*_t)^2 dM_t = - \dfrac{(v_t^*)^2}{\vol_t} dt + (Q^*_t)^2 dM_t .
\end{align*}
Since $\int (Q^*)^2 dM$ is a martingale, 
\begin{align*}
Y_t (Q^*_t)^2 &=\mathbb E \left[  \int_t^T \dfrac{(v_s^*)^2}{\vol_s} ds + \xi (Q^*_T)^2 \bigg| \mathcal F_t \right].
\end{align*}
In other words $Y_t (Q^*_t)^2$ is the value function of the control problem.
\end{proof}

We now give
\begin{proof}[Proof of Proposition  \ref{prop:constrained_control_pb}]
For $N= 2\Ymin Q^*$,
one can use arguments parallel to those in the proof of the previous lemma to show that $N$ is a martingale on $[0,T)$. 
Since $(\Ymin)^-$ is bounded (by $K$) we have
$$ \exp \left( - \int_0^t \Ymin_s \vol_s ds \right) \leq \exp \left(  \int_0^t (\Ymin_s)^- \vol_s ds \right) \leq  \exp \left( K \int_0^t  \vol_s ds \right).$$
	This and Assumption \ref{a:vol} imply that $Q^*$ is also bounded. Since $(\Ymin)^-$ is also bounded, the martingale $N$ is bounded from below. 
Therefore, the limit at time $T$ of $N$ exists in $\mathbb R$ and
$$Q^*_t = \dfrac{N_t}{2\Ymin_t}$$
tends to zero a.s. on the set $ \bm S = \{\xi = +\infty\}$, since $\lim_{t\to T} \Ymin_t {\bm 1}_{ \bm  S} = +\infty$.

Now we apply It\^o's formula to $\Ymin (Q^*)^2$: for any $0\leq t \leq r < T$
\begin{align*}
\Ymin_t (Q^*_t)^2 & = \Ymin_r (Q^*_r)^2 - \int_t^r (Q^*_s)^2\vol_s (\Ymin_s)^2 ds  \\
& - \int_t^r \Ymin_s 2 (Q^*_s) (-\Ymin_s \vol_s Q^*_s) ds -\int_t^r (Q^*_s)^2 d\Mmin_s  \\
& = \Ymin_r (Q^*_r)^2 + \int_t^r \dfrac{(v^*_s)^2}{\vol_s} ds-\int_t^r (Q^*_s)^2 d\Mmin_s 
\end{align*}
with $v^*_s= -\vol_s\Ymin_s  Q^*_s $. Taking the conditional expectation we get 
\begin{align*}
\Ymin_t (Q^*_t)^2 & = \mathbb E \left[ \Ymin_r (Q^*_r)^2 + \int_t^r \dfrac{(v^*_s)^2}{\vol_s} ds \bigg| \mathcal F_t\right].
\end{align*}
By the the monotone convergence theorem
$$\liminf_{r\to T} \mathbb E \left[ \int_t^r \dfrac{(v^*_s)^2}{\vol_s} ds \bigg| \mathcal F_t\right] =   \mathbb E \left[ \int_t^T \dfrac{(v^*_s)^2}{\vol_s} ds \bigg| \mathcal F_t\right].$$
And by Fatou's lemma ($(\Ymin)^-$ is bounded) 
$$\liminf_{r\to T} \mathbb E \left[\Ymin_r (Q^*_r)^2  \bigg| \mathcal F_t\right] \geq  \mathbb E \left[\liminf_{r\to T} (\Ymin_r (Q^*_r)^2)  \bigg| \mathcal F_t\right] .$$
Recall the definition of $N=2\Ymin Q^*$ and that the limit of $N$ at time $T$ exists in $\mathbb R$. Moreover $\lim_{r\to T} Q^*_r = 0 = Q^*_T$, when $\lim_{r\to T} \Ymin_r = +\infty$. Therefore, if $\lim_{r\to T} \Ymin_r = +\infty$, then 
$$\liminf_{r\to T} (\Ymin_r (Q^*_r)^2) = \liminf_{r\to T} (N_r)\lim_{r\to T} Q^*_r  = 0 = \xi (Q^*_T)^2$$ 
(with the convention $\infty \cdot 0 = 0$). If  $\liminf_{r\to T} \Ymin_r < +\infty$, then 
$$\liminf_{r\to T} (\Ymin_r (Q^*_r)^2) \geq \xi (Q^*_T)^2 .$$ 
In both cases we obtain
\begin{align*}
\Ymin_t (Q^*_t)^2 & \geq \mathbb E \left[  \int_t^T \dfrac{(v^*_s)^2}{\vol_s} ds + \xi (Q^*_T)^2 \bigg| \mathcal F_t\right].
\end{align*}
Thus $\Ymin q^2$ dominates the value function $V(\cdot,q)$ of the constrained control problem.

Now if $Q$ is in $\contrset$, it is in $\mathcal A_{I,\emptyset}$ (no terminal constraint on $Q_T$). Therefore, the value function $V$ dominates the value function of the unconstrained control problem with terminal penalty $\xi \wedge n$, for any $n$. Denote by $Y^n$ the solution of the BSDE \eqref{e:convex_BSDE} 
	with bounded terminal value $Y^n_T = \xi \wedge n$. From Lemma \eqref{lem:unconstrainted_control_pb}, $Y^n q^2$ is the value function of the unconstrained control problem. We deduce that for any $n$
$$Y^n_t q^2 \leq V(t,q) \leq \Ymin_t q^2  .$$
Since $Y^n$ converges to $\Ymin$, we obtain that $\Ymin q^2$ is the value function of the constrained control problem and that $Q^*$ is the optimal state process. 

Note that the value function is finite at time $0$, that is
$$\mathbb E \left[  \int_0^T \dfrac{(v^*_s)^2}{\vol_s} ds + \xi (Q^*_T)^2 \right] \leq \Ymin_0 (Q^*_0)^2 < +\infty.$$
Using \eqref{eq:cond_eta_sing} and the proof of Proposition \ref{prop:sovlv_sing_BSDE}, for $t < T$, $\Ymin_t$ belongs to $L^\varsigma(\Omega)$ and we can also deduce that 
$$\mathbb E \left[  \left( \int_0^T \dfrac{(v^*_s)^2}{\vol_s} ds + \xi (Q^*_T)^2 \right)^\varsigma \ \right] < +\infty.$$
Therefore $Q^*$ belongs to $\contrset$ and this achieves the proof. 
\end{proof}

\subsection{Reduction to random time interval $[\![0, \tau_\ell \wedge T]\!]$ for $I_t = 1_{\{t \leq \tau_\ell \}}$}\label{ss:red}

In general, $\mathbb P(T-\epsilon < \tau_\ell \leq T) >0$ for any $\epsilon > 0$. 
therefore,
for $I = I^{(2)} = \mathbf 1_{\{t \le \tau_\ell \}}$, 
Assumption \ref{as:singular} in general doesn't hold regardless of what $\Vol$ is. For $I^{(2)}$ a natural way to deal with
this is to consider the problem
on the random interval $[\![0,\tau_\ell  \wedge T]\!]$:
for $Q \in {\mathcal A}_{I^{(2)},{\bm S}}$  
we have
$Q_t = Q_{\tau_\ell}$, $Q'_t = 0$ for $t > \tau_\ell.$ 
Therefore, the stochastic
optimal control problem \eqref{e:oc1m2} can also be expressed as
\begin{equation}\label{e:oc1m2red}
	\min_{Q \in {\mathcal A}} {\mathbb E}\left[ \int_0^{\tau_\ell \wedge T} 
	\frac{(Q'_t)^{\widehat{p}}}{{\Vol}_t}  dt +\left(-K {\bm 1}_{\{\tau_\ell < T \}}+ \infty \cdot {\bm 1}_{\{\tau_\ell \geq T\}} \right)Q_T^{\widehat{p}}\right].
\end{equation}
The corresponding BSDE is again \eqref{e:convex_BSDE} but
with terminal condition
\begin{equation}\label{e:xi2}
	Y_{\tau_l \wedge T} =\xi^{(2)} =  -K {\bm 1}_{\{\tau_\ell < T \}}+ \infty \cdot {\bm 1}_{\{\tau_\ell \geq T\}},
\end{equation}
and we have a BSDE with random terminal time $\tau_\ell \wedge T$. 
Note that now the dynamics of $Y$ is considered only on the random interval $[\![0,\tau_\ell  \wedge T]\!]$ and
$I_t^{(2)} = 1$ on this random interval. 

The next proposition (formulated in terms of a general stopping time $\tau$)
states that the BSDE \eqref{e:convex_BSDE} has a minimal supersolution for terminal conditions
of the form
\begin{equation}\label{e:xi2general}
	Y_{\tau \wedge T} = \xi = \zeta \mathbf 1_{\{\tau \geq T\}} + \psi_\tau \mathbf 1_{\{\tau < T\}},
\end{equation}
where $\tau$ is a stopping time and $\zeta$, $\psi$ are bounded from below by $-K$ and $\psi$ is also bounded from above;
\eqref{e:xi2} is a special case of \eqref{e:xi2general}.

\begin{proposition}  \label{prop:BSDE_stop_time} 
Assumption \ref{a:vol} holds. Let $\xi$ be as in \eqref{e:xi2general}.
If for some $\varrho > 1$
	\begin{equation}\label{e:intzeta}
\mathbb E \int_0^T \mathrm{vol}_s (\mathbb E [\zeta | \mathcal F_s])^\varrho ds <+\infty,
	\end{equation}
	then the BSDE \eqref{e:convex_BSDE} with terminal condition \eqref{e:xi2general} has a unique solution $(Y,M)$ such that 
	\begin{itemize}
\item  $\mathbb P$-a.s., on the set $\{t\geq \tau \wedge T \}$, $Y_t =\xi$ and $M_t=0$,
\item $\mathbb P$-a.s., for all $0\leq t \leq r \leq T$,
$$Y_{t\wedge \tau}  =  Y_{r\wedge \tau} + \int_{t\wedge \tau}^{r\wedge \tau} \mathrm{vol}_s |Y_s|^p ds - \int_{t\wedge \tau}^{r\wedge \tau} dM_s,$$
\item and for some constant $C$ depending on $\mathrm{vol}$, $K$ and $\varrho$  
$$
\mathbb E \left[  |Y_{ \tau\wedge T}|^\varrho + \int_{0}^{ \tau\wedge T}  |Y_{s}|^{\varrho}  ds + [M]_{ \tau \wedge T}^{\varrho/2}\right]  \leq C \mathbb E  |\xi|^\varrho.$$
\end{itemize}
	If, instead of the integrability condition \eqref{e:intzeta} on $\zeta$, 
	the following modified version of condition \eqref{eq:cond_eta_sing} holds
\begin{equation} \label{eq:cond_eta_sing_2}
	\mathbb E \left[ \int_{\tau\wedge (T-\epsilon)}^{\tau \wedge T} \dfrac{1}{(\mathrm{vol}_s)^{(\widehat p-1)\varsigma}} ds \right] <+\infty ,
\end{equation}
then
	there exists a minimal supersolution $(\Ymin,\Zmin)$ for the BSDE \eqref{e:convex_BSDE} with
	terminal condition \eqref{e:xi2general}.
\end{proposition} 
\begin{proof}
The first part of the claim comes from \cite[Proposition 6]{krus:popi:14} on BSDEs with random terminal time. Note that our terminal time $\tau \wedge T$ is bounded. As in the proof of Proposition \ref{thm:solve_aux_BSDE}, since the negative part of the data is bounded, we can modify the generator in order to have a monotone generator (see Equation \eqref{eq:def_mono_gene}). 

In the singular case, we proceed by truncation. For $n$ sufficiently large:
$$\xi \wedge n = ( \zeta\wedge n) \mathbf 1_{\tau \geq T} + \psi_\tau \mathbf 1_{\tau < T}.$$
Again we construct a sequence of solutions $(Y^n,M^n)$, such that $(Y^n)$ is non-decreasing. It only remains to control $Y^n$, uniformly in $n$. 

 We define $(\widehat U,\widehat V)$ as the solution with terminal condition $\psi_\tau \mathbf 1_{\tau < T} - \zeta^- \mathbf 1_{\tau \geq T}$. Note that $\widehat U$ is bounded.
Now for any $t\leq s \leq T$ we have
\begin{align*}
Y^n_{t\wedge \tau} - \widehat U_{t\wedge \tau} & = Y^n_{s\wedge \tau} - \widehat U_{s\wedge \tau} - \int_{t\wedge \tau} ^{s\wedge \tau} \vol_r \left[ \left|Y^n_r-\widehat U_r+\widehat U_r \right|^p - \left| \widehat U_r \right|^p \right]dr \\
& -  \int_{t\wedge \tau} ^{s\wedge \tau} (dM^n_r - d\widehat V_r).
\end{align*}
Hence $\Upsilon^n = Y^n - \widehat U$ solves the BSDE with terminal condition $( \zeta^+ \wedge L) \mathbf 1_{\tau \geq T} $ and generator 
$$g(s,y) = - \vol_s \left[ \left|y+\widehat U_s\right|^p - |\widehat U_s|^p \right].$$
Once again, the map $y\mapsto g(s,y)$ is not monotone on $\mathbb R$. But since the negative part of the data is bounded, we can modify $g$ on $(-\infty,-c_1)$, such that $g$ becomes monotone (see Equation \eqref{eq:def_mono_gene}). And again by the comparison principle, $\Upsilon^n$ is non-negative. 

To control $\Upsilon^n$, we again follow the arguments of \cite[Lemma 1]{krus:popi:15}. Young's inequality implies that for any $y\geq 0$ and $c\geq 0$
$$ \left|y+\widehat U_s \right|^p\geq \dfrac{\widehat p}{\widehat p-1} c^{ p -1}(y+\widehat U_s) - \dfrac{c^{p}}{\widehat p-1}.$$
Recall that $\widehat p$ is the H\"older conjugate of $p$. Thus with $c^{p-1}=(\widehat p-1)/(\vol_s(T-s))$
$$g(s,y) \leq - \dfrac{\widehat p}{T-s}y + \left(\dfrac{\widehat p-1}{\vol_s} \right)^{\widehat p-1}\left( \dfrac{1}{T-s}\right)^{\widehat p} - \dfrac{\widehat p}{T-s}\widehat U_s+ \vol_s |\widehat U_s|^p.$$ 
Since $\Upsilon^n_{T\wedge \tau} = 0$ if $\tau < T$, explicit solution for linear BSDE and comparison principle imply that for $t < T$
\begin{align*}
\Upsilon^n_{\tau\wedge t} & \leq \dfrac{1}{(T-\tau\wedge t)^{\widehat p}} \mathbb E \left[ \int_{t\wedge \tau}^{T\wedge \tau}  \left(\dfrac{\widehat p-1}{\vol_s} \right)^{\widehat p-1}ds \bigg| \mathcal F_{\tau\wedge t}\right] \\
& +\dfrac{1}{(T-\tau\wedge t)^{\widehat p}} \mathbb E \left[ \int_{t\wedge \tau}^{T\wedge \tau}  (T-s)^{\widehat p-1} [  \vol_s(T-s)|\widehat U_s|^p- \widehat p \widehat U_s] ds \bigg| \mathcal F_{\tau\wedge t}\right]. 
\end{align*}
This uniform bound on $\Upsilon^n$, thus on $Y^n$, allows us to define the solution of the BSDE with terminal time $\tau \wedge T$ and a singular terminal condition. 
\end{proof}

For $\tau = \tau_\ell$ we have $\vol = \Vol$ and for the existence of a minimal supersolution $\Vol$ must satisfy \eqref{eq:cond_eta_sing_2}.
The next proposition connects the value function of \eqref{e:oc1m2red} to the minimal supersolution $(\Ymin,\Mmin)$ whose existence was derived above.
\begin{proposition} \label{prop:control_pb_stop_time}
	Suppose ${\rm Vol}$ satisfies \eqref{eq:cond_eta_sing_2} with $\widehat p = 2$. 
	Let $(\Ymin,\Mmin)$ be the minimal supersolution of \eqref{e:convex_BSDE} and \eqref{e:xi2}. For any $t \in [0,T)$, 
	\[
Q^*_t = Q_0 \exp \left( - \int_0^t \Ymin_s \mathrm{vol}_s ds \right),
	\]
	(equivalently, $(Q^*)'_t = -  \Ymin_t \mathrm{vol}_t Q^*_t$)
is the optimal control for the stochastic optimal control problem \eqref{e:oc1m2red}.
Moreover the value function of the same control problem at time $t$ equals $\Ymin_t q^2.$
\end{proposition}
\begin{proof}
Following the arguments of the proof of Proposition \ref{prop:constrained_control_pb}, we have that for any $0 \leq t \leq r < T$:
\begin{align*}
\Ymin_{\tau_\ell \wedge t} (Q^*_{\tau_\ell \wedge t})^2 & = \Ymin_{\tau_\ell \wedge r} (Q^*_{\tau_\ell \wedge r})^2 + \int_{\tau_\ell \wedge t}^{\tau_\ell \wedge r} \dfrac{(v^*_s)^2}{\vol_s} ds-\int_{\tau_\ell \wedge t}^{\tau_\ell \wedge r} (Q^*_s)^2 d\Mmin_s 
\end{align*}
with $v^* = (Q^*)'$. 
Taking the conditional expectation and passing to the limit on $r$, we obtain
\begin{align*}
\Ymin_{\tau_\ell \wedge t} (Q^*_{\tau_\ell \wedge t})^2 & \geq \mathbb E \left[  \int_{\tau_\ell \wedge t}^{\tau_\ell \wedge T} \dfrac{(v^*_s)^2}{\vol_s} ds + \xi (Q^*_{\tau_\ell \wedge T})^2 \bigg| \mathcal F_t\right].
\end{align*}
	The rest of the proof continues as that of Proposition 
	\ref{prop:constrained_control_pb}.
\end{proof}

\subsection{Reduction to time interval $[0, T-\delta]$ for $(I,\vect{S}) = (I^{(3)},\vect{S}^{(3)})$ and $(I^{(4)},\vect{S}^{(4)})$} \label{ssect:cases_3_4}

As already noted, for $I = I^{(3)}$ and $I=I^{(4)}$, Assumption \ref{as:singular} doesn't hold in general since both $I^{(3)}_t$ and $I^{(4)}_t$ can be zero for $t$ 
arbitrarily close to $T$. 
However, in both of these cases the problem can
be reduced to the time interval $[0,T-\delta]$ 
where this assumption is no longer needed.

Both $(I,{\bm S}) = (I^{(3)},{\bm S}^{(3)})$ and $(I,{\bm S}) = (I^{(4)},{\bm S}^{(4)})$ consist of two phases: before and after time $T-\delta$, the
reason for this was explained in the paragraph following \eqref{e:stopbelowL}. 
In both cases 
the trading process for $I^{(3)}$ and $I^{(4)}$ proceeds exactly as in $I^{(2)}$ after time $T-\delta$: if the algorithm
is in trading mode at time $T-\delta$, the position is fully closed only when the price remains above $\ell$ throughout the
interval $[T-\delta, T]$; trading stops (and doesn't restart) if the price hits $\ell$. This implies that the stochastic optimal control problem
\eqref{e:soc1} can be written as
\begin{equation}\label{e:soc34}
	\inf_{Q \in {\mathcal A}} {\mathbb E}\left[ \int_0^{T-\delta} \frac{|Q'_t|^{\widehat{p}}}{\vol_t}  dt + 
	I^{(j)}_{T-\delta} 
	V^\infty_{T-\delta}
	 -(1-I^{(j)}_{T-\delta})K |Q_{T-\delta}|^{\widehat{p}}
	\right]
\end{equation}
where $V^\infty_{T-\delta}$ is the value function of the stochastic optimal control problem corresponding
to $(I^{2},{\bm S}^{(2)})$
on the time interval $[\![T-\delta,T \wedge \tau_{T-\delta,\ell}]\!]$ with initial position $Q_{T-\delta}.$
For $\hat{p} =2$, we know by Proposition \ref{prop:control_pb_stop_time} that 
 $V^\infty_{T-\delta}= Q_{T-\delta}^{\widehat{p}} Y^{\infty,T-\delta}_{T-\delta}$ 
 where
$Y^{\infty,T-\delta}$
is
the minimal supersolution of the BSDE \eqref{e:convex_BSDE} on the time interval $[\![T-\delta,\tau_{T-\delta,\ell} \wedge T]\!]$ with
terminal condition
\[\zeta= -K {\bm 1}_{\{ \tau_{T-\delta,\ell} < T \}} + \infty \cdot {\bm 1}_{\{ \tau_{T-\delta,\ell} \geq T \}}.\]
Existence of $Y^{\infty,T-\delta}$ follows from Proposition \ref{prop:BSDE_stop_time}.
Then for $\widehat{p}=2$, \eqref{e:soc34} can be written as
\begin{equation*}
	\inf_{Q \in {\mathcal A}} {\mathbb E}\left[ \int_0^{T-\delta} \frac{|Q'_t|^{\widehat{p}}}{\vol_t}  dt + \xi |Q_{T-\delta}|^{\widehat{p}}\right]
\end{equation*}
with terminal cost factor
\begin{equation}\label{e:termcdred34}
	\xi = 
I^{(j)}_{T-\delta} 
	Y^{\infty,T-\delta}_{T-\delta}- (1-I^{(j)}_{T-\delta})K.
\end{equation}

$Y^{\infty,T-\delta}_{T-\delta}$ is $ \mathcal F_{T-\delta}$-measurable and belongs to $L^\varsigma( \Omega)$. 
Now consider the solution $(Y,M)$ of the BSDE \eqref{e:convex_BSDE} on the time interval $[0,T-\delta]$ with the $\mathcal F_{T-\delta}$-measurable terminal condition $\xi$ of \eqref{e:termcdred34}.
Define 
$$\Ymin_t = \begin{cases} Y_t & \mbox{ if } t < T-\delta \\
Y^{\infty,T-\delta}_t   I_{T-\delta} &  \mbox{ if } t \geq T-\delta .
\end{cases}$$
\begin{lemma} \label{lem:reduc_time_interval}
	Suppose $\widehat{p}=p=2$.
For any $t \in [0,T)$, 
	\[
Q^*_t = Q_0 \exp \left( - \int_0^t \Ymin_s \mathrm{vol}_s ds \right),
	\]
	(equivalently, $v^*_t = -  \Ymin_t \mathrm{vol}_t Q^*_t$)
	is the optimal control for the stochastic optimal control problem corresponding to $(I,S) = (I^{(3)},S^{(3)})$ or $(I,S) = (I^{(4)},S^{(4)})$.
Moreover, the value function of the same control problem at time $t$ equals $q^2 \Ymin_t.$
\end{lemma}
\begin{proof}
From our previous arguments, we know that if $v^*_t = - Y_t \vol_t Q^*_t$ is the derivative of $Q^*_t$, with starting point $Q_0$, then for any other strategy $v$ we have
\begin{align}\label{eq:optim_before_T_delta}
Y_0(Q_0)^2 &=\mathbb E \left[  \int_0^{T-\delta} \dfrac{(v^*_s)^2}{\vol_s} ds + \xi(Q^*_{T-\delta})^2 \right]  \leq  \mathbb E \left[  \int_0^{T-\delta} \dfrac{(v_s)^2}{\vol_s} ds + \xi(Q^v_{T-\delta})^2 \right] .
\end{align}
Now let us define for $t\geq T-\delta$
$$Q^{\infty,T-\delta}_t = Q^*_{T-\delta} - \int_{T-\delta}^t Y^{\infty,T-\delta}_s \vol_s Q^{\infty,T-\delta}_s ds$$
and $v^{\infty,T-\delta} = - Y^{\infty,T-\delta} \vol Q^{\infty,T-\delta}$. From Proposition \ref{prop:control_pb_stop_time}, $Y^{\infty,T-\delta} (Q^{\infty,T-\delta})^2$  is the value function of the related control problem starting at time $T-\delta$ from $Q^*_{T-\delta}$:
\begin{align}\nonumber
Y^{\infty,T-\delta}_{T-\delta} (Q^*_{T-\delta})^2 &=\mathbb E \left[  \int_{T-\delta}^{\tau_{T-\delta,\ell}\wedge T} \dfrac{(v^{\infty,T-\delta}_s)^2}{\vol_s} ds + \zeta (Q^{\infty,T-\delta}_{\tau_{T-\delta,\ell}\wedge T})^2 \bigg| \mathcal F_{T-\delta}\right] \\ \label{eq:optim_after_T_delta}
& \leq \mathbb E \left[  \int_{T-\delta}^{\tau_{T-\delta,\ell}\wedge T} \dfrac{(v^{T-\delta}_s)^2}{\vol_s} ds + \zeta (Q^{v,T-\delta}_{\tau_{T-\delta,\ell}\wedge T})^2 \bigg| \mathcal F_{T-\delta}\right] 
\end{align}
for any process $Q^{v,T-\delta}$ with derivative $v^{T-\delta}$ starting at time $T-\delta$ from $Q^*_{T-\delta}$. Multiplying this equality by $I_{T-\delta}$, we obtain 
\begin{align*}
& \xi (Q^*_{T-\delta})^2 = I_{T-\delta} Y^{\infty,T-\delta}_{T-\delta} (Q^*_{T-\delta})^2  - (1-I_{T-\delta})k/2(Q^*_{T-\delta})^2  \\
& = \mathbb E \left[  \int_{T-\delta}^{\tau_{T-\delta,\ell}\wedge T} \dfrac{(I_{T-\delta}v^{\infty,T-\delta}_s)^2}{\vol_s} ds +I_{T-\delta} \zeta (Q^{\infty,T-\delta}_{\tau_{T-\delta,\ell}\wedge T})^2 \bigg| \mathcal F_{T-\delta}\right] - K (1-I_{T-\delta})(Q^*_{T-\delta})^2.
\end{align*}
If for $t\geq T-\delta$
$$Q^{*}_t = Q^*_{T-\delta} - \int_{T-\delta}^t I_{T-\delta} Y^{\infty,T-\delta}_s \vol_s Q^{*}_s ds,$$
we have 
$$Q^{*}_t = \begin{cases}
Q^{\infty,T-\delta}_t & \mbox{ if } I_{T-\delta} = 1\\
Q^*_{T-\delta} & \mbox{ if } I_{T-\delta} = 0
\end{cases}$$
and therefore 
\begin{align*}
& \xi (Q^*_{T-\delta})^2= \mathbb E \left[  \int_{T-\delta}^{\tau_{T-\delta,\ell}\wedge T} \dfrac{(v^{*}_s)^2}{\vol_s} ds + \zeta (Q^{*}_{\tau_{T-\delta,\ell}\wedge T})^2 \bigg| \mathcal F_{T-\delta}\right] .
\end{align*}
From \eqref{eq:optim_before_T_delta} and \eqref{eq:optim_after_T_delta}, we deduce the optimality of the defined process $Q^*$. 
\end{proof}

\section{PDE Analysis}\label{s:pde}

In this section, we will assume the price process to be Markovian and the cost structure to be a function of the price process; under these assumptions,
our goal is to relate the value function of the stochastic optimal control problem \eqref{e:oc1m2} to a PDE version of the BSDE
for the four choices of $I$ and ${\bm S}$ given in \eqref{e:I1}-\eqref{e:stopbelowLb}.

As noted in the introduction, we assume the price process $\bar{S}$ to be driven by a stochastic volatility model:
\[
\bar{S}_t =\int_0^t \sqrt{\nu_t} dW^{(1)}_t,
\]
where $\nu_t$ is the stochastic volatility process:
\begin{equation}\label{e:dynamics_nu}
d\nu_t = \alpha (\theta - \nu_t) dt + c \sqrt{\nu_t} dW^{(2)}_t;
\end{equation}
$W=(W^{(1)},W^{(2)})$ is a Brownian motion in ${\mathbb R}^2$ ($W^{(i)}_1$, $i=1,2$, have unit variance, but they are correlated with coefficient $\rho$).
We assume that the Feller condition, $2\alpha  \theta > c^2$, ensuring a positive process $\nu$ holds. Let $\mathcal L$ denote the second-order differential operator corresponding to these dynamics:
\begin{equation}\label{e:diff_operator}
{\mathcal L} u = \frac{1}{2}  \nu \partial^2_{ss} u+ \frac{1}{2}  \nu c^2\partial^2_{\nu\nu}u
+ \alpha (\theta- \nu) \partial_{\nu} u +  c \nu \rho  \partial^2_{s\nu}u.
\end{equation}
The variables $(\nu,s)$ belong to $D = (0,+\infty) \times  \mathbb R$. 

To get a PDE representation we take $\eta$ and $\Vol$ processes to be functions of $(\nu,\bar{S})$. With a slight abuse of notation,
we assume the market volume process to be $t\mapsto \Vol(t,\nu_t,\bar{S}_t)$ where $\Vol :[0,T]\times D \to \mathbb R_+$
is a non negative valued function; similarly the transaction cost process is $t\mapsto \eta(t,\nu_t,\bar{S}_t)$ where 
$\eta :[0,T]\times D \to (0,+\infty)$ is a strictly positive valued function. 
In this section we still use \eqref{e:Voltilde} and \eqref{e:vol}, we denote by $\Vol$ the quantity $\widetilde \Vol$, and will use the assumption \eqref{e:a1v2} on $\vol$:
\begin{equation}\label{eq:cond_prop_5_1_PDE}
0\leq \vol(t,s,\nu) \leq \overline{\vol},\qquad K^{p-1}T\overline{\vol} < \dfrac{1 }{p-1}
\end{equation}
where $K = k / \widehat p$ and $k$ is the constant in the permanent impact given by \eqref{e:kappavp}.

\subsection{PDE representation for $I=1,  \vect{S}=\{\bar{S}_T \geq \ell \}$}\label{ss:I1pde}

In this subsection we assume $I=1$, i.e., $\vol_t = (\widehat p-1) \Vol_t.$
Let us consider terminal values of the form
\begin{equation}\label{e:xipm}
	\xi = \Phi(\nu_T,\bar{S}_T);
\end{equation}
where 
\begin{equation}\label{e:Phi}
	\Phi : D \to [-K,\infty]
\end{equation}
is a measurable function. 
By \eqref{e:xip},
the choice ${\bm S} = \{\bar{S}_T \geq \ell \}$ corresponds to the $\Phi$ function
\[
	\Phi = \Phi^{(1)}(s) = -\frac{k}{\widehat p} 1_{(-\infty,\ell)}(s) +  \infty \cdot 1_{[\ell,\infty)}(s).
\]
Under the Markovian assumptions of the present section, and for $I=1$ the BSDE \eqref{e:convex_BSDE} and the terminal condition \eqref{e:xipm} correspond to the following PDE:
for any $(\nu,s) \in D$ and $t\in [0,T)$
\begin{equation}\label{e:HJB1}
\partial_t u + {\mathcal L}u - \mathrm{vol}_t |u|^p= 0. 
\end{equation}
with the terminal constraint
\begin{equation}\label{e:HJB1_term_cond}
u(T,\cdot,\cdot) = \Phi. 
\end{equation}
Our goal, under the Markovian assumptions of the present section, is to prove 
that the value function of the stochastic optimal control problem
\eqref{e:oc1m2} with $I=1$ and ${\bm S} = {\bm S}^{(1)}$ can be expressed
as a multiple of the unique solution of this PDE with terminal condition 
$\Phi=\Phi^{(1)}.$
To express the value function as a solution to the above PDE, 
we first extend the stochastic optimal control problem \eqref{e:oc1m2} to allow it to start from any time point $t$. Accordingly,
we define
\begin{align*}
\bar{S}^{t,\nu,s}_r & = s + \int_t^r \sqrt{\nu^{t,\nu,s}_u} dW^{(1)}_u, \\
\nu^{t,\nu,s}_r &=\nu +\int_t^r \alpha (\theta - \nu^{t,\nu,s}_u) du +\int_t^r c \sqrt{\nu^{t,\nu,s}_u} dW^{(2)}_u,\\
	\xi &= \Phi( \nu^{t,\nu,s}_T,\bar S^{t,\nu,s}_T),\qquad \vol^{t,\nu,s}_u = \vol(u, \nu^{t,\nu,s}_u,\bar S^{t,\nu,s}_u).
\end{align*}
Under our assumptions on $\Vol$ and $\Phi$, 
Proposition \ref{prop:sovlv_sing_BSDE} implies that the following BSDE
has a unique minimal supersolution:
\[
Y_r^{t,\nu,s} = \xi - \int_r^T\vol^{t,\nu,s}_v \left| Y^{t,\nu,s}_v \right|^p dv - \int_r^T Z^{t,\nu,s}_v dW_v ;
\]
set
\begin{equation}\label{e:valuefunI1}
	u^{\Phi}(t,\nu ,s) := Y_t^{t,\nu,s}. 
\end{equation}

As we already noted, our goal in this subsection is to prove that  $u^{\Phi}$ is the minimal supersolution (or the unique solution if $\Phi$ is finite) of the PDE \eqref{e:HJB1} for any $\Phi$ of the form \eqref{e:Phi} (and in particular for $\Phi = \Phi^{(1)}$). 
Our first step in this direction is the following:
\begin{lemma}
If $\Phi$ is continuous and with polynomial growth on $D$ and if ${\rm vo}l$ is also continuous on $[0,T]\times D$, 
then $u^\Phi$ is a continuous function of $(t,\nu ,s) \in [0,T] \times D$ and is the unique viscosity solution of the PDE \eqref{e:HJB1} 
with polynomial growth on $D$. 
\end{lemma}
\begin{proof}
See \cite[Theorem 5.37]{pard:rasc:14} or \cite[Theorems 3.4 and 3.5]{barl:buck:pard:97} (see also \cite{popi:17}). If $\Phi$ is bounded, the solution $Y^{t,\nu,s}$ and thus $u$ are also bounded. Hence our generator is Lipschitz continuous. 

If $\Phi$ satisfies $-k/\widehat p = -K \leq \Phi(s,\nu) \leq K (1 + |\nu|^m + |s|^m)$, then $\xi^+$ satisfies the condition imposed in Proposition \ref{thm:solve_aux_BSDE}. Thus $Y^{t,\nu,s}$ and thus $u$ are bounded from below, and we can modify our generator such that it becomes monotone (see the proof of Proposition \ref{thm:solve_aux_BSDE}). Then existence and uniqueness follow from \cite[Theorem 5.37]{pard:rasc:14}
	(this result is stated for $(\nu,s) \in \mathbb R^2$ in \cite{pard:rasc:14} but all arguments continue to work when $(\nu,s) \in D$).
\end{proof}

Now suppose that $\Phi$ is a continuous function from $D$ to $[-k/\widehat p,+\infty]$. We use the proof of Proposition \ref{prop:sovlv_sing_BSDE}. For any $n\geq 0$, we consider the bounded function $\Phi^n = \Phi \wedge n$. 
By the previous lemma
there exists a unique bounded viscosity solution $u^{\Phi \wedge n }$ and by comparison principle, 
\[
	u^{\Phi}(t,\nu,s) = \lim_{n\to +\infty} u^{\Phi\wedge n}(t,\nu,s)
\]
is well-defined with a bounded negative part. 
Suppose now that for some $m\geq 1$ and some $\epsilon > 0$: 
$$\forall (t,\nu,s)\in [T-\epsilon,T]\times D, \quad \dfrac{1}{\mathrm{vol}(t,\nu,s)} \leq C (1 +|\nu|^m+ |s|^m).$$
Then Condition \eqref{eq:cond_eta_sing} holds and we have on $[T-\epsilon,T]\times D$, 
\begin{align}  \nonumber
u^{\Phi \wedge n}(t,\nu,s) & \leq  \dfrac{1}{(T-t)^{\widehat p}} \mathbb E \left[ \int_t^T \left(  \dfrac{\widehat p-1}{\vol^{t,\nu,s}_\rho} \right)^{\widehat p-1} d\rho \right] \\\label{e:a_priori_estim_sol_PDE}
&  \leq  \dfrac{C}{(T-t)^{\widehat p-1}} (1+|\nu|^m+|s|^m)^{\widehat p-1}.
\end{align}
On the rest of the interval $[0,T]$, the bound of the solution $u^{\Phi\wedge n}$ is controlled by the previous estimate with $t = T-\epsilon$. In other words we have a bound on $u^{\Phi\wedge n}$ which does not depend on $n$. 
Hence $u$ is lower semi-continuous on $[0,T]\times D$ and finite (even locally bounded) on $[0,T)\times D$. These considerations
give us the following result:
\begin{lemma}
Suppose
$\Phi$ is a continuous function from $D$ to $[-k/\widehat p,+\infty]$; then
$u^{\Phi}$ 
is the minimal viscosity solution of the PDE \eqref{e:HJB1} on $[0,T)\times D$ (among all viscosity solutions with bounded negative part). 
\end{lemma}	
\begin{proof}
See \cite[Theorem 1]{popi:17}. 
\end{proof}
The concept of a viscosity solution
allows even a discontinuous solution and doesn't address the issue of smoothness/regularity of the solution; therefore the previous result doesn't
say anything about the regularity of $u^\Phi$.
The properties of the operator $\mathcal L$  (and smoothness assumptions on $\Vol$ and $\eta$)
allow us to establish the smoothness of $u^\Phi$, with a regularization bootstrap argument for parabolic PDE.
(see \cite[Lemma 5]{krus:popi:seze:18} for a similar argument).
\begin{lemma} \label{lem:regular_sol}
Suppose that $\mathrm{vol}$ is continuously differentiable with respect to all of its arguments. 
Assume that $\Phi_n$ is a sequence of continuous functions, converging to $\Phi$, such that the related (viscosity) solutions $u_n$ of the PDE  \eqref{e:HJB1} converge pointwise to $u$. Then $u$ belongs to $C^{1,2}([0,T) \times D)$ and is a classical solution of the PDE \eqref{e:HJB1}. 
\end{lemma}
\begin{proof}
Fix some $\epsilon > 0$ and $\mathfrak K$ a compact subset of $D$. 
First note that from \eqref{e:a_priori_estim_sol_PDE}, the bound of $u_n$ on $[0,T-\epsilon]\times \mathfrak K$ does not depend on the terminal value, that is on $n$, but only on $\epsilon$ and $\mathfrak K$. 

Moreover, the operator $\mathcal L$ can be written as follows:
\begin{align*}
{\mathcal L} u & = \frac{1}{2}  \nu \partial^2_{ss} u+ \frac{1}{2}  \nu c^2\partial^2_{\nu\nu}u
+ \alpha(\theta- \nu) \partial_{\nu} u +  c \nu \rho  \partial^2_{s\nu}u \\
& =\dfrac{1}{2} \mbox{div} \left( a( \nu, s) \nabla u\right) + b(\nu,s) \nabla u, \end{align*}
with
$$\nabla u(\nu,s) =  \left( \begin{array}{c}\partial_s u \\  \partial_\nu u \end{array} \right) $$
and 
$$a(\nu,s) =\nu  \left( \begin{array}{cc} 1 & c \rho \\  c \rho& c^2 \end{array} \right),\quad b(\nu,s) =  \left( \begin{array}{c}\dfrac{1}{2} +\dfrac{c\rho }{2} \\  \alpha(\theta-\nu) +\dfrac{c^2}{2}+\dfrac{c\rho }{2} \end{array} \right).$$
Our coefficients $a$ and $b$ are bounded on $\mathfrak K$, and $a$ is uniformly elliptic on $\mathfrak K$. Since $\vol$ is also continuously differentiable with respect to all of their arguments, then we can easily check that all conditions called a) (uniform ellipticity and boundedness condition), b) (growth condition on the derivatives, take $m=2$) and c) (regularity condition) of \cite[Theorem VI.4.4]{lady:solo:ural:68} hold. From this theorem, if there exists a function $\psi$ continuous on $[0,T-\epsilon] \times \mathfrak K$ and is of class $H^{1+\beta/2, 2+\beta}(]0,T-\epsilon[\times \overset{\circ}{\mathfrak K})$ for some $\beta > 0$  (space of functions which are $C^2$ with $\beta$-H\"older continuous second derivatives in the space variable $x$ and $C^1$ with $\beta/2$-H\"older continuous in the times variable $t$), then the PDE 
$$\partial_t v + \mathcal L v - \vol_t |v|^p = 0$$
with the boundary condition $u=\psi$, 
has a unique solution $v$ with the same regularity as $\psi$. Now our viscosity solutions $u_n$ are continuous and bounded on $[0,T-\epsilon]\times \mathfrak K$. Let us consider a sequence of smooth mollifiers $\zeta_m$ and define $\psi_m = u_n \star \zeta_m$. There exists a classical smooth solution $u_{n,m}$ of the PDE \eqref{e:HJB1} with boundary condition $\psi_m$ and pointwise $u_{n,m}$ converges to $u_n$ as $m$ goes to $+\infty$. As $u_n$, the bound of $u_{n,m}$ on $[0,T-\epsilon]\times \mathfrak K$ does not depend on $n$ and $m$.

In the next step, we prove that $u_{n,m}$ belongs to $H^{\beta/2,\beta}([0,T-\epsilon]\times \mathfrak K)$ (space of functions which $\beta$-H\"older continuous in the space variable $x$ and $\beta/2$-H\"older continuous in the times variable $t$) and that the norm of $u_{n,m}$ in this space does not depend on $n$ and $m$. 
Note that Conditions (1.2) and (7.1) of \cite[Chapter 3]{lady:solo:ural:68} are satisfied. Suppose that $v$ is a smooth solution of the PDE \eqref{e:HJB1} on $[0,T)\times D$: for any $(t,\nu,s)\in [0,T) \times D$
$$\partial_t v + \mathcal L v - \vol_t |v|^p = 0,$$
such that for any $\epsilon >0$ and any compact subset $\mathfrak K$ of $D$ , $v$ is bounded on $[0,T-\epsilon] \times \mathfrak K$. Thus $v$ solves on $[0,T-\epsilon] \times \mathfrak K$ the PDE
$$\partial_t v + \mathcal L v = \vol_t |v|^p = f,$$
where $f$ is a bounded function. We can apply \cite[Theorem III.10.1]{lady:solo:ural:68}. Hence $v$ is in $H^{\beta/2,\beta}([0,T-\epsilon]\times \mathfrak K)$. The value of $\beta > 0$ and the H\"older norm of $v$ depend on $\epsilon$, $\mathfrak K$ and the bound on $v$. In other words $\beta$ does not depend on the terminal value. 

In the last step, we get the desired regularity. Now we know that $u_{n,m}$ belongs to $H^{\beta/2,\beta}([0,T-\epsilon]\times \mathfrak K)$ and solves the PDE
$$\partial_t v + \mathcal L v = \vol_t |u_{n,m}|^p  = f_{n,m}.$$
Since $\vol$ is also in $C^1(D)$, then from \cite[Theorem IV.10.1]{lady:solo:ural:68}, $u_{n,m}$ is in $H^{1+\beta/2,2+\beta}([0,T-\epsilon']\times \mathfrak K')$ for any $\epsilon' > \epsilon$ and $\mathfrak K' \subset \mathfrak K$, and the norm depends only on the $H^\beta$-norm of $f_{n,m}$. 
Therefore, $u_n$, and thus $u$, belong to the same space\footnote{The Arzela-Ascoli theorem implies that $u_{n,m}$ (up to a subsequence) converges to some function $\widetilde u_n \in H^{\beta/2,\beta}( [0,T-\epsilon]\times \mathfrak K)$. Here $\widetilde u_n = u_n$ since pointwise convergence has been proved before.}, that is on any subset $[0,T-\epsilon]\times \mathfrak K$, $u_n$ and $u$ are in $C^{1,2}$. 
\end{proof}
The proof also shows that the regularity of any solution does not depend on the terminal value. In other words, far from $t=T$ and $\nu = 0$, the solutions are smooth and classical solutions. 

Let us summarize the foregoing results:
\begin{proposition}\label{prop:related_PDE}
Suppose that $\mathrm{vol}$ is continuously differentiable on $[0,T]\times D$, that $\Phi$ is bounded from below by $-k/\widehat p$ and that \eqref{e:a1v2} holds.
If one of the next conditions holds:
\begin{itemize}
\item $\Phi$ is continuous and with polynomial growth on $D$,
	\item $\Phi$ is continuous from $D$ to $[-k/\widehat p,+\infty]$ and for some $m\geq 1$ and $\epsilon > 0$  
$$\forall (t,\nu,s)\in [T-\epsilon,T]\times D,\quad \dfrac{1}{\mathrm{vol}(t,\nu,s)} \leq C (1 +|\nu|^m+ |s|^m),$$
\end{itemize}
then there exists a viscosity solution $u$ of the PDE \eqref{e:HJB1} with terminal value $\Phi(\cdot,\cdot)$. Moreover $u$ is of class $C^{1,2}([0,T) \times D)$ and is the minimal viscosity solution (among all viscosity solutions with bounded negative part). 
\end{proposition} 

\begin{remark}
The particular dynamics for $\bar S$ is not important, as soon as the related operator $\mathcal L$ regularizes the solution of the PDE. 
\end{remark}

\begin{remark}
Here we use the existence of a solution for the BSDE to deduce the existence of a viscosity solution for the PDE and
then we prove the regularity of the solution. 
The existence of a solution for the PDE could be directly proved (without probabilistic arguments) starting from 
\cite[Theorem V.8.1]{lady:solo:ural:68}. 
The proof would involve arguments similar to those given above because of the lack of monotonicity 
(see Assumption b) in  \cite[Theorem V.8.1]{lady:solo:ural:68} or Condition (4.17) in \cite[Chapter 7]{friedman}). 
Hence the proof of existence of solutions using PDE arguments will be as involved as the BSDE based approach given above.
\end{remark}

Let us next comment on the smoothness of $u$ on the boundary $\nu = 0$, called the hyperbolic part of the boundary. Recall that our operator $\mathcal L$ is defined by \eqref{e:diff_operator} and that the dynamics of $\nu$ is given by \eqref{e:dynamics_nu}. The Feller condition ensuring a positive process $\nu$ is $2\alpha \theta > c^2$. Under this condition, the Fichera function 
$$b(\nu) = \alpha(\theta- \nu) -\dfrac{1}{2} c^2 $$
is positive when $\nu$ goes to zero. Hence no boundary condition has to be supplied on $\nu = 0$ (see for example \cite{buck:zuza:gunt:16}). 

For the control problem, we have to consider the terminal condition
\[
\Phi_\ell(T,\nu,\cdot) = \Phi_\ell(\cdot)= \infty \cdot {\bm 1}_{[\ell,\infty)}(\cdot) - \frac{k}{\widehat p} {\bm 1}_{(-\infty, \ell)}(\cdot) .
\]
Note that we cannot directly apply Proposition \ref{prop:related_PDE}, since $\Phi_\ell$ is not continuous. Nonetheless
\begin{lemma} \label{lem:regul_arg_Phi1}
There exists a minimal viscosity solution $u_\ell$, which is of class $C^{1,2}$ on $[0,T) \times D$. 
\end{lemma}
\begin{proof}
Indeed let us define
$$\Phi^n(\cdot) = n \cdot {\bm 1}_{[\ell,\infty)}(\cdot) - \frac{k}{\widehat p} {\bm 1}_{(-\infty, \ell)}(\cdot)$$
and
\begin{align*}
\Phi^{n,m}(\cdot) & = n \cdot {\bm 1}_{[\ell,\infty)}(\cdot) - \frac{k}{\widehat p} {\bm 1}_{(-\infty, \ell-1/m]}(\cdot) \\
& + \left[(n+k/\widehat p) m(\cdot-\ell + 1/m) - k/\widehat p\right] {\bm 1}_{(\ell-1/m,\ell)}(\cdot) .
\end{align*}
$\Phi^{n,m}$ is continuous and non-decreasing with respect to $m$ and converges to $\Phi^n$. Therefore, the related continuous viscosity solutions $u^{n,m}$ converge to $u^n$. Arguing as in the proof of Proposition \ref{prop:related_PDE}, we obtain a uniform norm of $u^n$ in the space $H^{1+\beta,2+\beta}([0,T-\epsilon]\times \mathfrak K)$, for any compact subset $\mathfrak K$ of $D$. Then we pass on the limit on $n$ to obtain the desired result. 
Minimality can be obtained as for Proposition \ref{prop:related_PDE}. 
\end{proof}

Next we connect the solutions of PDE obtained above to the value function of the control problem.

\subsubsection*{Extended control problem}

The $L$ function corresponding to the permanent impact $\kappa$ given by \eqref{e:kappavp} is $L_t(\rho) = \eta_t |\rho|^{\widehat p}$ and the control problem \eqref{e:oc1m2} is
\begin{equation}\label{e:oc1m2p}
\inf_{Q \in {\mathcal A}_{I,\bm S}}  {\mathbb E}\left[ \int_0^T  \frac{1}{\Vol^{\widehat p-1}_t} |Q'_t|^{\widehat p}dt + \Phi(\nu_T,\bar S_T) |Q_T|^{\widehat p} \right].
\end{equation}
Here $I\equiv 1$ on $[0,T]$, with \eqref{e:Voltilde} and \eqref{e:vol} with the notation $\Vol =\widetilde \Vol$. 

For $p=2$, by Proposition \ref{prop:constrained_control_pb}, $u^{\Phi}(t,\nu,s) q^2$ is the value function of the extended version of the stochastic optimal control problem. With a slight abuse of language, we will refer to $u^\Phi$ simply as the value function of the extended stochastic optimal control problem with terminal cost $\Phi |Q_T|^{\widehat p}.$

These give us the main result of this subsection:
\begin{proposition} \label{prop:extended_control_pb_p}
Suppose \eqref{eq:cond_prop_5_1_PDE} and the conditions of Proposition \ref{prop:related_PDE} hold. $u^\Phi$ is the solution of the PDE \eqref{e:HJB1} with terminal condition \eqref{e:HJB1_term_cond}. 
If $\Phi$ is bounded and continuous or if $\Phi = \Phi^{(1)}$, then the value function of \eqref{e:oc1m2p} is given by $|q|^{\widehat p} u^\Phi$ and an optimal control is given by: 
\begin{equation} \label{e:optim_control_PDE}
v(r,\nu,s,q)= - (p-1)\mathrm{vol}(r,\nu,s) |u(r,\nu,s)|^{p-1}{\rm sgn}(u(r,\nu,s)) q.
\end{equation}
\end{proposition}
\begin{proof}
Since the proof is quite standard, let us provide an outline of the main arguments. First standard computations show that  for any $x$: 
$$\inf_v \left[ x v + \frac{1}{\Vol^{\widehat p-1}_t} |v|^{\widehat p} \right] = - (\widehat p -1) \Vol_t  \left(  \dfrac{|x|}{\widehat p}\right)^p$$
and that the infimum is attained at 
$$v^* = - \dfrac{\Vol_t }{(\widehat p )^{p-1}} |x|^{p-1} \sgn(x).$$

If $\Phi$ is continuous and bounded, by Proposition \ref{prop:related_PDE}, $u=u^\Phi$ is of class $C^{1,2}([0,T) \times D)$ and bounded and continuous on $[0,T] \times D$. 
Now define $V(t,\nu,s,q) = |q|^{\widehat p} u(t,\nu,s).$
This function has the same regularity as $u$ with respect to $(t,\nu,s)$ and is of class $C^{1}$ with respect to $q\in \mathbb R$. Take $Q=Q^w \in \contrset$ and $w=(Q^w)'$ is its derivative. Then It\^o's formula leads to: for any $t < T$
\begin{align*}
V(0,\nu,s,Q_0) & = V(t,\nu_t,\bar{S}_t,Q^w_t) - \int_0^t \left( \partial_t V + \mathcal L V  \right) (r,\nu_r,\bar{S}_r, Q^w_r) + \partial_q V(r,\nu_r,\bar{S}_r, Q^w_r) w_r dr \\
& - \int_0^t \sqrt{\nu_r} \left[ \partial_s V (r,\nu_r,\bar{S}_r, Q^w_r)   dW^{(1)}_r + \partial_\nu V (r,\nu_r\bar{S}_r,Q^w_r) c  dW^{(2)}_r \right] \\
& =  V(t,\nu_t,\bar{S}_t,Q^w_t) + \int_0^t \left[ - \vol_r | u  (r,\nu_r\bar{S}_r)|^p |Q^w_r|^{\widehat p} - \partial_q V(r,\nu_r,\bar{S}_r, Q^w_r) w_r \right] dr \\
& - \int_0^t \sqrt{\nu_r} \left[ \partial_s V (r,\nu_r,\bar{S}_r, Q^w_r)   dW^{(1)}_r + \partial_\nu V (r,\nu_r\bar{S}_r,Q^w_r) c  dW^{(2)}_r \right].
\end{align*}
Note that $ \partial_q V (r,\nu,s,q) = u(r,\nu,s) \widehat p |q|^{\widehat p -1} \sgn(q)$. Thus
\begin{align*}
-\vol_r  | u  (r,\nu_r\bar{S}_r)|^p |Q^w_r|^{\widehat p}  & = -\dfrac{\vol_r }{(\widehat p)^p} \left| \partial_q V (r,\nu_r\bar{S}_r,Q^w_r) \right|^p \\
& =- (\widehat p -1) \Vol_r  \dfrac{|\partial_q V (r,\nu_r\bar{S}_r,Q^w_r)|^{p}}{(\widehat p)^p} \\
& =  \inf_v \left[ \partial_q V (r,\nu_r\bar{S}_r,Q^w_r) v + \frac{1}{\Vol^{\widehat p-1}_t} |v|^{\widehat p} \right] \\
& \leq  \partial_q V (r,\nu_r\bar{S}_r,Q^w_r) v + \frac{1}{\Vol^{\widehat p-1}_t} |v|^{\widehat p}.
\end{align*}
Hence we obtain
\begin{align*}
V(0,\nu,s,Q_0) &  \leq  V(t,\nu_t,\bar{S}_t,Q^w_t) + \int_0^t  \frac{1}{\Vol^{\widehat p-1}_t} |w_r|^{\widehat p}  dr \\
& - \int_0^t \sqrt{\nu_r} \left[ \partial_s V (r,\nu_r,\bar{S}_r, Q^w_r)   dW^{(1)}_r + \partial_\nu V (r,\nu_r\bar{S}_r,Q^w_r) c  dW^{(2)}_r \right].
\end{align*}
Letting $t$ tend to $T$ and taking the expectation, we have
\begin{align*}
V(0,\nu,s,Q_0) &  \leq \mathbb E \left[ |Q^w_T|^{\widehat p} \Phi(\nu_T,\bar S_T) + \int_0^T  \frac{1}{\Vol^{\widehat p-1}_t} |w_r|^{\widehat p}  dr \right]. 
\end{align*}
Moreover if we take 
$$v^*_r = -(p-1) \vol_r  |u(r,\nu,s)|^{p-1}\sgn(u(r,\nu,s)) Q^*_r,$$
we have an equality. The first part of the proposition is proved. 

If $\Phi = \Phi^{(1)}$, then $u=u^{\Phi^{(1)}}$ is not continuous at time $T$, but $C^{1,2}$ on $[0,T) \times D$. Hence for any $t < T$:
\begin{align*}
V(0,\nu,s,Q_0) & = V(t,\nu_t,\bar{S}_t,Q^*_t) + \int_0^t  \frac{1}{\Vol^{\widehat p-1}_t} |v^*_r|^{\widehat p}  dr \\
& - \int_0^t \sqrt{\nu_r} \left[ \partial_s V (r,\nu_r,\bar{S}_r, Q^*_r)   dW^{(1)}_r + \partial_\nu V (r,\nu_r\bar{S}_r,Q^*_r) c  dW^{(2)}_r \right].
\end{align*}
We take the expectation and we use Fatou's lemma:
\begin{align*}
V(0,\nu,s,Q_0) & \geq \mathbb E \left[ \liminf_{t\to T} V(t,\nu_t,\bar{S}_t,Q^*_t) + \int_0^T  \frac{1}{\Vol^{\widehat p-1}_t} |v^*_r|^{\widehat p}  dr \right] \\
& = \mathbb E \left[ \liminf_{t\to T} \left( u(t,\nu_t,\bar{S}_t) |Q^*_t|^{\widehat p} \right) + \int_0^T  \frac{1}{\Vol^{\widehat p-1}_t} |v^*_r|^{\widehat p}  dr \right] .
\end{align*}
A direct computation shows that $N_t = \Ymin_t |Q^*_t|^{\widehat p-1} = u(t,\nu_t,\bar{S}_t) |Q^*_t|^{\widehat p -1}$ is a martingale on $[0,T)$. This martingale being bounded from below, it has a limit at time $T$. Arguing as at the end of the proof of Proposition \ref{prop:constrained_control_pb}, we deduce that 
\begin{align*}
V(0,\nu,s,Q_0) & \geq \mathbb E \left[   \Phi(\nu_T,\bar{S}_T) |Q^*_T|^{\widehat p}  + \int_0^T  \frac{1}{\Vol^{\widehat p-1}_t} |v^*_r|^{\widehat p}  dr \right] .
\end{align*}
Using the proof of Lemma \ref{lem:regul_arg_Phi1}, we obtain a sequence of smooth functions $u^{n,m}$ such that $|q|^{\widehat p} u^{n,m} $ is the value function of \eqref{e:oc1m2p} with $\Phi^{n,m}$ as terminal condition. Hence 
\begin{align*}
|Q_0|^{\widehat p} u^{n,m}(0,\nu,s) & \leq \inf_{Q\in\contrset} \mathbb E \left[   \Phi(\nu_T,\bar{S}_T) |Q_T|^{\widehat p}  + \int_0^T  \frac{1}{\Vol^{\widehat p-1}_t} |Q'_r|^{\widehat p}  dr \right] \\
& \leq V(0,\nu,s,Q_0) = |Q_0|^{\widehat p -1} u(0,\nu,s).
\end{align*}
Passing through the limit on $m \to +\infty$ and then on $n\to +\infty$ achieves the proof of this proposition. 
\end{proof}

An optimal state process $Q^*$ given by \eqref{e:optim_control_PDE} can be written:
$$(Q^*_r)' = v^*_r  = -(p-1) \vol_r  |\Ymin_r|^{p-1}\sgn(\Ymin_r)  Q^*_r.$$
And our proof shows that this control is admissible and 
$$\inf_{Q\in\contrset} \mathbb E \left[  \xi |Q_T|^{\widehat p}  + \int_0^T  \frac{1 }{\Vol^{\widehat p-1}_t} |Q'_r|^{\widehat p}  dr \right] \leq |Q_0|^{\widehat p} \Ymin_0.$$
However, even for bounded $\xi$, and except in the case $p=\widehat p =2$, the lack of convexity prevents us proving equality. Therefore the Markovian setting and the fact that the solution $u$ of the HJB equation is smooth, are crucial to extend our result to any $\widehat p>1$.

\subsection{PDE representation for $I_t=\mathbf 1_{\{t \leq \tau_\ell \}}$ and $\vect{S} = \{\tau_\ell \geq T \}$}\label{ss:pdexi2}

To get a PDE representation of the BSDE \eqref{e:convex_BSDE}, \eqref{e:termcond} for $I_t= I^{(2)}_t = \mathbf 1_{\{t \le \tau_\ell \}}$, we consider the
problem in the interval $[\![0,\tau_\ell \wedge T]\!]$.  As discussed in subsection \ref{ss:red},
the corresponding BSDE is again \eqref{e:convex_BSDE} 
but with terminal condition \eqref{e:xi2}.
This reduced formulation of the problem is indeed Markovian.
The PDE is the same as before \eqref{e:HJB1} but solved over the domain
$D=(0,\infty)\times [\ell,\infty)$ and with boundary conditions
\begin{equation}\label{e:bcI2}
	u|_{[0,T) \times (0,\infty) \times \{\ell\}} = -k/\widehat p, \qquad	u|_{ \{T\}\times (0,\infty) \times [\ell,\infty)} = \infty.
\end{equation}

The value function $u^{(2)}$ of the extended version of the stochastic control problem is again
defined through \eqref{e:valuefunI1} and we have:
\begin{proposition} \label{prop:pde_case_2}
$u^{(2)}$ is the minimal viscosity solution of \eqref{e:HJB1} and boundary conditions \eqref{e:bcI2}. Moreover $u^{(2)}$ is of class $C^{1,2}$ on $[0,T) \times (0,\infty) \times (\ell,\infty)$. 
\end{proposition}
The proof proceeds parallel to the arguments given in the previous section. We therefore provide an outline.
We begin by considering the case where the boundary condition is given by
\begin{equation} \label{e:bcI2m1}
u|_{[0,T) \times (0,\infty) \times \{\ell\}} = \psi ,\qquad u|_{ \{T\}\times (0,\infty) \times [\ell,\infty)} = \Phi
\end{equation}
where $\psi$ is a continuous and bounded function. 
If $\Phi$ is bounded and if the compatibility constraint $\Phi(\nu,\ell) = \psi(T,\nu)$ is verified, we can directly apply \cite[Theorem 5.41]{pard:rasc:14} to obtain the existence of a unique viscosity bounded and continuous solution $u$ of the PDE \eqref{e:HJB1} with 
the boundary condition \eqref{e:bcI2m1}. The regularity inside the domain can be obtained by the arguments of Lemma \ref{lem:regular_sol}.
If the compatibility condition does not hold, the solution still exists but is not continuous up to the boundary. Finally,
the $\infty$ terminal condition can be handled via approximation from below (as was done in the previous section as well
as in Section \ref{s:bsde} in the treatment of the BSDE \eqref{e:convex_BSDE} and the singular terminal condition \eqref{e:termcond}).

Using Proposition \ref{prop:control_pb_stop_time}, the value function of the control problem \eqref{e:oc1m2red} is: $V(t,\nu,s,q) = |q|^2 u(t,\nu,s)$ for any $s \geq \ell$. But the regularity of $u^{(2)}$ also allows us to solve the control problem \eqref{e:oc1m2p} up to $\tau_\ell \wedge T$:
$$|Q_0|^{\widehat p} u (0,\nu,s) = \min_{Q \in {\mathcal A}_{I,\bm S}}  {\mathbb E}\left[ \int_0^{\tau_\ell \wedge T}  \frac{1}{\Vol^{\widehat p-1}_t} |Q'_t|^{\widehat p}dt + \Phi(\nu_{\tau_\ell \wedge T},\bar S_{\tau_\ell \wedge T}) |Q_{\tau_\ell \wedge T}|^{\widehat p} \right].$$
The arguments are the same as in the proof of Proposition \eqref{prop:extended_control_pb_p}, where $T$ is replaced by $\tau_\ell \wedge T$.

\subsection{PDE representation for $I_t = 1_{\{\bar{S}_t \geq \ell \}}$ and $\vect{S}= \{ \tau_{T-\delta,\ell} \geq  T \}$}\label{ss:pdeI4}

In Section \ref{ssect:cases_3_4}, we already explain how to reduce the problem on the interval $[0,T-\delta]$. Let us consider the PDE 
$$\partial_t u + \mathcal L u  - \vol_t  |u|^p= 0$$
on the set $[0,T] \times (0,+\infty) \times [\ell,+\infty)$, with terminal condition $+\infty$ at time $T$ and $-k/\widehat p$ on the lateral boundary $\{s=\ell\}$. From the previous section, there exists a unique solution $u^\infty$ defined on $[0,T)\times (0,+\infty) \times [\ell,+\infty)$. 

Following the representation given by \eqref{e:termcdred34}, we define
$$r(t,s,v) = |v|^p \vol_t \mathbf 1_{ s  >\ell },$$
and solve
\begin{equation}\label{e:pde3}
\partial_t u + \mathcal L u  - r(t,s,u) = 0
\end{equation}
over $[0,T-\delta] \times (0,+\infty) \times \mathbb R$
with terminal boundary condition
$$g(\nu,s) = u^\infty (T-\delta,\nu,s) \mathbf 1_{s >\ell} - (k/\widehat p) \mathbf 1_{s \leq \ell}.$$ 
Note that the terminal boundary condition $g$ is bounded and continuous.  Nonetheless the free term $r$ is not continuous at $s=\ell$. 
\begin{lemma} \label{lem:PDE_second_case}
There exists a function $v$ such that $v$ is bounded and continuous on $[0,T-\delta] \times (0,+\infty) \times \mathbb R$ and is a solution of class $C^{1,2}$ of the PDE \eqref{e:pde3} on $[0,T-\delta) \times (0,+\infty) \times (\mathbb R \setminus \{\ell \})$. 
\end{lemma}
\begin{proof}
To circumvent the discontinuity of $r$, let us introduce 
$$\phi^\epsilon(s) =  \left(  \mathbf 1_{ s  >\ell+\epsilon } + \dfrac{(s-\ell)}{\epsilon} \mathbf 1_{ \ell < s  \leq \ell+\epsilon } \right),\quad r^{\epsilon} (t,s,v) =  |v|^p \vol_t \phi^\epsilon(s).$$
This function is Lipschitz continuous with respect to $s$, satisfies $r^\epsilon \leq r$ and converges increasingly and pointwise to $r$ when $\epsilon$ tends to zero. 

From standard arguments (see \cite[Theorem 5.37]{pard:rasc:14}), there exists a unique bounded and continuous viscosity solution $v^\epsilon$ of the PDE 
$$\partial_t u + \mathcal L u  - r^\epsilon(t,s,u) = 0$$
with the same terminal condition $g$ as $v$. Note that the bounds on $v^\epsilon$ do not depend on $\epsilon$. Thus arguing as in Lemma \ref{lem:regular_sol}, we can prove that $v^\epsilon$ is of class $C^{1,2}$ on $[0,T-\delta) \times (0,+\infty) \times (\mathbb R \setminus \{\ell\})$ with a norm independent of $\epsilon$. 

The comparison principle shows that $v^\epsilon$ is a decreasing sequence and thus we can define $v^\star$ as the decreasing limit of $v^\epsilon$ as $\epsilon$ tends to zero. We obtain immediately that $v^\star$ is bounded and upper semi-continuous and is a viscosity subsolution of PDE \eqref{e:pde3} (well-known result on stability for viscosity solutions \cite{cran:ishi:lion:92}). 

The only remaining point concerns the continuity of $v^\star$ on the set $\{s = \ell\}$. Let us define another approximating sequence $w^\epsilon$ defined as the solution of PDE \eqref{e:pde3} where $r$ is replaced by $\widetilde r^\epsilon$:
$$\psi^\epsilon(s) =  \left(  \mathbf 1_{ s  >\ell } + \dfrac{(s-\ell)}{\epsilon} \mathbf 1_{ \ell-\epsilon < s  \leq \ell } \right),\quad \widetilde r^{\epsilon} (t,s,v) =  |v|^p \vol_t \psi^\epsilon(s).$$
As $v^\epsilon$, $w^\epsilon$ converges to $v_\star$, which is lower semi-continuous and is a viscosity supersolution of PDE \eqref{e:pde3}. Moreover by comparison principle, $w^\epsilon \leq v_\star \leq v^\star \leq v^\epsilon$. Comparing sub- and supersolution implies that $v_\star = v^\star$ (standard result for viscosity solution). 
Let us prove this statement in our case. For any $(\nu,s)$ we have
\begin{align*}
& v^\epsilon(t,\nu,s) - w^\epsilon(t,\nu,s) = Y^{\epsilon,t,\nu,s}_t - \widetilde Y^{\epsilon,t,\nu,s}_t \\
& = Y^{\epsilon,t,\nu,s}_{T-\delta} -  \widetilde Y^{\epsilon,t,\nu,s}_{T-\delta} - \int_t^{T-\delta} r^{\epsilon} (u,\bar S^{t,\nu,s}_u,Y^{\epsilon,t,\nu,s}_u)  - \widetilde r^{\epsilon} (u,\bar S^{t,\nu,s}_u, \widetilde Y^{\epsilon,t,\nu,s}_u)  du \\
& \quad - \int_t^{T-\delta} ( Z^{\epsilon,t,\nu,s}_u -  \widetilde Z^{\epsilon,t,\nu,s}_u ) dW_u \\
& = - \int_t^{T-\delta} \vol(u,\nu^{t,\nu,s}_{u},\bar S^{t,\nu,s}_{u}) \phi^\epsilon(\bar S^{t,\nu,s}_u) (Y^{\epsilon,t,\nu,s}_u  -  \widetilde Y^{\epsilon,t,\nu,s}_u )  h (Y^{\epsilon,t,\nu,s}_u , \widetilde Y^{\epsilon,t,\nu,s}_u )  du \\
& \quad - \int_t^{T-\delta} \vol(u,\nu^{t,\nu,s}_{u},\bar S^{t,\nu,s}_{u}) \left| \widetilde Y^{\epsilon,t,\nu,s}_u \right|^p \left(  \phi^\epsilon(\bar S^{t,\nu,s}_u)  - \widetilde \phi^\epsilon(\bar S^{t,\nu,s}_u)\right) du \\
& \quad - \int_t^{T-\delta} ( Z^{\epsilon,t,\nu,s}_u -  \widetilde Z^{\epsilon,t,\nu,s}_u ) dW_u.
\end{align*}
Here 
$$h(y,\widetilde y) =p \int_0^1 |\widetilde y + \alpha (y - \widetilde y) |^{p-1} \sgn(\widetilde y + \alpha (y - \widetilde y)) d\alpha.$$
Using the boundedness of $Y^{\epsilon,\cdot,\cdot,\cdot}$ and $ \widetilde Y^{\epsilon,\cdot,\cdot,\cdot}$ (uniformly with respect to $\epsilon$) and standard stability result for BSDEs, we obtain the existence of a constant $C$ independent of $\epsilon$ such that 
\begin{align*}
& |v^\epsilon(t,\nu,s) - w^\epsilon(t,\nu,s)|^2  \leq  C\mathbb E \left[  \int_t^{T-\delta} \left(  \phi^\epsilon(\bar S^{t,\nu,s}_u)  - \widetilde \phi^\epsilon(\bar S^{t,\nu,s}_u)\right)^2 du \right] \\
& \leq  2 C \int_t^{T-\delta}  \mathbb P \left( \ell-\epsilon \leq \bar S^{t,\nu,s}_u  \leq \ell+\epsilon  \right) du 
\end{align*}
Fix some $\eta > 0$. The uniform ellipticity of $\mathcal L$ implies that there exists $\epsilon_0$ such that for any $\epsilon < \epsilon_0$, 
$$ \mathbb P \left( \ell-\epsilon \leq \bar S^{t,\nu,s}_u  \leq \ell+\epsilon \right) \leq  \eta^2 / (2C).$$
Hence letting $\epsilon$ go to zero, we get for any $\eta > 0$
$$ |v^\star(t,\nu,s) - v_\star(t,\nu,s)| \leq \eta.$$
Thus $v^\star = v_\star$ and thus $v$ is continuous. 

Finally, by a regularization argument, it is a classical solution of the PDE \eqref{e:pde3} on $[0,T-\delta) \times (0,+\infty) \times (\mathbb R \setminus \{\ell\})$. 

\end{proof}

If we define $Y_t =  v(t,\nu_t,\bar S_t)$, using It\^o's formula (allowed since $\bar S_t \neq \ell$ a.s.) we check that $Y=\Ymin$ is the solution of BSDE  \eqref{e:convex_BSDE} with terminal condition \eqref{e:termcdred34}.

Now since $v$ is a smooth function, $|q|^{\widehat p} v(t,\nu,s)$ is the value function of the control problem \eqref{e:oc1m2p} with $I_t = 1_{\{\bar{S}_t \geq \ell \}}$ and $\vect{S}= \{ \tau_{T-\delta,\ell} \geq  T \}$ and 
\begin{align*}
Q^*_t & = Q_0 \exp\left( -\int_{0}^{t} (p-1)  \vol (s,\nu_s,\bar S_s)|v(s,\nu_s,\bar S_s)|^{p-1}\sgn(v(s,\nu_s,\bar S_s)) ds \right) \\
& =  Q_0 \exp\left( -\int_{0}^{t} (p-1)\vol (s,\nu_s,\bar S_s) |\Ymin_s|^{p-1}\sgn(\Ymin_s) ds \right) .
\end{align*}
The proof is a direct adaptation of the proof of Proposition \ref{prop:extended_control_pb_p} to the current setting.

\subsection{PDE representation for $I^{(4)}$ and $\vect{S}^{(4)}$}\label{ss:pdeI3}

Finally, let us derive a PDE representation for
the control problem \eqref{e:oc1m2p} with $I = I^{(4)}$ and ${\bm S}={\bm S}^{(4)}$. 
From Section \ref{ssect:cases_3_4}, we know that the problem can be reduced to the time interval $[0,T-\delta]$;
as opposed to what happens with $I=I^{(3)}$ and ${\bm S} = {\bm S}^{(3)}$, the problem is not Markovian after this reduction. This is because,
the choice $I = I^{(4)}$ and ${\bm S}={\bm S}^{(4)}$ introduces an additional state variable, which is $I^{(4)}$ itself; $I^{(4)}$
keeps track of whether the system is in the trading state or in the waiting state (for $I_t=I^{(3)}_t = \mathbf 1_{[\ell,\infty)}(W_t)$
$I$ is directly a function of $W$ and $W$ serves as the state of the system).
Correspondingly, 
we expect a value function of the form
$V = u_i(t,\nu,s)|q|^{\widehat p}$ where
the additional variable $i$ (shown as a subscript) takes values in 
$\{1,0\}$
($1$ for the trading state and $0$ for no trading) 
and $u$ satisfies the coupled PDE system
\begin{equation}\label{e:pde1intro}
\partial_t u_0 +{\mathcal L } u_0 = 0 
\end{equation}
where $\mathcal L$ is defined by \eqref{e:diff_operator}
	solved in the region $[0,T-\delta] \times (0,\infty) \times (-\infty, \ell + b)$ with boundary conditions
$$ u_{0}(T-\delta,\nu,s) = -\frac{k}{\widehat p}, \quad u_{0}(t,\nu,\ell+b) = u_{1}(t,\nu,\ell+b)$$
and
	\begin{equation}
	\partial_t u_{1} + {\mathcal L } u_{1}  - \vol_t |u_{1}|^{p}=0, \label{e:pde2intro}
\end{equation}
		solved in the region $[0,T] \times (0,\infty) \times (\ell,\infty)$ with boundary conditions
\begin{align*}
	u_{1}(T,\nu,s) &= \infty, \notag \\
	u_{1}(t,\nu,\ell) &= u_{0}(t,\nu,\ell), \ t < T-\delta, \notag \\
	u_{1}(t,\nu,\ell) &= -\frac{k}{\widehat p}, \ t \geq T-\delta. \notag
\end{align*}
Note that the first equation (corresponding to the waiting state) is linear since in this state no control is applied
and the PDE is determined only by the underlying diffusion.
We think that a solution to this system can be obtained from the minimal supersolution of the corresponding BSDE as we did
in earlier sections.
In the following
subsection we provide an alternative solution based on the number of switches between trading
and waiting states. The sequence of control problems where ${\mathcal N}$ is bounded by $n$ can be solved recursively,
letting $n\rightarrow \infty$ gives a solution to the above system. 
An advantage of this approach is that it also gives a numerical algorithm to compute the value function. The difficulty in trying
to numerically solve the above system directly is that $u_1$ and $u_0$ appear as boundary conditions in the equations that the other satisfies.

\subsubsection*{Finite trading approximation}

Define
\[
	{\mathcal N}_t \doteq (\sup\{k: \bar\tau_{b,k} < t \} \vee (-1)) + \mathbf 1_{\tau_\ell > 0}
\]
The number of trading intervals realized up to terminal time $t$ is equal to ${\mathcal N}_t +1 $. 
Recall that for $I = I^{(4)}$ the set of all trading times before terminal time $T$ is
\[
	{\mathcal I} = \{ t \in [0,T], \ I_t = 1\} = \bigcup_{k=1}^{{\mathcal N}_T+1} [\![ \underline{\rho}_k ,\overline{\rho}_k ]\!], 
\]
and ${\bm S}^{(4)} = \{ T \in \mathcal I\}$. Define 
$${\mathcal I}_n= \bigcup_{k=1}^{n}  [\![ \underline{\rho}_k ,\overline{\rho}_k ]\!], \quad \varpi_n = \sup {\mathcal I}_n, \qquad n=1,2,3,...$$
Note that $\varpi_n$ is one of the stopping times $\tau_{\ell,k}$. 
The control problem with a limit on the number of active intervals is:
\begin{equation*}
V_{n}(Q_0,\nu,s) = \inf_{Q \in \contrset^n} {\mathbb E}\left[ \int_0^T
\frac{|Q'_s|^{\widehat p}}{\vol_s} ds -\frac{k}{\widehat p} |Q_{T}|^{\widehat p}\right], 
\end{equation*}
where
\begin{align*}
\contrset^n&= \{Q: {\mathbb F}\text{-adapted, absolutely continuous}, \\
&~~~~~~Q_T = 0\text{ if } T \in {\mathcal I}_n,\ \quad Q_t' = 0, \  \text{ if }t\not \in {\mathcal I}\text{ or } t > \varpi_n \}.
\end{align*}

\begin{lemma} \label{lem:approx_val_fct}
The sequence $V_n(Q_0,\nu,s) $ is non-increasing and tends to $V(Q_0,\nu,s)$, the value function of the
	control problem \eqref{e:oc1m2p} with $(I,{\bm S}) = (I^{(4)}, {\bm S}^{(4)})$, as $n$ tends to $\infty.$
\end{lemma}
\begin{proof}
Note that $\contrset^n$ is included in $\contrset^{n+1}$. Hence $V_n$ is a non-increasing sequence. Moreover taking $Q \in {\mathcal A}_{I,{\bm S}}$ and defining $\widetilde Q$ equal to $Q$ on the random interval $[\![ 0  ,\varpi_n]\!]$ and $\widetilde Q'$ equal to zero after $\varpi_n$, the strategy $\widetilde Q$ belongs to $\contrset^n$ and 
\begin{align*} 
\int_0^{T} \frac{|Q'_s|^{\widehat p}}{\vol_s}ds -\frac{k}{\widehat p}|Q_{T}|^{\widehat p} & = \int_0^{\varpi_n} \frac{|Q'_s|^{\widehat p}}{\vol_s}ds  +  \int_{\varpi_n}^T \frac{|Q'_s|^{\widehat p}}{\vol_s}ds  -\frac{k}{\widehat p}|Q_{T}|^{\widehat p} \\
& =  \int_0^{\varpi_n} \frac{|\widetilde{Q}'_s|^{\widehat p}}{\vol_s}ds -\frac{k}{\widehat p}|\widetilde Q_{T}|^{\widehat p}  +  \int_{\varpi_n}^T \frac{|Q'_s|^{\widehat p}}{\vol_s}ds -\frac{k}{\widehat p}( |Q_{T}|^{\widehat p}- |\widetilde Q_{T}|^{\widehat p}) .
\end{align*}
Taking the expectation, we have
\begin{align*} 
\mathbb E \left[ \int_0^{T} \frac{|Q'_s|^{\widehat p}}{\vol_s}ds -\frac{k}{\widehat p}|Q_{T}|^{\widehat p} \right] &  = \mathbb E \left[ \int_0^{\varpi_n} \frac{|\widetilde{Q}'_s|^{\widehat p}}{\vol_s}ds -\frac{k}{\widehat p}|\widetilde Q_{T}|^{\widehat p}  \right] \\
& + \mathbb E \left[ \int_{\varpi_n}^T \frac{|Q'_s|^{\widehat p}}{\vol_s}ds -\frac{k}{\widehat p}( |Q_{T}|^{\widehat p}- |\widetilde Q_{T}|^{\widehat p}) \right] \\
& \geq V_n(Q_0,\nu,s) +  \mathbb E \left[ \int_{\varpi_n}^T \frac{|Q'_s|^{\widehat p}}{\vol_s}ds -\frac{k}{\widehat p}( |Q_{T}|^{\widehat p} -| \widetilde Q_{T}|^{\widehat p}) \right].
\end{align*}
Since $ \widetilde Q_{T} = Q_{\varpi_n}$, we have
\begin{align*} 
 \mathbb E \left[ |Q_{T}|^{\widehat p} - |\widetilde Q_{T}|^{\widehat p} \right] & = \mathbb E \left[ \left(  \int_{\varpi_n}^T Q'_s ds \right) \int_0^1 \widehat p |\widetilde Q_{T} + \alpha (Q_T - \widetilde Q_{T}) |^{\widehat p-1} \sgn(\widetilde Q_{T} + \alpha (Q_T - \widetilde q_{T})) d\alpha \right] \\
 & \leq \left[   \mathbb E \left( \int_0^1 \widehat p |\widetilde Q_{T} + \alpha (Q_T - \widetilde Q_{T}) |^{\widehat p-1}  d\alpha \right)^p \right]^{\frac{1}{p}}\left[   \mathbb E\left(  \int_{\varpi_n}^T Q'_s ds \right)^{\widehat p}  \right]^{\frac{1}{\widehat p}}  \\
 & \leq C T^{\frac{1}{\widehat p}}  \left[  \mathbb E  \left(  |\widetilde Q_{T}  |^{\widehat p} +|Q_T|^{\widehat p} \right) \right]^{\frac{1}{p}} \left[   \mathbb E\left(  \int_{\varpi_n}^T \dfrac{|Q'_s|^{\widehat p}}{\vol_s} ds \right) \right]^{\frac{1}{\widehat p}} 
 \end{align*}
 under Assumption  \ref{a:vol} on $\vol$.  Let us notice that
 $$   \mathbb E  \left(  |\widetilde Q_{T}  |^{\widehat p} +|Q_T|^{\widehat p} \right)  \leq 2^{\widehat p}   \mathbb E  \left( |Q_0|^{\widehat p} + \int_0^T |Q'_s|^{\widehat p} ds \right).$$
 Thereby we obtain for any strategy $Q \in \contrset$ and any $n$:
 \begin{align*} 
& \mathbb E \left[ \int_0^{T} \frac{|Q'_s|^{\widehat p}}{\vol_s}ds -\frac{k}{\widehat p}|Q_{T}|^{\widehat p} \right]  \geq V_n(Q_0,\nu,s) +  \mathbb E \left[ \int_{\varpi_n}^T \frac{|Q'_s|^{\widehat p}}{\vol_s}ds \right] \\
&\qquad  - kC T^{\frac{1}{\widehat p}}   \left[   \mathbb E  \left( |Q_0|^{\widehat p} + \int_0^T |Q'_s|^{\widehat p} ds \right) \right]^{\frac{1}{p}} \left[   \mathbb E\left(  \int_{\varpi_n}^T \dfrac{|Q'_s|^{\widehat p}}{\vol_s} ds \right) \right]^{\frac{1}{\widehat p}}  \\
& \qquad \geq V(Q_0,\nu,s)  +  \mathbb E \left[ \int_{\varpi_n}^T \frac{|Q'_s|^{\widehat p}}{\vol_s}ds \right] \\
&\qquad  - kC T^{\frac{1}{\widehat p}}   \left[   \mathbb E  \left( |Q_0|^{\widehat p} + \int_0^T |Q'_s|^{\widehat p} ds \right) \right]^{\frac{1}{p}} \left[   \mathbb E\left(  \int_{\varpi_n}^T \dfrac{|Q'_s|^{\widehat p}}{\vol_s} ds \right) \right]^{\frac{1}{\widehat p}}.
\end{align*}
From the monotone convergence theorem we deduce
 \begin{align*} 
&  \mathbb E \left[ \int_0^{T} \frac{|Q'_s|^{\widehat p}}{\vol_s}ds -\frac{k}{\widehat p}|Q_{T}|^{\widehat p} \right] 
 \geq \lim_{n} V_n(Q_0,\nu,s) \geq V(Q_0,\nu,s) .
\end{align*}
Since this holds for any $Q \in \contrset$, we obtain the desired result. 
\end{proof}

Now we consider the dynamical version of the control problem with $n$ trading intervals: for any $t \in [0,T]$:
\[
	V_{i,n}(t,q,\nu,s) = \inf_{Q \in \contrset^{n,t}} {\mathbb E}\left[ \int_t^T
\frac{|Q'_s|^{\widehat p}}{\vol_s} ds -\frac{k}{\widehat p}|Q_{T}|^{\widehat p} \right], 
\]
the value function of the control problem starting at an arbitrary time $t < T$.
The additional variable $i$ indicates the starting value of the process $I$; $I_t = 1$ means that the problem starts from
a trading state and $I_t = 0$ means that the problem starts from a waiting state.
The definition of $ {\mathcal A}^{n,t}$ is the same as ${\mathcal A}^{n}$: there are at most $n$ trading intervals during the time interval $[t,T]$.

Considering the problem until the first transition from one state to the other (from trading to no trading or vice versa)
we can write the above optimal control problems recursively as follows:
\begin{align}  \label{e:rec1}
	V_{0,n}(t,q,\nu,s) & = {\mathbb E}\left(
		V_{1,n}(\bar\tau_{b,j_t},q,\nu_{\bar\tau_{b,j_t}},\ell+b) {\bm 1}_{\{\bar\tau_{b,j_t} < T\}}
- \dfrac{k}{\widehat p} |q|^{\widehat p} {\bm 1}_{\{\bar\tau_{b,j_t} \ge  T\}}\right) 
\end{align}
where $j_t = \inf \{j, \ \bar\tau_{b,j} > t\}$ and
\begin{align} \nonumber 
	V_{1,n}(t,q,\nu,s) & =
	\inf_{Q \in \contrset^{n,t} } {\mathbb E}
\left(
\int_t^{\tau_{\ell,j_t} \wedge T} 
\frac{|Q'_s|^{\widehat p}}{\vol_s}ds
	- \dfrac{k}{\widehat p} |Q_T|^{\widehat p} {\bm 1}_{\{\tau_{\ell,j_t} \ge  T\}} \right. \\ \label{e:rec2}
&\qquad\qquad\qquad \left.  +V_{0,n-1}(\tau_{\ell,j_t},Q_{\tau_{\ell,j_t}},\nu_{\tau_{\ell,j_t}},\ell) {\bm 1}_{\{\tau_{\ell,j_t} < T\}}\right) 
\end{align}
where $\tau_{\ell,j_t}$ is the first time after $t$ when the price $\bar S$ goes below $\ell$. 
In \eqref{e:rec1} the problem starts in a waiting interval; the controller waits until the first time after $t$ when $\bar S$ goes above $\ell+b$
($\bar\tau_{b,j_t}$)
or
$T$, whichever comes first.
If it is $\bar\tau_{b,j_t}$,
 trading starts; if it is $T$, the controller pays the terminal cost.
Note that the recursion \eqref{e:rec1}  involves no control since the liquidation process starts in the waiting state (i.e., no trading, $Q'=0$) 
and remains in that state until $\bar\tau_{b,j_t}.$ In \eqref{e:rec2}
the agent already uses one trading possibility, after $\tau_{\ell,j_t}$, there are at most $n-1$ trading intervals.

\bigskip

As before the homogeneous cost structure suggests
\[
V_{i,n}(t,q,\nu,s) = |q|^{\widehat p} u_{i,n}(t,\nu,s),
\]
and we define
$u_{0,0}$ to be the constant function $-\frac{k}{\widehat{p}}$: there is no trading (i.e., no control and $Q'= 0$) 
and at terminal time the trader pays the terminal cost.
The above recursions imply the following sequence of PDE to compute $u_{0,n}$ and $u_{1,n}$ for $n\geq 1$:
\begin{equation}\label{e:pde1}
\partial_t u_{0,n} +{\mathcal L } u_{0,n} = 0,
\end{equation}
where $\mathcal L$ is defined by \eqref{e:diff_operator}
solved in the region $[0,T-\delta] \times (0,\infty) \times (-\infty, \ell + b)$ with boundary conditions
\begin{equation}\label{e:bcu0n}
u_{0,n}(T-\delta,\nu,s) = -\frac{k}{\widehat p}, \quad u_{0,n}(t,\nu,\ell+b) = u_{1,n}(t,\nu,\ell+b).
\end{equation}
And
\begin{equation}\label{e:pde2}
\partial_t u_{1,n} + {\mathcal L } u_{1,n}  - \vol_t |u_{1,n}|^{p} = 0,
\end{equation}
solved in the region $[0,T] \times (0,\infty) \times (\ell,\infty)$ with boundary conditions
\begin{align*}
& u_{1,n}(T,\nu,s) = \infty, \\
& u_{1,n}(t,\nu,\ell) = u_{0,n-1}(t,\nu,\ell), \ t < T-\delta, \\
& u_{1,n}(t,\nu,\ell) = -\frac{k}{\widehat p}, \ t \geq T-\delta.
\end{align*}
The PDE for $u_{1,n}$ and its boundary conditions on the set $[T-\delta,T]\times (0,+\infty)\times[\ell,\infty)$ do not depend on $n$. 
Hence on this set, we have $u_{1,n} = u^\infty$ (see Section \ref{ss:red} and the beginning of Section \ref{ss:pdeI4}). 

For $n=1$ we have: $u_{1,1} = u^{(2)}$ where $u^{(2)}$ is the value function in Section \ref{ss:pdexi2}, Proposition \ref{prop:pde_case_2}. 
The arguments of Lemma \ref{lem:regular_sol} show that $u_{1,1}$ is of class $C^{1,2}$ on $[0,T-\delta] \times (0,\infty) \times (\ell,\infty)$ and 
continuous on $[0,T-\delta] \times (0,\infty) \times [\ell,\infty)$. 
Once $u_{1,1}$ is available, the rest of the value functions can be computed recursively by solving \eqref{e:pde1} and \eqref{e:pde2} 
in the following order:
\begin{equation*}
u_{1,1} \rightarrow u_{0,1} \rightarrow u_{1,2} \rightarrow u_{0,2}
\cdots
\end{equation*}
Note that $u_{0,1}$ solves the linear PDE \eqref{e:pde1} with smooth and bounded boundary conditions $-k/\widehat p$ and $u_{1,1}(\cdot,\cdot,\ell+b)$. 
Therefore, $u_{0,1}$ is also of class $C^{1,2}$ on $[0,T-\delta] \times (0,\infty) \times (-\infty,\ell+b)$ and continuous on $[0,T-\delta) \times (0,\infty) \times (-\infty,\ell+b]$. In particular $u_{0,1}(T-\delta,\nu,\ell) = -k/\widehat p$. 
Hence the boundary condition for $u_{1,2}$ is continuous. Recursively all functions $u_{1,n}$ (resp. $u_{0,n}$) are well-defined and bounded on $[0,T-\delta] \times (0,\infty) \times [\ell,\infty)$ (resp.  $[0,T-\delta] \times (0,\infty) \times (-\infty,\ell+b]$), continuous on $[0,T-\delta] \times (0,\infty) \times [\ell,\infty)$ (resp.  $[0,T-\delta) \times (0,\infty) \times (-\infty,\ell+b]$) and of class $C^{1,2}$ on $[0,T-\delta] \times (0,\infty) \times (\ell,\infty)$ (resp.  $[0,T-\delta] \times (0,\infty) \times (-\infty,\ell+b)$). 

Recall from \eqref{e:vol} that
$$ \vol_t = (\widehat p -1)I_t  \Vol_t.$$
To conclude we give
\begin{lemma}[Verification] \label{lem:verif_case_4}
	For $n\geq 1$, the representation $V_{i,n}(q,\nu,s) = |q|^{\widehat p} u_{i,n}(0,\nu,s)$, $i\in \{0,1\}$ holds, 
	and the optimal strategy is given by 
	\[
		Q^*_t  =  Q_0 \exp\left( -\int_{0}^{t} I_s^{(4,n)} \Vol_s |\Ymin_s|^{p-1}\sgn(\Ymin_s) ds \right), t \in [0,T],
\]
where
	\[
		\Ymin_t = u_{I_t, (n - {\mathcal N}_t )^+}(t,\nu_t,\bar{S}_t),~~~~ I_t^{(4,n)} = I^{(4)}_t 1_{(-\infty, n]} ({\mathcal N}_t + 1 )
	\]
\end{lemma}
Note that $(n -{\mathcal N}_t)^+$ is the number of remaining trading
intervals at time $t$ and $1_{(-\infty, n]} ({\mathcal N}_t + 1 )$ indicates whether $n$ trading intervals have been
used by time $t$.
\begin{proof}
The definition of $\vol$ implies 
$ \vol (s,\nu_s,\bar S_s) = 0$
when $s \not\in \mathcal I$. Hence the optimal control $(Q^*)'$ is equal to zero when $I = 0$. 
For $n=0$, there is no trading, $\varpi_0 = 0$, $Q_t = Q_0$ and 
	$$V_{0,0}(Q_0,0,\nu,s)  = -(k/\widehat p) |Q_0|^{\widehat p} = |Q_0|^{\widehat p} u_{0,0}(0,\nu,s).$$ 

For $n=1$, there are two cases. For $\tau_\ell > 0$ ($I_0 = 1$) and $\varpi_1 = \tau_\ell$: the control problem starts from the trading
state and the result follows from Section \ref{ss:pdexi2}: $V_{1,1}(q,t,\nu ,s) = u_{1,1}(t,\nu,s) |q|^{\widehat p}$. Moreover for any $n \geq 2$, $V_{1,n}(q,t,\nu,s) \leq u_{1,1}(t,\nu,s) |q|^{\widehat p}$: indeed after $\varpi_1$, there could be some opportunity for trading. 
	The second case is when the trader starts from a waiting state ($I_0 = 0$):
the trader has to wait until $\tau_{b,0}$ before trading and $\varpi_1 = \tau_{\ell,1}$. Hence for any strategy, 
$$Q_t = Q_0,\quad t \in [\![0,\tau_{b,0}]\!].$$
Then due to the Markovian structure, for $t \in [\![\tau_{b,0},\varpi_1]\!]$
	$$Q_t = Q_0\exp\left( -\int_{\tau_{b,0}}^{t\wedge \varpi_1} (p-1)\vol (s,\nu_s,\bar S_s)  |u_{1,1}(s,\nu_s,\bar S_s)|^{p-1}
	\sgn( u_{1,1}(s,\nu_s,\bar S_s)) ds \right),$$
 is the optimal state process on the interval $ [\![\tau_{b,0},T]\!]$
 and 
\begin{align*}
	V_{1,1} (Q_0,\tau_{b,0},\nu_{\tau_{b,0}},\bar S_{\tau_{b,0}}) & = |Q_0|^{\widehat p} u_{1,1}(\tau_{b,0},\nu_{\tau_{b,0}},\bar S_{\tau_{b,0}}) \\
	& = |Q_0|^{\widehat p} u_{1,1}(\tau_{b,0},\nu_{\tau_{b,0}},\ell+b) \\
	& =  |Q_0|^{\widehat p} u_{0,1}(\tau_{b,0},\nu_{\tau_{b,0}},\ell+b),
\end{align*}
where the last equality comes from the boundary condition \eqref{e:bcu0n} of $u_{0,n}$. The value function at time 0 is thus given by:
$$V_{0,1} (Q_0,0,\nu,s) =\mathbb E V_{1,1} (Q_0,\tau_{b,0},\nu_{\tau_{b,0}},\bar S_{\tau_{b,0}}) = |Q_0|^{\widehat p}  \mathbb E u_{0,1}(\tau_{b,0},\nu_{\tau_{b,0}},\ell+b).$$ 
The PDE that $u_{0,1}$ satisfies and It\^o's formula give
	$$\mathbb E u_{0,1}(\tau_{b,0},\nu_{\tau_{b,0}},\ell+b) = \mathbb E u_{0,1}(\tau_{b,0},\nu_{\tau_{b,0}},\bar S_{\tau_{b,0}})  = u_{0,1}(0,\nu,s).$$
Thus we get for any $n \geq 2$
	$$V_{0,1}(Q_0,0,\nu,s) =  |Q_0|^{\widehat p} u_{0,1}(0,\nu,s) ,$$
which achieves the proof for the case $n=1$. 

The rest of the proof proceeds by induction on $n$. Let us detail the case $n=2$ when $I_0=1$. The trader starts by following the strategy $u_{1,2}$. It\^o's formula gives:
\begin{align*}
	& |Q_{\tau_{\ell}}|^{\widehat p}  u_{1,2}(\tau_{\ell},\nu_{\tau_{\ell}},\bar S_{\tau_{\ell}})  = |Q_0|^{\widehat p}  
	u_{1,2}(0,\nu,s) + \int_0^{\tau_{\ell}} \widehat p Q'_r  |Q_r|^{\widehat p -1} \sgn(Q_r) u_{1,2}(r,\nu_r,\bar S_r) dr \\
	&\qquad +  \int_0^{\tau_{\ell}} |Q_r|^{\widehat p}  \left( \dfrac{\partial}{\partial t} + \mathcal L \right)(u_{1,2}(r,\nu_r,\bar S_r) dr  +   
	\int_0^{\tau_{\ell}} |Q_r|^{\widehat p}  \nabla u_{1,2}(r,\nu_r,\bar S_r) dW_r \\
&\quad  = |Q_0|^{\widehat p}  u_{1,2}(0,\nu,s) +   \int_0^{\tau_{\ell}} |Q_r|^{\widehat p}  \nabla u_{1,2}(r,\nu_r,\bar S_r) dW_r \\
	&\qquad +  \int_0^{\tau_{\ell}} \left[ \widehat p  Q'_r  |Q_r|^{\widehat p -1} \sgn(Q_r)  u_{1,2}(r,\nu_r,\bar S_r) + |Q_r|^{\widehat p}  \vol_r |u_{1,2}(r,\nu_r,\bar S_r)|^p \right] dr .
\end{align*}
If $Q'_r= - (p-1)\vol_r Q_r |u_{1,2}(r,\nu_r,\bar S_r)|^{p-1} \sgn(u_{1,2}(r,\nu_r,\bar S_r) )   $, we obtain that 
\begin{align*}
	& |Q_{\tau_{\ell}}|^{\widehat p}  u_{1,2}(\tau_{\ell},\nu_{\tau_{\ell}},\bar S_{\tau_{\ell}})  = |Q_0|^{\widehat p}  u_{1,2}(0,\nu,s) +   
	\int_0^{\tau_{\ell}} |Q_r|^{\widehat p}  \nabla u_{1,2}(r,\nu_r,\bar S_r) dW_r \\
&\qquad -  \int_0^{\tau_{\ell}} \dfrac{(\widehat p - 1)^{\widehat p -1}}{(\vol_r)^{\widehat p-1}} |Q'_r|^{\widehat p} dr \\
& = |Q_0|^{\widehat p}  u_{1,2}(0,\nu,s)  -  \int_0^{\tau_{\ell}} \dfrac{1}{\Vol_r^{\widehat p-1}} |Q'_r|^{\widehat p} dr \\
& \qquad +   \int_0^{\tau_{\ell}} |Q_r|^{\widehat p} \sqrt{\nu_r} \left[ \partial_s u_{1,2}(r,\nu_r,\bar S_r)  dW^{(1)}_r + c \partial_\nu u_{1,2}(r,\nu_r,\bar S_r)  dW^{(2)}_r\right] 
\end{align*}
from the definition of $\vol$. 

If $\tau_\ell \geq T$, the problem stops and since the used strategy has the same dynamics as $u_{1,1}$ with the same terminal condition at time $T$, from the case $n=1$, we know that it is the best strategy on $[0,T]$. 

If $\tau_\ell < T$, then the boundary condition connecting $u_{1,2}$ and $u_{0,1}$ over  $\bar S=\ell$, the PDE satisfied by $u_{0,1}$ and It\^o's formula give
\begin{align*}
& u_{1,2}(\tau_{\ell},\nu_{\tau_{\ell}},\bar S_{\tau_{\ell}})  =  u_{0,1}(\tau_{\ell},\nu_{\tau_{\ell}},\bar S_{\tau_{\ell}})\\
& =  u_{0,1}(\tau_{b,0},\nu_{\tau_{b,0}},\bar S_{\tau_{b,0}}) - \int_{\tau_{\ell}}^{\tau_{b,0}} \sqrt{\nu_r} \left[ \partial_s u_{0,1}(r,\nu_r,\bar S_r)  dW^{(1)}_r + c \partial_\nu u_{0,1}(r,\nu_r,\bar S_r)  dW^{(2)}_r\right] .
\end{align*}
And since there is no trading between the times $\tau_\ell$ and $\tau_{b,0}$, $Q_t = Q_{\tau_\ell}$ for any $t \in  [\![\tau_\ell,\tau_{b,0}]\!]$. The boundary conditions and $Q'_t = 0$ for $t \in [\![ \tau_l, \tau_{b,0}]\!]$ imply 
\begin{align*}
	& |Q_{\tau_{b,0}}|^{\widehat p}  u_{1,1}(\tau_{b,0},\nu_{\tau_{b,0}},\bar S_{\tau_{b,0}}) 
	 = |Q_{\tau_{\ell}}|^{\widehat p}  u_{0,1}(\tau_{b,0},\nu_{\tau_{b,0}},\bar S_{\tau_{b,0}}) \\
	& = |Q_{\tau_{\ell}}|^{\widehat p}  u_{1,2}(\tau_{\ell},\nu_{\tau_{\ell}},\bar S_{\tau_{\ell}}) \\
	& \qquad +  \int_{\tau_{\ell}}^{\tau_{b,0}} |Q_{\tau_{\ell}}|^{\widehat p} \sqrt{\nu_r} \left[ \partial_s u_{0,1}(r,\nu_r,\bar S_r)  dW^{(1)}_r + c \partial_\nu u_{0,1}(r,\nu_r,\bar S_r)  dW^{(2)}_r\right].
\end{align*}
Again the strategy is optimal and there are again two cases. If $\tau_{b,0} \geq T$, the trading is finished and the agent has traded only one time. If $\tau_{b,0} < T$, then the agent starts again to trade until $\tau_{\ell,1}\wedge T =\varpi_2$. We use the step $n=1$, starting at time $\tau_{b,0}$ from the value $Q_{\tau_{b,0}}$. The best strategy is to follow $u_{1,1}$ and the value function of the trader is given by:
	$$ V_{1,1}(Q_{\tau_{b,0}},\tau_{b,0},\nu_{\tau_{b,0}},\bar S_{\tau_{b,0}}) = 
	|Q_{\tau_{b,0}} |^{\widehat p} u_{1,1}(\tau_{b,0},\nu_{\tau_{b,0}},\bar S_{\tau_{b,0}}).$$ 
Gathering all steps together leads to 
\begin{align*}
	 |Q_0|^{\widehat p}  u_{1,2}(0,\nu,s)  & =  |Q_{\tau_{\ell}}|^{\widehat p}  u_{1,2}(\tau_{\ell},\nu_{\tau_{\ell}},\bar S_{\tau_{\ell}}) +  \int_0^{\tau_{\ell}} \dfrac{1}{\Vol_r^{\widehat p-1}} |Q'_r|^{\widehat p} dr \\
& \qquad +   \int_0^{\tau_{\ell}} |Q_r|^{\widehat p} \sqrt{\nu_r} \left[ \partial_s u_{1,2}(r,\nu_r,\bar S_r)  dW^{(1)}_r + c \partial_\nu u_{1,2}(r,\nu_r,\bar S_r)  dW^{(2)}_r\right] \\
& =\left[ |Q_{\tau_{\ell}}|^{\widehat p}  u_{1,2}(\tau_{\ell},\nu_{\tau_{\ell}},\bar S_{\tau_{\ell}}) +  \int_0^{\tau_{\ell}} \dfrac{1}{\Vol_r^{\widehat p-1}} |Q'_r|^{\widehat p} dr \right] \mathbf 1_{\tau_\ell \geq T} \\
& +\left[  \int_0^{\tau_{\ell}} \dfrac{1}{\Vol_r^{\widehat p-1}} |Q'_r|^{\widehat p} dr \right] \mathbf 1_{\tau_\ell < T}
+  |Q_{\tau_{b,0}}|^{\widehat p}  u_{0,1}(\tau_{b,0},\nu_{\tau_{b,0}},\bar S_{\tau_{b,0}}) \mathbf 1_{\tau_\ell < T ,\tau_{b,0} \geq T} \\
& + V_{1,1}(Q_{\tau_{b,0}},\tau_{b,0},\nu_{\tau_{b,0}},\bar S_{\tau_{b,0}})  \mathbf 1_{\tau_\ell \leq \tau_{b,0} < T} \\
& \quad +   \int_0^{\tau_{\ell}\wedge T} |Q_r|^{\widehat p} \sqrt{\nu_r} \left[ \partial_s u_{1,2}(r,\nu_r,\bar S_r)  dW^{(1)}_r + c \partial_\nu u_{1,2}(r,\nu_r,\bar S_r)  dW^{(2)}_r\right]  \\
& \quad -   \int_{\tau_{\ell}\wedge T}^{\tau_{b,0}\wedge T} |Q_{\tau_{\ell}}|^{\widehat p} \sqrt{\nu_r} \left[ \partial_s u_{0,1}(r,\nu_r,\bar S_r)  dW^{(1)}_r + c \partial_\nu u_{0,1}(r,\nu_r,\bar S_r)  dW^{(2)}_r\right].
\end{align*}
Since we cover all possible scenarios and since the strategies are optimal on each (random) time intervals, taking the expectation, we conclude that 
$$V_{1,2} (0,Q_0,\nu,s) = |Q_0|^{\widehat p} u_{1,2}(0,\nu,s),$$
and that the optimal strategy is $u_{1,2} \to u_{0,1} \to u_{1,1}$,
which achieves the proof of this particular case. 

The proof for the other scenarios uses very similar arguments. 
\end{proof}

\begin{remark}
Our sequence of solutions $u_{1,n}$ and $u_{0,n}$ satisfies for any $b>0$ and $n$
$$\forall s < \ell, \ u_{0,n}(T-\delta,\nu,s) = - k /\widehat p,  \qquad \forall t > T-\delta, \ u_{1,n}(t,\nu,\ell) = -k/\widehat p.$$
The boundary condition imposed on $v$ in \eqref{e:pde3} is coherent with these previous properties. Moreover $u_{0,n}$ (resp. $u_{1,n}$) solves the PDE \eqref{e:pde1} on  $[0,T-\delta) \times (0,+\infty) \times (-\infty,\ell)$ (resp. \eqref{e:pde2} on $[0,T-\delta) \times (0,+\infty) \times (\ell,+\infty)$). In other words they solve \eqref{e:pde3} on their own domain of definition. 

A natural but non-trivial question is: do $u_{1,n}$ and $u_{0,n}$ provide some approximation of $v$ given by Lemma \ref{lem:PDE_second_case}, when $b$ tends to zero and $n$ to $\infty$ ? This question is left for further research. 
\end{remark}

\subsubsection*{When the numbers of active intervals tends to $\infty$}

With Lemma \ref{lem:approx_val_fct} we saw that the value function of the control problem with a finite number of trading intervals
converge to the value function of the same control problem with no bounds on the number of trading intervals.
We now argue that 1) a similar result holds for the $u_n$ 
2) the limit determines the value function of the control problem with unbounded number of trading intervals
and 3) the limit is a solution of the coupled PDE system \eqref{e:pde1intro} and
\eqref{e:pde2intro}. Since the arguments follow the same lines as those above we only provide an outline.

The argument of the previous verification lemma gives
\[V_{i,n}(t,q,\nu,s) = |q|^{\widehat p} u_{i,n}(t,\nu,s).\]
Together with Lemma \ref{lem:approx_val_fct}, this shows that $u_{i,n}$ is a non-increasing sequence. 
Hence it converges pointwise to some function $u_i$ defined on $[0,T-\delta]\times (0,\infty) \times (-\infty,\ell+b]$
for $i=0$ and 
$[0,T-\delta]\times (0,\infty) \times [\ell,\infty)$ for $i=1.$
Moreover we have the relation:
\[
	V_i(t,q,\nu,s) = |q|^{\widehat p} u_i(t,\nu,s),i \in \{0,1\}.
\]

Since $u_i$ is bounded from below and bounded from above by $u_{i,1}$ and is upper semi-continuous, 
$u_0$ (resp. $u_1$) is a viscosity subsolution of \eqref{e:pde1} (resp. \eqref{e:pde2}) with terminal condition 
$u_0(\cdot,\cdot,T-\delta)=-k/\widehat p$ (resp. $u_1(\cdot,\cdot,T-\delta)=u^\infty(\cdot,\cdot,T-\delta)$). 
Since the terminal conditions are bounded and continuous, the arguments of Lemma \ref{lem:regular_sol} imply that $u_0$ and $u_1$ 
are smooth solution of the PDE \eqref{e:pde1} and \eqref{e:pde2} on the set $[0,T-\delta]\times (0,\infty) \times (-\infty,\ell+b)$ and $[0,T-\delta]\times (0,\infty) \times (\ell,\infty)$. 

Let us consider the lateral boundary conditions. For $s=\ell+b$, $u_0$ is equal to $\psi=u_1(\cdot,\cdot,\ell+b)$, 
which is a smooth function. From \cite{lady:solo:ural:68}, the PDE \eqref{e:pde1} with terminal condition $-k/\widehat p$ at time $T-\delta$ and lateral condition $\psi$ has a unique smooth solution. 
A similar argument shows that $u_1$ is the unique smooth solution of \eqref{e:pde2} with terminal condition $u^\infty$ at time $T-\delta$ and lateral 
condition $\phi =u_0(\cdot,\cdot,\ell) $.
These results give us the required regularity to proceed as in Lemma \ref{lem:verif_case_4} and to prove that for $t \in [0,T]$
\[Q_t = Q_0 \exp\left( -\int_0^{t} (p-1)  \mathrm{vol} (s,\nu_s,\bar S_s)  |u_{I_t}(s,\nu_s,\bar S_s) |^{p-1}\sgn(u_{I_t}(s,\nu_s,\bar S_s) )ds \right)\]
is an optimal control for the control problem \eqref{e:oc1m2p} with $I = I^{(4)}$ and ${\bm S}={\bm S}^{(4)}$.

\section{A partial analysis of the output of the algorithm and numerical examples}\label{s:numerical}

Recall $A$ of \eqref{d:Sp}, which is the percentage deviation from the target price $S_0$ of the average price at which the position is (partially) closed in the time interval $[0,T].$ 
The actual output of the trading algorithm defined by the optimal control $Q^*$
of \eqref{e:optimalcontrol1} or \eqref{e:optim_control_PDE}, 
is the random pair $(Q^*_T,A)$, where $A$ is computed for $Q=Q^*$. 
An important question is the distribution of this pair and the dependence of this
distribution on model parameters. In this section we would like to give a partial analysis of this question including some numerical examples.
Compared to the
original IS order, the modified IS order considered in the present work
has two additional parameters: the process $I$ that determines when
trading takes place and the event ${\bm S}$ that determines when full
liquidation takes place. 
In the numerical examples we will limit ourselves to 
$I = I^{(1)} =1$, ${\bm S} = {\bm S}^{(1)} = \{\bar{S}_T \geq \ell \}$. 
To further simplify the presentation and the calculations we take $p=2$, $\bar{S} =  \sigma W$, where
$W$ is a standard Brownian motion, $\sigma > 0$ a constant and $\Vol_t= V > 0$; these are
also the choices made for these parameters in the standard Almgren-Chriss framework \cite[Chapter 3]{gueant2016financial}.
Under these assumptions we will compute $Q^*$ by discretizing and numerically solving the corresponding PDE,
which becomes:
\begin{equation}\label{e:pdesimple}
	u_t  + \frac{1}{2}\sigma^2 u_{xx} - \frac{V}{\eta}I_t u^{2} = 0, 
\end{equation}
where the domain of the equation and its boundary conditions depend on
${\bm S}$ and $I$.

To better understand how $Q^*$ and $(Q^*_T/q_0,A)$ change with the model parameters,
we factor out as many parameters as possible from the calculations.
If we let 
\begin{equation}\label{e:rescale}
	v(t,x) =  \frac{V}{\eta} u(t,\sigma x) 
\end{equation}
the equation \eqref{e:pdesimple} reduces to
\begin{equation}\label{e:pdesimple2}
	v_t  + \frac{1}{2}v_{xx} - I_t v^{2} = 0.
\end{equation}

To see how $A$ depends on model parameters let us reduce the expression \eqref{d:Sp} as much as possible
(remember that $Q_0 = q_0$ is the initial position size):
\begin{align*}
	A &= \frac{X_T - (Q_0-Q_T) S_0}{(Q_0-Q_T) S_0}
	\intertext{By \eqref{e:XTp} (the expression for $X_T$) and the assumption $L(v) = \eta v^2$:}
	&= \frac{- \int_0^T S_t Q'_t dt  -\frac{\eta}{V} \int_0^T  (Q'_t)^2 dt  - S_0(Q_0 - Q_T)}{S_0(Q_0 - Q_T)}
	\intertext{By the definition \eqref{e:midprice} of $S_t$ and the assumptions $\kappa(v) = kv$, $\bar{S}_t = \sigma W_t$:}
	& =  \frac{ -k/2 ( Q_T - Q_0)^2 - \sigma\int_0^T W_t Q'_t dt  -\frac{\eta}{V}\int_0^T  (Q'_t)^2 dt}{(Q_0 - Q_T)S_0}.
\end{align*}
Simplifying the last expression we get
\begin{equation}\label{simpleSp}
	A= 
	-\frac{k Q_0}{2S_0 } \left( 1- \frac{Q_T}{Q_0} \right) - 
	\frac{\sigma}{S_0}\frac{1}{(Q_0-Q_T)}\int_0^T W_t Q'_t dt -
	\frac{\eta}{VS_0} \frac{1}{(Q_0- Q_T)} 
	\int_0^T (Q'_t)^2 dt.
\end{equation}
From this expression we see that $A$ consists of three components: 1) one due to the permanent price impact
2) one due to random fluctuations in price and 3) one due to transaction costs.
All components consist of a coefficient term and a term depending on $Q$ or its
derivative $Q'$:
\begin{align*}
	&\text{Permanent market impact term: }A_1 = 1- \frac{Q_T}{Q_0},
	& 
	&\hspace{-1cm}\text{coefficient: } -\frac{k Q_0}{2S_0 }, \\
	&\text{Random fluctuations term: } 
	 A_2= \frac{1}{(Q_0-Q_T)}\int_0^T W_t Q'_t dt,  &
	&\hspace{-1cm}\text{coefficient: }
	\frac{\sigma}{S_0},\\
	&\text{Transaction costs term: }
	A_3=
	\frac{1}{1-\fq_T}
	\int_0^T (Q'_t/Q_0)^2 dt,&
	&\hspace{-1cm}\text{coefficient: }
 \frac{\eta Q_0}{VS_0},
\end{align*}
where
\[
\fq_t = Q_t/Q_0.
\]
The permanent impact term
$1 - Q_T/Q_0=1-\fq_T$ is the portion of the initial position that is closed;
$A$ depends linearly on
this portion with coefficient $\frac{kQ_0}{2S_0}$.
Secondly note that
if $S_0$, $k$ and $\eta$ are parameterized as multiples of $\sigma$
then none of the coefficients appearing in
$A$ depend on $\sigma.$ 
We will comment on the behavior of the other two
terms below.

Before we move on let us note the following for comparison. The case $I=1$ and ${\bm S} = \Omega$
corresponds to the standard Almgren Chriss liquidation algorithm for which the optimal 
control is known to be 
\begin{equation}\label{e:qstars}
	Q^{*,S}_t = q_0 \frac{T-t}{T},
\end{equation}
i.e., closing the position with uniform speed over the time interval $[0,T].$
Then $(Q')^{*,S}_t/q_0 = -1/T$ and $\fq_T = 0$.
These reduce $A_3$ to
\begin{equation*}
A_{3}^S = \dfrac{1}{T}
\end{equation*}
for the standard IS algorithm.
Similarly, for $A_2$ we have 
\begin{equation}\label{e:A2s}
	A_2^{S}= -\dfrac{1}{T}\int_0^T W_t dt,
\end{equation}
which is normally distributed with $0$ mean by the iid normal increments of $W$.

\bigskip
We continue our analysis with the choices
$I = I^{(1)} = 1$ and ${\bm S} =  {\bm S}^{(1)} = \{ \bar{S}_T \geq \ell \}$
for $I$ and ${\bm S}$ given in \eqref{e:I1};
these choices correspond to: no restriction on trading
and closing the position fully is required only when the terminal
price $\bar{S}_T$ is above a given threshold $\ell$.
Parallel to the change of variable in \eqref{e:rescale} we assume
$\ell$ is given as a multiple of $\sigma >0$; with this convention
and the assumption $\bar{S}_t = \sigma W_t$,
${\bm S}$ becomes ${\bm S} = \{ W_T \geq \ell \}.$ For $I_t = 1$, the
PDE \eqref{e:pdesimple2} is
\begin{equation}\label{e:pdesimplec1}
	v_t  + \frac{1}{2}v_{xx} - v^{2} = 0;
\end{equation}
for ${\bm S} = \{ W_T \geq \ell \}$ the domain and the boundary conditions
for this PDE
are: $(t,x) \in [0,T] \times {\mathbb R}$ and 
\begin{equation}\label{e:bc1}
	v(T,x) = \infty \cdot {\bm 1}_{[\ell,\infty)}(x) 
	-\frac{kV}{2 \eta} \cdot{\bm 1}_{(-\infty,\ell)}(x),
\end{equation}
$x \in {\mathbb R}$,
where we again use the scaling \eqref{e:rescale}.

Recall our convention that $k$ and $\eta$ are specified as multiples
of $\sigma$; it follows that PDE \eqref{e:pdesimplec1} and
its boundary condition \eqref{e:bc1} are independent of $\sigma.$
The optimal control $Q^*$ is computed 
from $v$ via the formula \eqref{e:optim_control_PDE}
\begin{equation}\label{e:qstarv}
	Q^*_t = q_0 \fq_t = q_0 \exp\left(-\int_0^T v(t,W_t) dt\right)
\end{equation} we note that $Q^*$ is
independent of $\sigma.$
We have already noted that the coefficients in \eqref{e:rescale}
are independent of $\sigma$. We have observed above that the same
is true also for $Q^*$, therefore all of $A_1$, $A_2$ and $A_3$
are independent of $\sigma$ as well. %
This gives us the following result:
\begin{lemma}
	Suppose all of $S_0$, $k$, $\eta$ and $\ell$ are parameterized
	as multiples of $\sigma.$ Then $Q^*$ and $A$ do not depend on
	$\sigma$. 
\end{lemma}
Note that the same analysis in fact holds for all of  $I = I^{(i)}$, ${\bm{S}} = {\bm{S}}^{(i)}$, $i=2,3,4$ treated in the previous sections.

As already noted $Q^*$ is computed via the solution of the PDE
\eqref{e:pdesimplec1} which obviously doesn't have an explicit solution.
To see how $Q^*$ behaves we will solve \eqref{e:pdesimplec1}
numerically; for the parameter values we begin by
considering those used
in \cite[Chapter 3]{gueant2016financial}: $T=1$, $\eta =0.1$,
$V = 4 \times 10^6$, $S_0=45$, $\sigma =0.6$.  
Recall that the permanent price impact parameter $k$ doesn't appear
in the control problem corresponding to the original IS order,
so no value for $k$ is specified in \cite[Chapter 3]{gueant2016financial}.
A $k$ value of $k=2 \times 10^{-7}$ accompanying these parameter values is given in
\cite[Chapter 8]{gueant2016financial} in the context of block
trade pricing.
The assumption \eqref{a:vol} in the present case reduces to
\[
k V/2\eta < 1;
\]
for the above parameter values we have $kV/2\eta = 4$, therefore
the above parameter values do not satisfy Assumption \ref{a:vol}.
To continue with our numerical example, we take
$\eta = 0.3$, $V = 4 \times 10^6$ and $k= 10^{-7}$ for these values
we have $kV/2\eta = 2/3$ which satisfies \eqref{a:vol}.
In addition to these, we need to provide a value for the $\ell$ parameter,
which we choose as $\ell = 1.4\times \sigma $.
The graph of $u$ for the parameter values above are shown in Figures \ref{f:pdesolI1.png} and \ref{f:pdesolI1f2.png}.

\ninseps{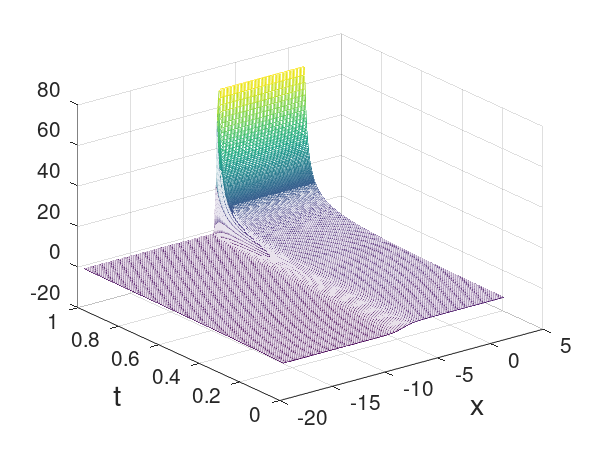}{Graph of $u$ for $I=1$ and ${\bm S} = \{ W_T \geq \ell \}$}
{0.5}
\ninseps{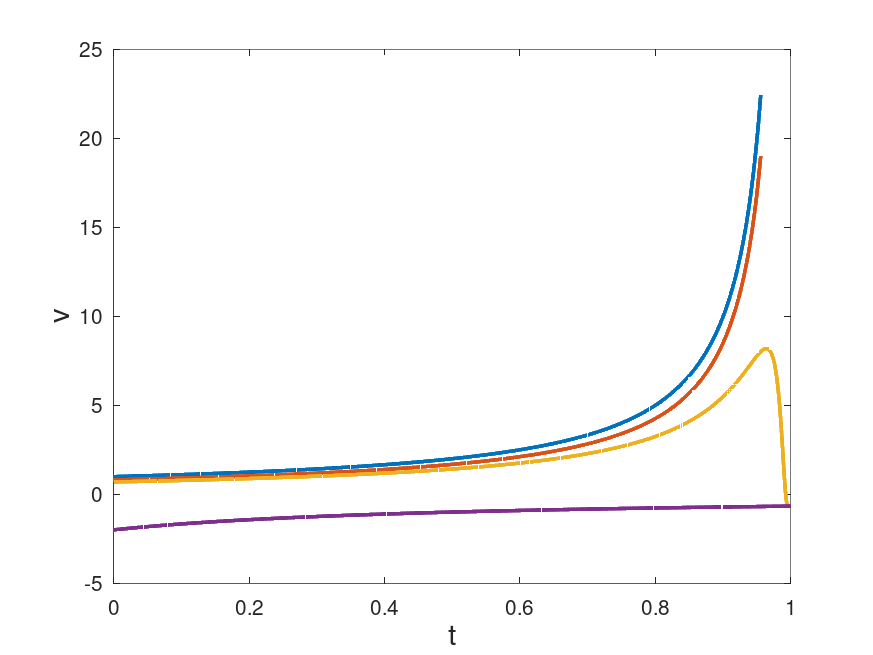}{Graph of $u(x,\cdot)$ for $x \in \{0, \ell,1.2\ell,5\ell\}$}{0.3}
We note that for $x >\ell$ and $x$ away from $\ell$, $u(x,\cdot)$ behaves like $t \mapsto y_t=1/(T-t)$ (the solution of
\eqref{e:pdesimple} with terminal condition $y_T = \infty$).
and for $x < \ell$ and $x$ away from $\ell$, $u(x,\cdot)$ behaves like $t \mapsto z_t$ (given by \eqref{e:z1})
The negative boundary condition for $u$ means that $u(x,t)$ takes negative values  for $x < \ell$;
\eqref{e:XTp} implies that whenever $u$ is negative, the corresponding $Q^*$ is actually buying the underlying stock.

\begin{figure}[H]
\begin{center}
	\scalebox{0.6}{ \input{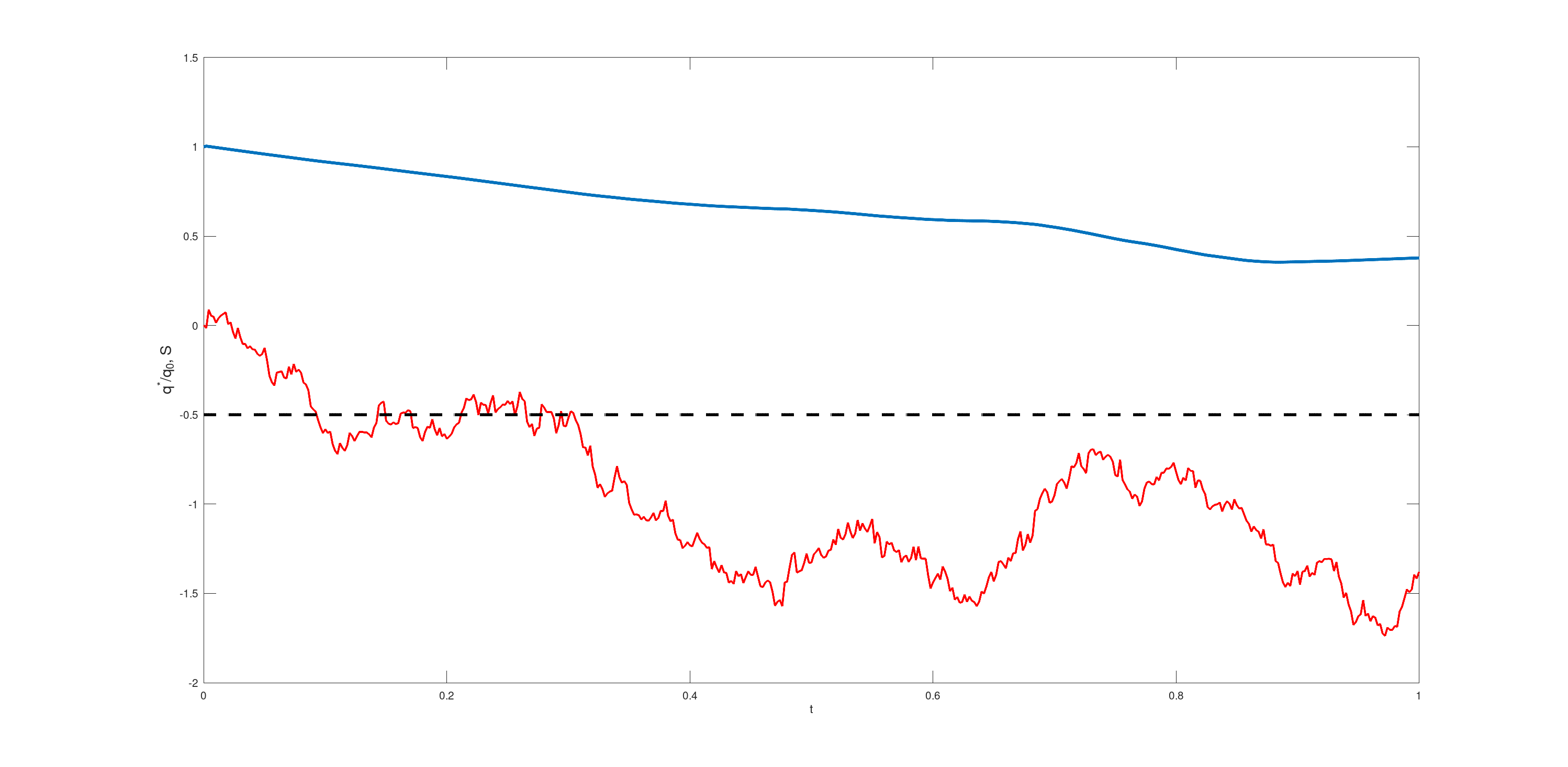}}
	\scalebox{0.6}{ \input{path4I1}}
	\scalebox{0.6}{ \input{path2I1}}
	\scalebox{0.6}{ \input{path3I1}}
	\caption{Sample paths of $\bar{S}$ and $Q^*$ for $I = 1$ and ${\bm S} = \{\bar{S}_T \geq \ell \}$; the dashed
	line shows $\ell =-1.4\sigma$}
	\label{f:samplepaths1}
\end{center}
\end{figure}

Figure \ref{f:samplepaths1} shows four sample paths of $\bar{S}$ and $Q^*$. 
In the first two examples $\bar{S}$ stays above $\ell$ at all times and the corresponding $Q^*$ goes parallel to $Q^{*,S}$ of \eqref{e:qstars},
the optimal liquidation path for the standard IS order.
In the third example $\bar{S}$ is below $\ell$ approximately in the time interval
$[0.6,0.8]$ when trading slows down, it goes above $\ell$ around $0.8$ and closes
above $\ell$;
correspondingly $Q^*$ speeds up trading after $0.8$ and closes the position
at terminal time.
In the fourth example, $\bar{S}$ hits $\ell$ around the middle of the trading interval and remains below $\ell$ till the end;
correspondingly $Q^*$ slows down and stops trading and the position is only partially closed at terminal time. 
In the last example $Q^*$ is in fact slightly increasing near $t=T=1$
(i.e., $Q^*$ buying the underlying asset) ; this is due to the negative value that the
terminal value takes for $x < \ell.$
These examples suggest that $Q^*$ behaves approximately as follows:
when $\bar{S}$ is above $\ell$, it behaves like the standard IS algorithm $Q^{*,S}$, linearly closing the remaining position; when $\bar{S}$
goes below $\ell$, $Q^*$ slows down/ stops trading. The negative boundary condition implies that the algorithm can in fact execute buy trades
especially when the price is below $\ell$ near terminal time $T$.

\subsubsection*{Distribution of $(Q^*_T/q_0, A)$} 
For ${\bm S} = \{ \bar{S}_T \geq \ell \}$, the position fully closes when the closing price is above the lower-bound $l$, therefore, the probability
that the algorithm closes the position at terminal time is:
\[
	{\mathbb P}( Q^*_T = 0) = {\mathbb P}(W_T \geq \ell ) = 1- N_{0,1}(\ell/\sqrt{T})
\]
where $N_{0,1}$ denotes the standard normal distribution.

\begin{figure}[H]
\begin{center}
	\scalebox{0.7}{ \input{Ptail}}
	\scalebox{0.7}{ \input{logP}}
\end{center}
	\caption{On the left: graph of $p_x = {\mathbb P}( Q_T^*/q_0 > x | Q_T^* > 0)$, on the right:
	graphs of $-\log(p_x)$ and $\frac{x q_0}{m_T}$, $m_T = {\mathbb E}[Q_T^* | Q_T^* > 0]$}\label{f:q0qTg0}
\end{figure}
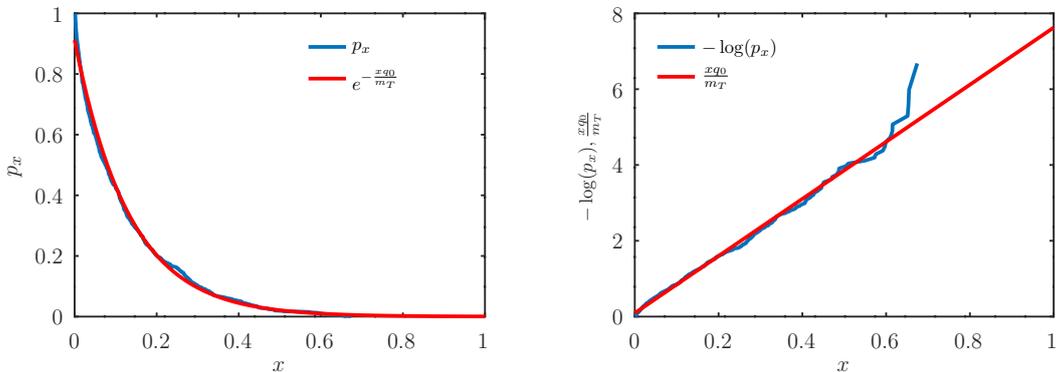

A random variable $E$ is said to be exponentially distributed with rate $\lambda$ if 
	${\mathbb P}( E > x) =e^{-\lambda x}$, i.e.,
\begin{equation}\label{e:expdist}
	-\log({\mathbb P}( E > x)) =\lambda x;
\end{equation}
a well known fact is
\begin{equation}\label{e:expofexp}
	{\mathbb E}[E] = \frac{1}{\lambda}.
\end{equation}
The distribution of $\fq_T=Q^*_T/q_0$ over $(0,\infty)$ depends on $u$ via \eqref{e:qstarv} (or \eqref{e:optim_control_PDE}) and it obviously doesn't have an explicit formula.
Figure \ref{f:q0qTg0} shows graphs
 of $x \mapsto {\mathbb P}( \fq_T > x | \fq_T > 0)$,
 $x\mapsto -\log({\mathbb P}( \fq_T > x | \fq_T > 0))$ and $x\mapsto \frac{x q_0}{{\mathbb E}[Q_T^* |Q_T^* > 0]}$
 (all estimated via simulating $10^4$ sample paths). 
These graphs, \eqref{e:expdist} and \eqref{e:expofexp} suggest that
 the exponential distribution provides a rough approximation
 for the conditional distribution of $Q_T^*/q_0$ given the
 event $\{Q_T^* > 0\}.$
 An exponentially distributed random variable satisfies
${\mathbb E}[E] =\frac{1}{\lambda}$ and $\var(E) = \frac{1}{\lambda^2}$.
In the case of $Q_T^*/q_0$ conditioned over $\{Q_T^* > 0\}$ we have
the Monte Carlo
estimates ${\mathbb E}[\fq_T  | \fq_T > 0] = 0.1218$ and 
$\var(\fq_T | \fq_T > 0)^{1/2} = 0.1387$ for the parameter values specified above.

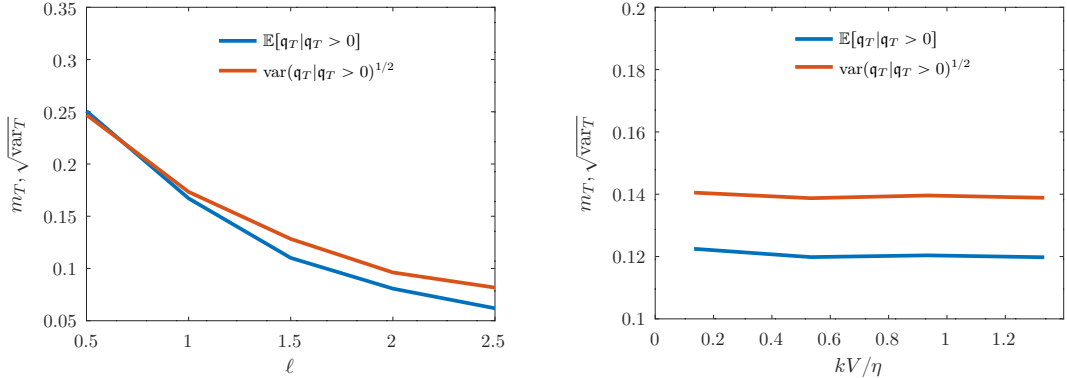
\begin{figure}[h]
\begin{center}
	\scalebox{0.7}{ \input{varlI1}}
	\scalebox{0.7}{ \input{varkI1}}
\end{center}
	\caption{Graph of $m_T = {\mathbb E}[\fq_T  | \fq_T > 0]$ and 
	$\sqrt{\var_T} = \var(\fq_T   | \fq_T > 0)^{1/2}$ as a function of $\ell$ }\label{f:varlI1}
\end{figure}

$Q^*$ depends on $\ell$ via the domain of the PDE \eqref{e:pdesimplec1}
and on $kV/2\eta$ via the terminal condition \eqref{e:bc1}. Figure
\ref{f:varlI1} shows how 
${\mathbb E}[\fq_T  | \fq_T>0]$ and 
$\var(\fq_T  | \fq_T > 0)^{1/2}$ vary with these parameters.

We have already noted that the permanent impact factor term  
$A_1$ of \eqref{simpleSp} is fully determined by $\fq_T$.
We now consider the joint distribution of 
$(A_2, \fq_T)$.
This distribution consists of two parts:
the distribution of $A_2$ conditioned on $\fq_T = 0$ (i.e., the cases where
the algorithm closes the initial position $q_0$ fully) and
the conditional distribution of $A_2$ given 
$\fq_T$ for $\fq_T> 0$
(the cases where the algorithm closes the initial position $q_0$ partially).
If $Q^*$ were a deterministic function (as in the case of the standard 
IS order), $A_2$ would be normally distributed by
the normal and independent increments of $W$ (see \eqref{e:A2s}).
The $q$-$q$ plot of the conditional distribution of $A_2$ given $Q_T^* = 0$ 
and $Q_T^* =x$ for several values of $x$ is shown Figure \ref{f:Wqqplots};
(for $x > 0$ we approximate $P(A_2\in A | Q_T^*/q_0 = x)$
with $P(A_2 \in A | Q_T^*/q_0 \in (x,x+\delta))$
where $\delta > 0$ is small and we estimate the latter by simulating
$2 \times 10^5$ sample paths of $W$ and $Q^*$). These plots suggest that the conditional distribution $A_2$
given $Q_T^*/q_0$ is approximately normal even though $Q^*$ is random and a function of $W$.

\ninseps{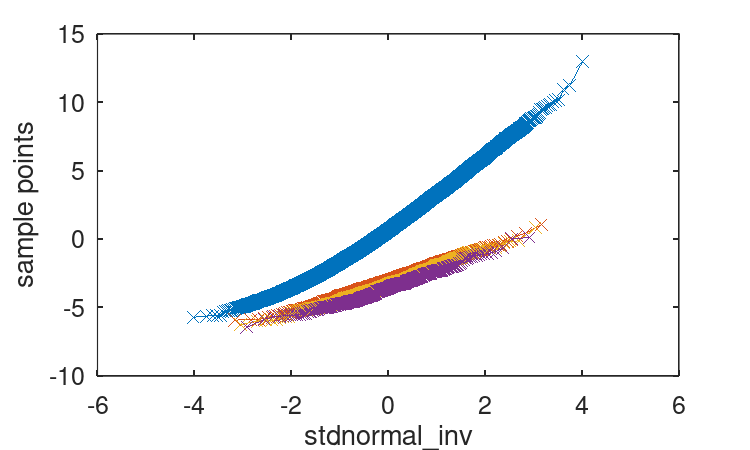}{$q$-$q$ plots of the conditional distribution
of $A_2$ given $\fq_T=x$ for $x=0$, $x=0.06 + j 0.07$, $j \in \{1,2,3,4\}$}{0.3}

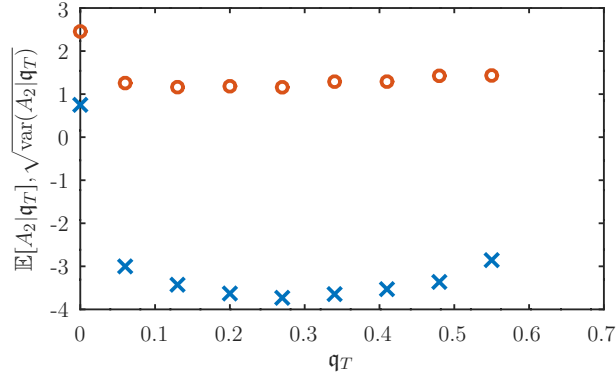
\begin{figure}[H]
\begin{center}
	\scalebox{0.7}{\input{S2primegivenqTmeanvar}}
\end{center}
	\caption{Graphs of 
	${\mathbb E}[A_2 | \fq_T=x]$ and
	$\sqrt{\var(A_2 | \fq_T=x)}$ for $x=0$ and $x = 0.06 + j 0.07$, $j \in \{1,\cdots,8\}$} \label{f:S2primegivenqTmeanvar}
\end{figure}

Figure \ref{f:S2primegivenqTmeanvar} shows the graphs of 
	${\mathbb E}[A_2 | \fq_T=x]$ and
	$\sqrt{\var( A_2| \fq_T=x)}$ (using the same approximation as above and then simulation).

Figures \ref{f:Wvaryl2} and \ref{f:Wvaryk} show how the distribution of $A_2$ change with $\ell$ and $kV/\eta.$
\begin{figure}[h]
\begin{center}
	\scalebox{0.55}{\input{Wvaryl1}}
	\scalebox{0.55}{\input{Wvaryl2}}
\end{center}
	\caption{Graphs of 
	$P(A_2 \le x) | \fq_T=0)$ 
	and
	$P(A_2 \le x) | \fq_T>0)$ 
	for $\ell=-\sigma$
	and $\ell=-1.4\sigma$}\label{f:Wvaryl2}
\end{figure}

\begin{figure}[h]
\begin{center}
	\scalebox{0.55}{\input{WvarykqT0I1}}
	\scalebox{0.55}{\input{WvarykqTg0I1}}
\end{center}
	\caption{Graphs of 
	$P(A_2 \le x) | \fq_T=0)$ 
	and
	$P(A_2 \le x) | \fq_T>0)$ 
	for $kV/\eta =4/3$
	and $kV/\eta=4/30$}\label{f:Wvaryk}
\end{figure}

\begin{figure}[H]
\begin{center}
	\scalebox{0.55}{\input{EtavarykqT0I1}}
	\scalebox{0.55}{\input{EtavarykqTg0I1}}
\end{center}
	\caption{Graphs of 
	$P(A_3 \le x) | \fq_T=0)$ 
	and
	$P(A_3 \le x) | \fq_T>0)$ 
	for $kV/\eta = 4/3$
	and $kV/\eta=4/30$}\label{f:Wvaryk1}
\end{figure}
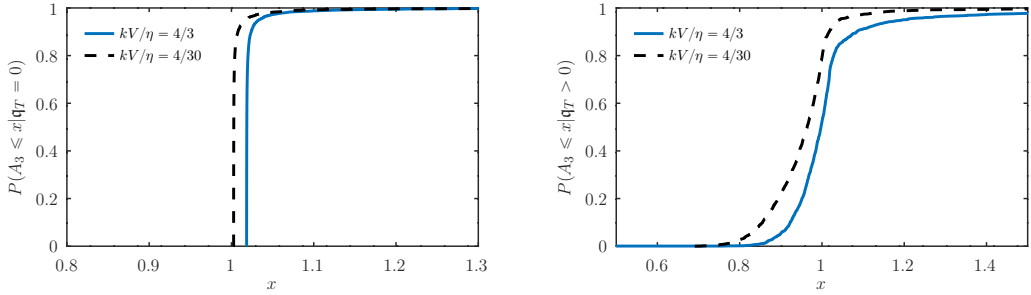

Lastly, we consider the distribution of $A_3.$ Like $A_2$ this distribution consists of two parts:
over the event $\fq_T =0$ and over the event $\fq_T > 0.$
\begin{figure}[H]
\begin{center}
	\scalebox{0.55}{\input{A3qT0dist}}
	\scalebox{0.55}{\input{A3vsqT}}
\end{center}
	\caption{On the left: distribution of $A_3$ conditioned on $\fq_T=0$,
	on the right:  graphs of ${\mathbb E}[A_3 | \fq_T = x]$~ ($\times$) and 
	$\var(A_3|\fq_T = x)^{1/2}$~($o$)}
	\label{f:A3qT0dist}
\end{figure}
The first part of Figure \ref{f:A3qT0dist} shows the distribution of $A_2$ over $\fq_t = 0$; this graphs suggests 
that $A_3$ behaves approximately like $A_3^S = 1$ of the standard IS algorithm: most of the mass is concentrated around
a constant near $1$. The second part of Figure \ref{f:A3qT0dist} shows how the conditional mean and variance of $A_3$ change with
$\fq_T$ for $\fq_T > 0.$ This graph suggests that the near constant behavior of $A_3$ persists for $\fq_T > 0.$
Finally Figures \ref{f:Wvaryl1} and \ref{f:Wvaryk1} show how the distribution of $A_3$ changes with model parameters
$\ell$ and $kV/\eta.$
\begin{figure}[H]
\begin{center}
	\scalebox{0.55}{\input{EtavarylqT0I1}}
	\scalebox{0.55}{\input{EtavarylqTg0I1}}
\end{center}
	\caption{Graphs of 
	$P(A_3 \le x) | \fq_T=0)$ 
	and
	$P(A_3 \le x) | \fq_T>0)$ 
	for $\ell=-\sigma$
	and $\ell=-1.4\sigma$}\label{f:Wvaryl1}
\end{figure}

\section{Comments on the continuity problem}
Remember that the only condition on the terminal behavior of the
minimal supersolution $\Ymin$ is \eqref{e:minimalitycond}:
\begin{equation*}
	\liminf_{t\rightarrow T} \Ymin_t \ge \xi.
\end{equation*}
For the class of stochastic optimal control problems studied in this paper,
minimal supersolutions suffice to characterize the value function; this
is the reason for the focus of this paper on minimal supersolutions.
An interesting further 
question is whether $\Ymin$ satisfies the terminal condition
with equality, i.e., whether
\begin{equation}\label{e:continuity}
	\lim_{t\rightarrow T} \Ymin_t = \xi
\end{equation}
holds.
We call this question ``the continuity problem.''
Recall the constraint $\{Q_T = 0\}$ over the set ${\bm S}$ for the
stochastic optimal control problem.
Let $Q^*$ be the optimal control defined by $\Ymin.$
The relation \eqref{e:continuity} corresponds 
to $\{Q_T^*=0\} \supset {\bm S}$, i.e., the constraint is satisfied
over ${\bm S}$, but it may happen that it is satisfied beyond the event 
${\bm S}.$ The relation \eqref{e:continuity} on the other hand implies
that the optimal control satisfies the constraint tightly: $\{Q_T^* = 0\}
= {\bm S}.$ This further information about the optimal control
is the main significance of the continuity problem
from the control perspective.

In previous works \cite{popi:16,krus:popi:seze:18,mahd:popi:seze:21, popier2020continuity} 
we studied the continuity problem for 
a number of BSDE with singular
terminal conditions. The results in these works do not apply in the current
setting because they all assume monotone generators. We think that
the continuity problem can be resolved positively for the Markovian cases covered
in Section \ref{s:pde} using the ideas and methods presented in these prior works.
The continuity problem for the general BSDE studied in Section \ref{s:bsde} is more challenging:
results in all of the previously cited works concern ${\bm S}$ of special forms such as
${\bm S} = \{ \tau \le T\}$ where $\tau$ is a stopping time. We think that results and methods in
\cite{mahd:popi:seze:21, popier2020continuity} can be extended to cover the
BSDE in Section \ref{s:bsde} for ${\bm S}$ of the form studied in these works. To the best
of our knowledge, the continuity problem
with terminal conditions of the form $\lim_{t\rightarrow T} \Ymin_t = \infty \cdot {\bm 1}_{\bm S}$
for general ${\bm S} \in {\mathcal F}_T$ is an open problem even when the filtration is generated
only by a Brownian motion.

\section{Conclusion}\label{s:conc}

The present work studies the optimal control problem \eqref{e:oc1m2} of optimal liquidation where the key parameters are the measurable set ${\bm S}$ specifying
conditions for full liquidation and the process $I$ specifying when trading is allowed; we note that choosing ${\bm S} \neq \Omega$ introduces a negative
term into the terminal condition of the associated BSDE determined by the permanent price impact parameter.
We find the optimal control to be of the form
\begin{align*}
Q^*_t & =  q_0 \exp\left( -\int_{0}^{t} I_s \Vol_s |\Ymin_s|^{p-1}\sgn(\Ymin_s) ds \right)
\end{align*}
where $\Ymin$ is the first component of the minimal supersolution of the BSDE \eqref{e:convex_BSDE} with terminal condition 
$\xi = - \dfrac{k}{\widehat p} \mathbf 1_{\bm S^c} + \infty \mathbf 1_{\bm S}.$
In Section \ref{s:bsde} this is proved directly using the BSDE for the case $p=2$; Section \ref{s:pde} uses a PDE approach and focuses on the case $p \neq 2$
where the price dynamics are assumed Markovian and $I$ and ${\bm S}$ are chosen so that the resulting problem is either Markovian
or can be broken into Markovian pieces.

In the remaining  paragraphs we comment on some of our assumptions and on future research.
In the present work we use the identity function as our utility function (see \eqref{e:oc0}); 
this corresponds to setting the risk aversion parameter $\gamma$ in the exponential utility
function to $0$. Let us comment on the case $\gamma > 0.$ Recall that we proceed in two steps: we start with the
control problem \eqref{e:oc0} and then derive from it the secondary problem \eqref{e:oc1m2} and then work with this problem.
For $\gamma > 0$ the reduction from \eqref{e:oc0} to \eqref{e:oc1m2} is possible only in special cases (for example, when the
optimal control is known to be deterministic, as in the case of the standard IS order (see \cite[Chapter 3]{gueant2016financial}).
To circumvent this problem, many prior works directly start from the secondary problem and 
add a risk measure on the size of the portfolio directly to this problem
(see, e.g., \cite{almg:12,anki:jean:krus:13}). When the risk measure is added directly to the second problem, it can even be
allowed to be a random process,
an example is given in \eqref{e:sockp} (the term $\gamma_t |Q_t|^{\widehat p}$).
As is the case in these prior works, 
it is straight-forward to introduce a penalization term on the size of the portfolio
involving a positive process in the stochastic optimal control problem \eqref{e:oc1m2} and almost all of the
analysis presented in sections \ref{s:bsde} and \ref{s:pde} will have straightforward modifications. On the other hand,
an analysis of the exponential utility case  for $\gamma > 0$ (i.e., introducing $\gamma > 0$ and starting the analysis from the problem
$\sup_{Q \in {\mathcal A}_{I,{\bm S}}} {\mathbb E}[-e^{-\gamma \tilde{X}_T}]$)  would require significant changes from the analysis presented in the
current work. In our view, one benefit of 
starting the analysis
from \eqref{e:oc0} as we did in the present work, is that it provided
a relatively simple framework to understand
the role permanent price impact plays in the stochastic optimal control problem in the presence of the parameter ${\bm S}$.

An important assumption in the present work is the monetary representation of the terminal position $Q_T$, for which we used $m(q) = q S_T$. 
A simple modification that would keep the problem within the framework of the current work is to set
$m(q) =q S_T+ \xi' |q|^{\widehat{p}}$ where $\xi'$ is 
an  ${\mathcal F}_T$-measurable random variable.
This would introduce an additional $\xi'$ term into the terminal condition of the BSDE which would not impact the analysis as long as $\xi'^{-}$ is bounded
from above.
Other choices are obviously possible depending on how the events taking place after time $T$ are modeled. These choices will give rise to different
stochastic optimal control problems whose analysis will probably require
new/other tools and ideas. This is a natural direction for future research. 

Recall the random variable $A$ of \eqref{d:Sp}, the percentage deviation from the target price $S_0$ of the average price at which
the portfolio is (partially) closed at terminal time.
We provided a numerical study of the distribution of $\fq_T = Q_T^*/q_0$ (remaining
portion of the position) and $A$ in Section
\ref{s:numerical} for the case $I=1$ and ${\bm S} = \{ \bar{S}_T > \ell\}.$
An analytical study of these distributions for this choice of $(I, {\bm S})$
and for others seems interesting and challenging.
One idea in the study of the distribution of $Q_T/q_0$ is the use of Malliavin calculus (see \cite[Theorem 2.1.3]{nual:06}). 
This would require that the product 
$\Ymin \vol$ has a Malliavin derivative with some additional conditions. The Malliavin derivative of the solution of a BSDE has already been studied; see among others \cite{elka:peng:quen:97} for integrable terminal condition and \cite{popi:06} for singular terminal condition. 
In the Markovian framework, it requires some Malliavin regularity on the forward processes $S$ and $\nu$. 
A study of these ideas and problems remain for future work.
Another natural direction is to try to compute the performance of the algorithms of the present work (especially the joint distribution of $(\fq_T,A)$) 
on real trading data.

\bibliography{biblio}

\end{document}

%% file: path1I1.tex
\setlength{\unitlength}{1pt}
\begin{picture}(0,0)
\includegraphics[scale=1]{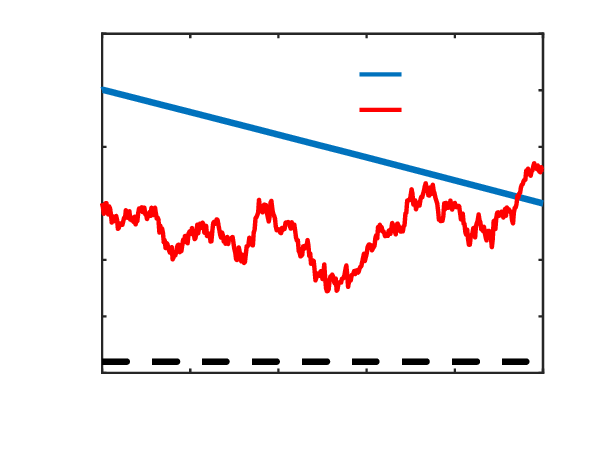}
\end{picture}%
\begin{picture}(288,216)(0,0)
\fontsize{12}{0}\selectfont\put(48.9901,28){\makebox(0,0)[t]{\textcolor[rgb]{0.15,0.15,0.15}{{0}}}}
\fontsize{12}{0}\selectfont\put(91.3201,28){\makebox(0,0)[t]{\textcolor[rgb]{0.15,0.15,0.15}{{0.2}}}}
\fontsize{12}{0}\selectfont\put(133.65,28){\makebox(0,0)[t]{\textcolor[rgb]{0.15,0.15,0.15}{{0.4}}}}
\fontsize{12}{0}\selectfont\put(175.98,28){\makebox(0,0)[t]{\textcolor[rgb]{0.15,0.15,0.15}{{0.6}}}}
\fontsize{12}{0}\selectfont\put(218.31,28){\makebox(0,0)[t]{\textcolor[rgb]{0.15,0.15,0.15}{{0.8}}}}
\fontsize{12}{0}\selectfont\put(260.64,28){\makebox(0,0)[t]{\textcolor[rgb]{0.15,0.15,0.15}{{1}}}}
\fontsize{12}{0}\selectfont\put(43,36.9895){\makebox(0,0)[r]{\textcolor[rgb]{0.15,0.15,0.15}{{-1.5}}}}
\fontsize{12}{0}\selectfont\put(43,64.1246){\makebox(0,0)[r]{\textcolor[rgb]{0.15,0.15,0.15}{{-1}}}}
\fontsize{12}{0}\selectfont\put(43,91.2597){\makebox(0,0)[r]{\textcolor[rgb]{0.15,0.15,0.15}{{-0.5}}}}
\fontsize{12}{0}\selectfont\put(43,118.395){\makebox(0,0)[r]{\textcolor[rgb]{0.15,0.15,0.15}{{0}}}}
\fontsize{12}{0}\selectfont\put(43,145.53){\makebox(0,0)[r]{\textcolor[rgb]{0.15,0.15,0.15}{{0.5}}}}
\fontsize{12}{0}\selectfont\put(43,172.665){\makebox(0,0)[r]{\textcolor[rgb]{0.15,0.15,0.15}{{1}}}}
\fontsize{12}{0}\selectfont\put(43,199.8){\makebox(0,0)[r]{\textcolor[rgb]{0.15,0.15,0.15}{{1.5}}}}
\fontsize{13}{0}\selectfont\put(17,118.395){\rotatebox{90}{\makebox(0,0)[b]{\textcolor[rgb]{0.15,0.15,0.15}{{$\frac{q^*}{q_0}, \bar{S}$}}}}}
\fontsize{13}{0}\selectfont\put(154.815,13){\makebox(0,0)[t]{\textcolor[rgb]{0.15,0.15,0.15}{{$t$}}}}
\fontsize{10}{0}\selectfont\put(196.624,180.284){\makebox(0,0)[l]{\textcolor[rgb]{0,0,0}{{$\frac{q^*}{q_0}$}}}}
\fontsize{10}{0}\selectfont\put(196.624,163.279){\makebox(0,0)[l]{\textcolor[rgb]{0,0,0}{{$\bar{S}$}}}}
\end{picture}

%% file: path4I1.tex
\setlength{\unitlength}{1pt}
\begin{picture}(0,0)
\includegraphics[scale=1]{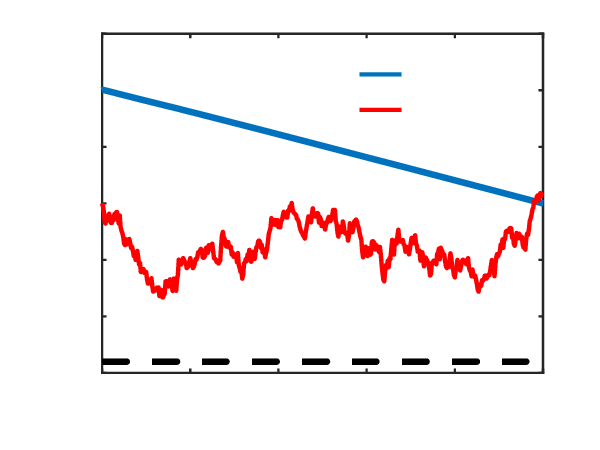}
\end{picture}%
\begin{picture}(288,216)(0,0)
\fontsize{12}{0}\selectfont\put(48.9901,28){\makebox(0,0)[t]{\textcolor[rgb]{0.15,0.15,0.15}{{0}}}}
\fontsize{12}{0}\selectfont\put(91.3201,28){\makebox(0,0)[t]{\textcolor[rgb]{0.15,0.15,0.15}{{0.2}}}}
\fontsize{12}{0}\selectfont\put(133.65,28){\makebox(0,0)[t]{\textcolor[rgb]{0.15,0.15,0.15}{{0.4}}}}
\fontsize{12}{0}\selectfont\put(175.98,28){\makebox(0,0)[t]{\textcolor[rgb]{0.15,0.15,0.15}{{0.6}}}}
\fontsize{12}{0}\selectfont\put(218.31,28){\makebox(0,0)[t]{\textcolor[rgb]{0.15,0.15,0.15}{{0.8}}}}
\fontsize{12}{0}\selectfont\put(260.64,28){\makebox(0,0)[t]{\textcolor[rgb]{0.15,0.15,0.15}{{1}}}}
\fontsize{12}{0}\selectfont\put(43,36.9895){\makebox(0,0)[r]{\textcolor[rgb]{0.15,0.15,0.15}{{-1.5}}}}
\fontsize{12}{0}\selectfont\put(43,64.1246){\makebox(0,0)[r]{\textcolor[rgb]{0.15,0.15,0.15}{{-1}}}}
\fontsize{12}{0}\selectfont\put(43,91.2597){\makebox(0,0)[r]{\textcolor[rgb]{0.15,0.15,0.15}{{-0.5}}}}
\fontsize{12}{0}\selectfont\put(43,118.395){\makebox(0,0)[r]{\textcolor[rgb]{0.15,0.15,0.15}{{0}}}}
\fontsize{12}{0}\selectfont\put(43,145.53){\makebox(0,0)[r]{\textcolor[rgb]{0.15,0.15,0.15}{{0.5}}}}
\fontsize{12}{0}\selectfont\put(43,172.665){\makebox(0,0)[r]{\textcolor[rgb]{0.15,0.15,0.15}{{1}}}}
\fontsize{12}{0}\selectfont\put(43,199.8){\makebox(0,0)[r]{\textcolor[rgb]{0.15,0.15,0.15}{{1.5}}}}
\fontsize{13}{0}\selectfont\put(154.815,13){\makebox(0,0)[t]{\textcolor[rgb]{0.15,0.15,0.15}{{$t$}}}}
\fontsize{13}{0}\selectfont\put(17,118.395){\rotatebox{90}{\makebox(0,0)[b]{\textcolor[rgb]{0.15,0.15,0.15}{{$\frac{q^*}{q_0}, \bar{S}$}}}}}
\fontsize{10}{0}\selectfont\put(196.624,180.284){\makebox(0,0)[l]{\textcolor[rgb]{0,0,0}{{$\frac{q^*}{q_0}$}}}}
\fontsize{10}{0}\selectfont\put(196.624,163.279){\makebox(0,0)[l]{\textcolor[rgb]{0,0,0}{{$\bar{S}$}}}}
\end{picture}

%% file: path2I1.tex
\setlength{\unitlength}{1pt}
\begin{picture}(0,0)
\includegraphics[scale=1]{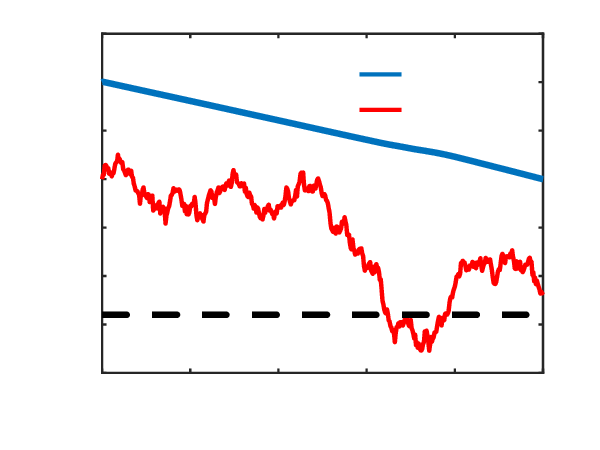}
\end{picture}%
\begin{picture}(288,216)(0,0)
\fontsize{12}{0}\selectfont\put(48.9901,28){\makebox(0,0)[t]{\textcolor[rgb]{0.15,0.15,0.15}{{0}}}}
\fontsize{12}{0}\selectfont\put(91.3201,28){\makebox(0,0)[t]{\textcolor[rgb]{0.15,0.15,0.15}{{0.2}}}}
\fontsize{12}{0}\selectfont\put(133.65,28){\makebox(0,0)[t]{\textcolor[rgb]{0.15,0.15,0.15}{{0.4}}}}
\fontsize{12}{0}\selectfont\put(175.98,28){\makebox(0,0)[t]{\textcolor[rgb]{0.15,0.15,0.15}{{0.6}}}}
\fontsize{12}{0}\selectfont\put(218.31,28){\makebox(0,0)[t]{\textcolor[rgb]{0.15,0.15,0.15}{{0.8}}}}
\fontsize{12}{0}\selectfont\put(260.64,28){\makebox(0,0)[t]{\textcolor[rgb]{0.15,0.15,0.15}{{1}}}}
\fontsize{12}{0}\selectfont\put(43,36.9895){\makebox(0,0)[r]{\textcolor[rgb]{0.15,0.15,0.15}{{-2}}}}
\fontsize{12}{0}\selectfont\put(43,60.2482){\makebox(0,0)[r]{\textcolor[rgb]{0.15,0.15,0.15}{{-1.5}}}}
\fontsize{12}{0}\selectfont\put(43,83.5068){\makebox(0,0)[r]{\textcolor[rgb]{0.15,0.15,0.15}{{-1}}}}
\fontsize{12}{0}\selectfont\put(43,106.765){\makebox(0,0)[r]{\textcolor[rgb]{0.15,0.15,0.15}{{-0.5}}}}
\fontsize{12}{0}\selectfont\put(43,130.024){\makebox(0,0)[r]{\textcolor[rgb]{0.15,0.15,0.15}{{0}}}}
\fontsize{12}{0}\selectfont\put(43,153.283){\makebox(0,0)[r]{\textcolor[rgb]{0.15,0.15,0.15}{{0.5}}}}
\fontsize{12}{0}\selectfont\put(43,176.541){\makebox(0,0)[r]{\textcolor[rgb]{0.15,0.15,0.15}{{1}}}}
\fontsize{12}{0}\selectfont\put(43,199.8){\makebox(0,0)[r]{\textcolor[rgb]{0.15,0.15,0.15}{{1.5}}}}
\fontsize{13}{0}\selectfont\put(154.815,13){\makebox(0,0)[t]{\textcolor[rgb]{0.15,0.15,0.15}{{$t$}}}}
\fontsize{13}{0}\selectfont\put(17,118.395){\rotatebox{90}{\makebox(0,0)[b]{\textcolor[rgb]{0.15,0.15,0.15}{{$\frac{q^*}{q_0}, \bar{S}$}}}}}
\fontsize{10}{0}\selectfont\put(196.624,180.284){\makebox(0,0)[l]{\textcolor[rgb]{0,0,0}{{$\frac{q^*}{q_0}$}}}}
\fontsize{10}{0}\selectfont\put(196.624,163.279){\makebox(0,0)[l]{\textcolor[rgb]{0,0,0}{{$\bar{S}$}}}}
\end{picture}

%% file: path3I1.tex
\setlength{\unitlength}{1pt}
\begin{picture}(0,0)
\includegraphics[scale=1]{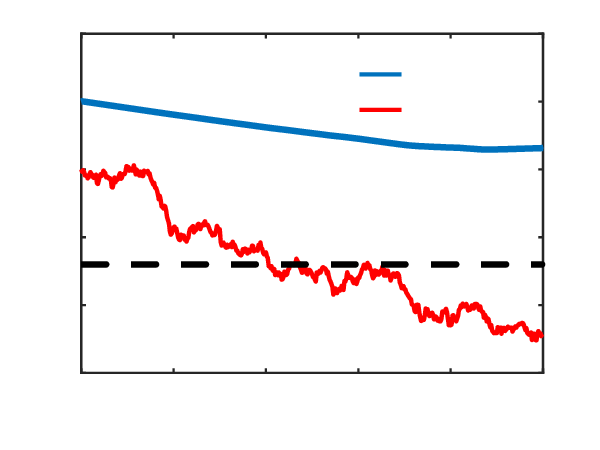}
\end{picture}%
\begin{picture}(288,216)(0,0)
\fontsize{12}{0}\selectfont\put(38.9905,28){\makebox(0,0)[t]{\textcolor[rgb]{0.15,0.15,0.15}{{0}}}}
\fontsize{12}{0}\selectfont\put(83.3204,28){\makebox(0,0)[t]{\textcolor[rgb]{0.15,0.15,0.15}{{0.2}}}}
\fontsize{12}{0}\selectfont\put(127.65,28){\makebox(0,0)[t]{\textcolor[rgb]{0.15,0.15,0.15}{{0.4}}}}
\fontsize{12}{0}\selectfont\put(171.98,28){\makebox(0,0)[t]{\textcolor[rgb]{0.15,0.15,0.15}{{0.6}}}}
\fontsize{12}{0}\selectfont\put(216.31,28){\makebox(0,0)[t]{\textcolor[rgb]{0.15,0.15,0.15}{{0.8}}}}
\fontsize{12}{0}\selectfont\put(260.64,28){\makebox(0,0)[t]{\textcolor[rgb]{0.15,0.15,0.15}{{1}}}}
\fontsize{12}{0}\selectfont\put(33,36.9895){\makebox(0,0)[r]{\textcolor[rgb]{0.15,0.15,0.15}{{-3}}}}
\fontsize{12}{0}\selectfont\put(33,69.5516){\makebox(0,0)[r]{\textcolor[rgb]{0.15,0.15,0.15}{{-2}}}}
\fontsize{12}{0}\selectfont\put(33,102.114){\makebox(0,0)[r]{\textcolor[rgb]{0.15,0.15,0.15}{{-1}}}}
\fontsize{12}{0}\selectfont\put(33,134.676){\makebox(0,0)[r]{\textcolor[rgb]{0.15,0.15,0.15}{{0}}}}
\fontsize{12}{0}\selectfont\put(33,167.238){\makebox(0,0)[r]{\textcolor[rgb]{0.15,0.15,0.15}{{1}}}}
\fontsize{12}{0}\selectfont\put(33,199.8){\makebox(0,0)[r]{\textcolor[rgb]{0.15,0.15,0.15}{{2}}}}
\fontsize{13}{0}\selectfont\put(17,118.395){\rotatebox{90}{\makebox(0,0)[b]{\textcolor[rgb]{0.15,0.15,0.15}{{$\frac{q^*}{q_0}, \bar{S}$}}}}}
\fontsize{13}{0}\selectfont\put(149.815,13){\makebox(0,0)[t]{\textcolor[rgb]{0.15,0.15,0.15}{{$t$}}}}
\fontsize{10}{0}\selectfont\put(196.624,180.284){\makebox(0,0)[l]{\textcolor[rgb]{0,0,0}{{$\frac{q^*}{q_0}$}}}}
\fontsize{10}{0}\selectfont\put(196.624,163.279){\makebox(0,0)[l]{\textcolor[rgb]{0,0,0}{{$\bar{S}$}}}}
\end{picture}

%% file: Ptail.tex
\setlength{\unitlength}{1pt}
\begin{picture}(0,0)
\includegraphics[scale=1]{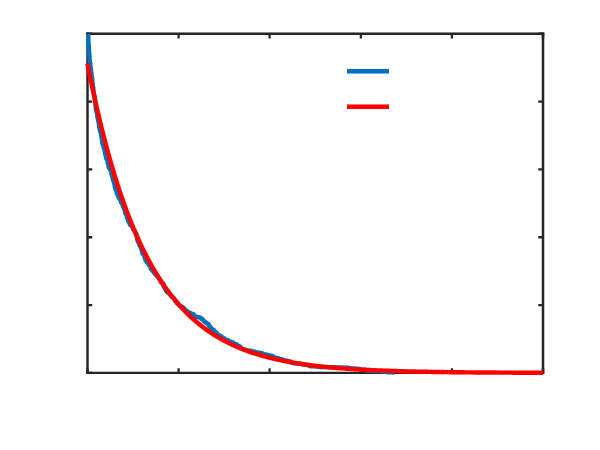}
\end{picture}%
\begin{picture}(288,216)(0,0)
\fontsize{12}{0}\selectfont\put(41.9904,28){\makebox(0,0)[t]{\textcolor[rgb]{0.15,0.15,0.15}{{0}}}}
\fontsize{12}{0}\selectfont\put(85.7203,28){\makebox(0,0)[t]{\textcolor[rgb]{0.15,0.15,0.15}{{0.2}}}}
\fontsize{12}{0}\selectfont\put(129.45,28){\makebox(0,0)[t]{\textcolor[rgb]{0.15,0.15,0.15}{{0.4}}}}
\fontsize{12}{0}\selectfont\put(173.18,28){\makebox(0,0)[t]{\textcolor[rgb]{0.15,0.15,0.15}{{0.6}}}}
\fontsize{12}{0}\selectfont\put(216.91,28){\makebox(0,0)[t]{\textcolor[rgb]{0.15,0.15,0.15}{{0.8}}}}
\fontsize{12}{0}\selectfont\put(260.64,28){\makebox(0,0)[t]{\textcolor[rgb]{0.15,0.15,0.15}{{1}}}}
\fontsize{12}{0}\selectfont\put(36,36.9895){\makebox(0,0)[r]{\textcolor[rgb]{0.15,0.15,0.15}{{0}}}}
\fontsize{12}{0}\selectfont\put(36,69.5516){\makebox(0,0)[r]{\textcolor[rgb]{0.15,0.15,0.15}{{0.2}}}}
\fontsize{12}{0}\selectfont\put(36,102.114){\makebox(0,0)[r]{\textcolor[rgb]{0.15,0.15,0.15}{{0.4}}}}
\fontsize{12}{0}\selectfont\put(36,134.676){\makebox(0,0)[r]{\textcolor[rgb]{0.15,0.15,0.15}{{0.6}}}}
\fontsize{12}{0}\selectfont\put(36,167.238){\makebox(0,0)[r]{\textcolor[rgb]{0.15,0.15,0.15}{{0.8}}}}
\fontsize{12}{0}\selectfont\put(36,199.8){\makebox(0,0)[r]{\textcolor[rgb]{0.15,0.15,0.15}{{1}}}}
\fontsize{13}{0}\selectfont\put(14,118.395){\rotatebox{90}{\makebox(0,0)[b]{\textcolor[rgb]{0.15,0.15,0.15}{{$p_x$}}}}}
\fontsize{13}{0}\selectfont\put(151.315,13){\makebox(0,0)[t]{\textcolor[rgb]{0.15,0.15,0.15}{{$x$}}}}
\fontsize{10}{0}\selectfont\put(190.624,181.784){\makebox(0,0)[l]{\textcolor[rgb]{0,0,0}{{$p_x$}}}}
\fontsize{10}{0}\selectfont\put(190.624,164.779){\makebox(0,0)[l]{\textcolor[rgb]{0,0,0}{{$e^{-\frac{x q_0}{m_T}}$}}}}
\end{picture}

%% file: logP.tex
\setlength{\unitlength}{1pt}
\begin{picture}(0,0)
\includegraphics[scale=1]{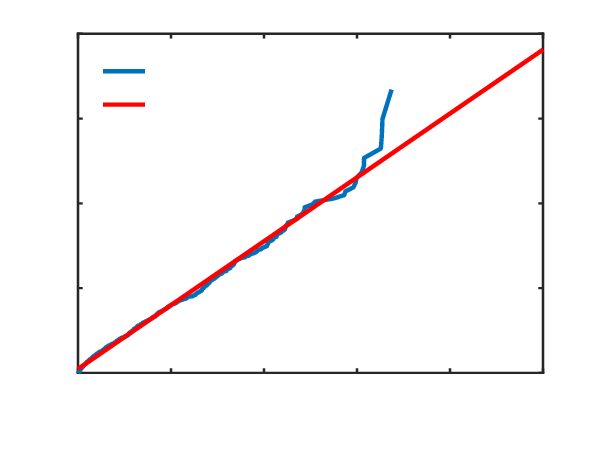}
\end{picture}%
\begin{picture}(288,216)(0,0)
\fontsize{12}{0}\selectfont\put(37.44,28){\makebox(0,0)[t]{\textcolor[rgb]{0.15,0.15,0.15}{{0}}}}
\fontsize{12}{0}\selectfont\put(82.08,28){\makebox(0,0)[t]{\textcolor[rgb]{0.15,0.15,0.15}{{0.2}}}}
\fontsize{12}{0}\selectfont\put(126.72,28){\makebox(0,0)[t]{\textcolor[rgb]{0.15,0.15,0.15}{{0.4}}}}
\fontsize{12}{0}\selectfont\put(171.36,28){\makebox(0,0)[t]{\textcolor[rgb]{0.15,0.15,0.15}{{0.6}}}}
\fontsize{12}{0}\selectfont\put(216,28){\makebox(0,0)[t]{\textcolor[rgb]{0.15,0.15,0.15}{{0.8}}}}
\fontsize{12}{0}\selectfont\put(260.64,28){\makebox(0,0)[t]{\textcolor[rgb]{0.15,0.15,0.15}{{1}}}}
\fontsize{12}{0}\selectfont\put(31.4346,36.9895){\makebox(0,0)[r]{\textcolor[rgb]{0.15,0.15,0.15}{{0}}}}
\fontsize{12}{0}\selectfont\put(31.4346,77.6922){\makebox(0,0)[r]{\textcolor[rgb]{0.15,0.15,0.15}{{2}}}}
\fontsize{12}{0}\selectfont\put(31.4346,118.395){\makebox(0,0)[r]{\textcolor[rgb]{0.15,0.15,0.15}{{4}}}}
\fontsize{12}{0}\selectfont\put(31.4346,159.097){\makebox(0,0)[r]{\textcolor[rgb]{0.15,0.15,0.15}{{6}}}}
\fontsize{12}{0}\selectfont\put(31.4346,199.8){\makebox(0,0)[r]{\textcolor[rgb]{0.15,0.15,0.15}{{8}}}}
\fontsize{13}{0}\selectfont\put(149.04,13){\makebox(0,0)[t]{\textcolor[rgb]{0.15,0.15,0.15}{{$x$}}}}
\fontsize{13}{0}\selectfont\put(19.4346,118.395){\rotatebox{90}{\makebox(0,0)[b]{\textcolor[rgb]{0.15,0.15,0.15}{{{\small $-\log(p_x), \frac{x q_0}{m_T}$}}}}}}
\fontsize{10}{0}\selectfont\put(73.4874,181.784){\makebox(0,0)[l]{\textcolor[rgb]{0,0,0}{{$-\log(p_x)$}}}}
\fontsize{10}{0}\selectfont\put(73.4874,165.779){\makebox(0,0)[l]{\textcolor[rgb]{0,0,0}{{$\frac{x q_0}{m_T}$}}}}
\end{picture}

%% file: varlI1.tex
\setlength{\unitlength}{1pt}
\begin{picture}(0,0)
\includegraphics[scale=1]{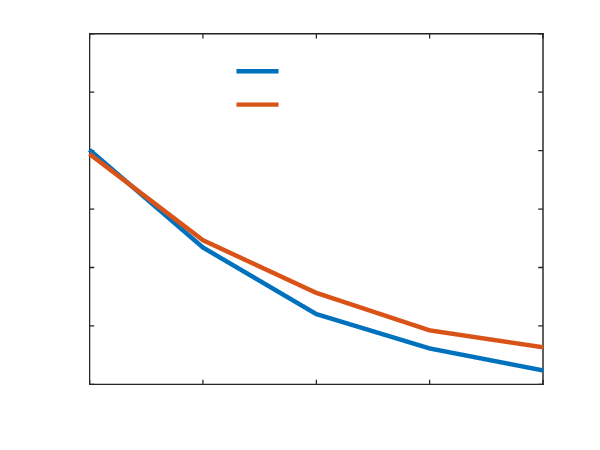}
\end{picture}%
\begin{picture}(288,216)(0,0)
\fontsize{10}{0}\selectfont\put(42.9919,24){\makebox(0,0)[t]{\textcolor[rgb]{0.15,0.15,0.15}{{0.5}}}}
\fontsize{10}{0}\selectfont\put(97.4039,24){\makebox(0,0)[t]{\textcolor[rgb]{0.15,0.15,0.15}{{1}}}}
\fontsize{10}{0}\selectfont\put(151.816,24){\makebox(0,0)[t]{\textcolor[rgb]{0.15,0.15,0.15}{{1.5}}}}
\fontsize{10}{0}\selectfont\put(206.228,24){\makebox(0,0)[t]{\textcolor[rgb]{0.15,0.15,0.15}{{2}}}}
\fontsize{10}{0}\selectfont\put(260.64,24){\makebox(0,0)[t]{\textcolor[rgb]{0.15,0.15,0.15}{{2.5}}}}
\fontsize{10}{0}\selectfont\put(38,31.5128){\makebox(0,0)[r]{\textcolor[rgb]{0.15,0.15,0.15}{{0.05}}}}
\fontsize{10}{0}\selectfont\put(38,59.5607){\makebox(0,0)[r]{\textcolor[rgb]{0.15,0.15,0.15}{{0.1}}}}
\fontsize{10}{0}\selectfont\put(38,87.6085){\makebox(0,0)[r]{\textcolor[rgb]{0.15,0.15,0.15}{{0.15}}}}
\fontsize{10}{0}\selectfont\put(38,115.656){\makebox(0,0)[r]{\textcolor[rgb]{0.15,0.15,0.15}{{0.2}}}}
\fontsize{10}{0}\selectfont\put(38,143.704){\makebox(0,0)[r]{\textcolor[rgb]{0.15,0.15,0.15}{{0.25}}}}
\fontsize{10}{0}\selectfont\put(38,171.752){\makebox(0,0)[r]{\textcolor[rgb]{0.15,0.15,0.15}{{0.3}}}}
\fontsize{10}{0}\selectfont\put(38,199.8){\makebox(0,0)[r]{\textcolor[rgb]{0.15,0.15,0.15}{{0.35}}}}
\fontsize{11}{0}\selectfont\put(12,115.656){\rotatebox{90}{\makebox(0,0)[b]{\textcolor[rgb]{0.15,0.15,0.15}{{$m_T, \sqrt{\var_T}$}}}}}
\fontsize{11}{0}\selectfont\put(151.816,11){\makebox(0,0)[t]{\textcolor[rgb]{0.15,0.15,0.15}{{$\ell$}}}}
\fontsize{9}{0}\selectfont\put(137.624,181.784){\makebox(0,0)[l]{\textcolor[rgb]{0,0,0}{{${\mathbb E}[\fq_T  | \fq_T > 0]$}}}}
\fontsize{9}{0}\selectfont\put(137.624,165.779){\makebox(0,0)[l]{\textcolor[rgb]{0,0,0}{{$\var(\fq_T   | \fq_T > 0)^{1/2}$}}}}
\end{picture}

%% file: varkI1.tex
\setlength{\unitlength}{1pt}
\begin{picture}(0,0)
\includegraphics[scale=1]{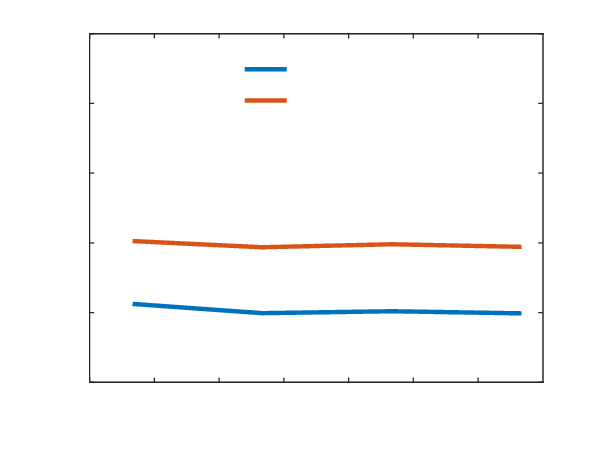}
\end{picture}%
\begin{picture}(288,216)(0,0)
\fontsize{10}{0}\selectfont\put(42.9919,25){\makebox(0,0)[t]{\textcolor[rgb]{0.15,0.15,0.15}{{0}}}}
\fontsize{10}{0}\selectfont\put(74.0845,25){\makebox(0,0)[t]{\textcolor[rgb]{0.15,0.15,0.15}{{0.2}}}}
\fontsize{10}{0}\selectfont\put(105.177,25){\makebox(0,0)[t]{\textcolor[rgb]{0.15,0.15,0.15}{{0.4}}}}
\fontsize{10}{0}\selectfont\put(136.27,25){\makebox(0,0)[t]{\textcolor[rgb]{0.15,0.15,0.15}{{0.6}}}}
\fontsize{10}{0}\selectfont\put(167.362,25){\makebox(0,0)[t]{\textcolor[rgb]{0.15,0.15,0.15}{{0.8}}}}
\fontsize{10}{0}\selectfont\put(198.455,25){\makebox(0,0)[t]{\textcolor[rgb]{0.15,0.15,0.15}{{1}}}}
\fontsize{10}{0}\selectfont\put(229.547,25){\makebox(0,0)[t]{\textcolor[rgb]{0.15,0.15,0.15}{{1.2}}}}
\fontsize{10}{0}\selectfont\put(38,32.5129){\makebox(0,0)[r]{\textcolor[rgb]{0.15,0.15,0.15}{{0.1}}}}
\fontsize{10}{0}\selectfont\put(38,65.9703){\makebox(0,0)[r]{\textcolor[rgb]{0.15,0.15,0.15}{{0.12}}}}
\fontsize{10}{0}\selectfont\put(38,99.4277){\makebox(0,0)[r]{\textcolor[rgb]{0.15,0.15,0.15}{{0.14}}}}
\fontsize{10}{0}\selectfont\put(38,132.885){\makebox(0,0)[r]{\textcolor[rgb]{0.15,0.15,0.15}{{0.16}}}}
\fontsize{10}{0}\selectfont\put(38,166.343){\makebox(0,0)[r]{\textcolor[rgb]{0.15,0.15,0.15}{{0.18}}}}
\fontsize{10}{0}\selectfont\put(38,199.8){\makebox(0,0)[r]{\textcolor[rgb]{0.15,0.15,0.15}{{0.2}}}}
\fontsize{11}{0}\selectfont\put(151.816,12){\makebox(0,0)[t]{\textcolor[rgb]{0.15,0.15,0.15}{{$kV/\eta$}}}}
\fontsize{11}{0}\selectfont\put(12,116.156){\rotatebox{90}{\makebox(0,0)[b]{\textcolor[rgb]{0.15,0.15,0.15}{{$m_T, \sqrt{\var_T}$}}}}}
\fontsize{9}{0}\selectfont\put(141.624,182.784){\makebox(0,0)[l]{\textcolor[rgb]{0,0,0}{{${\mathbb E}[\fq_T | \fq_T > 0]$}}}}
\fontsize{9}{0}\selectfont\put(141.624,167.779){\makebox(0,0)[l]{\textcolor[rgb]{0,0,0}{{$\var(\fq_T  | \fq_T > 0)^{1/2}$}}}}
\end{picture}

%% file: S2primegivenqTmeanvar.tex
\setlength{\unitlength}{1pt}
\begin{picture}(0,0)
\includegraphics[scale=1]{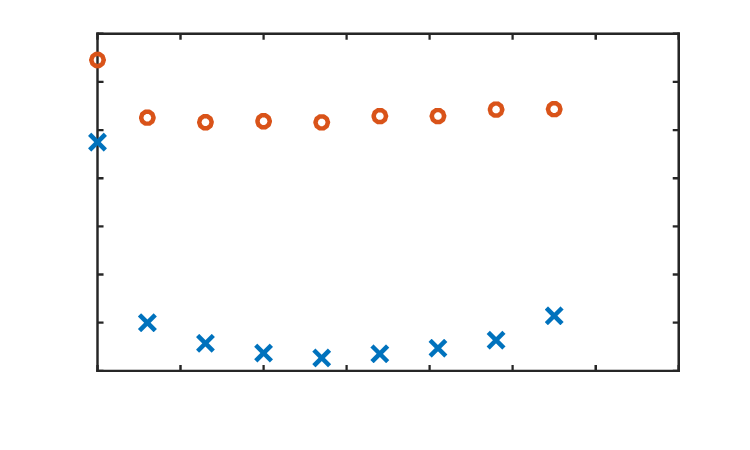}
\end{picture}%
\begin{picture}(360,216)(0,0)
\fontsize{12}{0}\selectfont\put(46.8,29){\makebox(0,0)[t]{\textcolor[rgb]{0.15,0.15,0.15}{{0}}}}
\fontsize{12}{0}\selectfont\put(86.6572,29){\makebox(0,0)[t]{\textcolor[rgb]{0.15,0.15,0.15}{{0.1}}}}
\fontsize{12}{0}\selectfont\put(126.514,29){\makebox(0,0)[t]{\textcolor[rgb]{0.15,0.15,0.15}{{0.2}}}}
\fontsize{12}{0}\selectfont\put(166.371,29){\makebox(0,0)[t]{\textcolor[rgb]{0.15,0.15,0.15}{{0.3}}}}
\fontsize{12}{0}\selectfont\put(206.229,29){\makebox(0,0)[t]{\textcolor[rgb]{0.15,0.15,0.15}{{0.4}}}}
\fontsize{12}{0}\selectfont\put(246.086,29){\makebox(0,0)[t]{\textcolor[rgb]{0.15,0.15,0.15}{{0.5}}}}
\fontsize{12}{0}\selectfont\put(285.943,29){\makebox(0,0)[t]{\textcolor[rgb]{0.15,0.15,0.15}{{0.6}}}}
\fontsize{12}{0}\selectfont\put(325.8,29){\makebox(0,0)[t]{\textcolor[rgb]{0.15,0.15,0.15}{{0.7}}}}
\fontsize{12}{0}\selectfont\put(40.8,37.9895){\makebox(0,0)[r]{\textcolor[rgb]{0.15,0.15,0.15}{{-4}}}}
\fontsize{12}{0}\selectfont\put(40.8,61.1052){\makebox(0,0)[r]{\textcolor[rgb]{0.15,0.15,0.15}{{-3}}}}
\fontsize{12}{0}\selectfont\put(40.8,84.221){\makebox(0,0)[r]{\textcolor[rgb]{0.15,0.15,0.15}{{-2}}}}
\fontsize{12}{0}\selectfont\put(40.8,107.337){\makebox(0,0)[r]{\textcolor[rgb]{0.15,0.15,0.15}{{-1}}}}
\fontsize{12}{0}\selectfont\put(40.8,130.453){\makebox(0,0)[r]{\textcolor[rgb]{0.15,0.15,0.15}{{0}}}}
\fontsize{12}{0}\selectfont\put(40.8,153.568){\makebox(0,0)[r]{\textcolor[rgb]{0.15,0.15,0.15}{{1}}}}
\fontsize{12}{0}\selectfont\put(40.8,176.684){\makebox(0,0)[r]{\textcolor[rgb]{0.15,0.15,0.15}{{2}}}}
\fontsize{12}{0}\selectfont\put(40.8,199.8){\makebox(0,0)[r]{\textcolor[rgb]{0.15,0.15,0.15}{{3}}}}
\fontsize{13}{0}\selectfont\put(24.8,118.895){\rotatebox{90}{\makebox(0,0)[b]{\textcolor[rgb]{0.15,0.15,0.15}{{${\mathbb E}[A_2| \fq_T ], \sqrt{\var(A_2|\fq_T)}$}}}}}
\fontsize{13}{0}\selectfont\put(186.3,14){\makebox(0,0)[t]{\textcolor[rgb]{0.15,0.15,0.15}{{$\fq_T$}}}}
\end{picture}

%% file: Wvaryl1.tex
\setlength{\unitlength}{1pt}
\begin{picture}(0,0)
\includegraphics[scale=1]{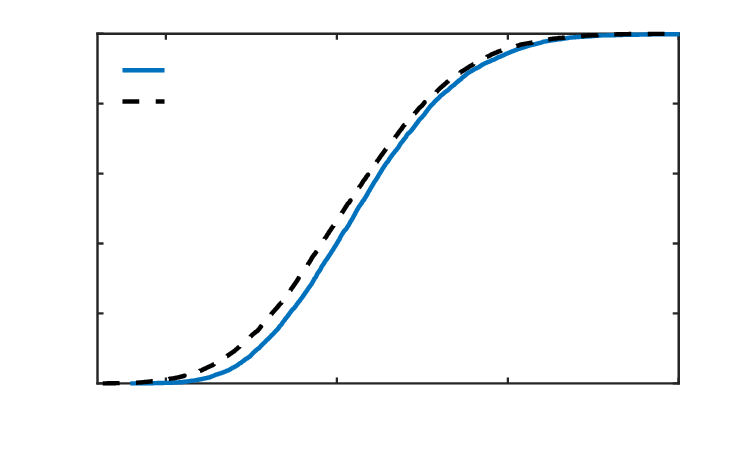}
\end{picture}%
\begin{picture}(360,216)(0,0)
\fontsize{12}{0}\selectfont\put(79.6235,23){\makebox(0,0)[t]{\textcolor[rgb]{0.15,0.15,0.15}{{-5}}}}
\fontsize{12}{0}\selectfont\put(161.682,23){\makebox(0,0)[t]{\textcolor[rgb]{0.15,0.15,0.15}{{0}}}}
\fontsize{12}{0}\selectfont\put(243.741,23){\makebox(0,0)[t]{\textcolor[rgb]{0.15,0.15,0.15}{{5}}}}
\fontsize{12}{0}\selectfont\put(325.8,23){\makebox(0,0)[t]{\textcolor[rgb]{0.15,0.15,0.15}{{10}}}}
\fontsize{12}{0}\selectfont\put(40.8,31.9898){\makebox(0,0)[r]{\textcolor[rgb]{0.15,0.15,0.15}{{0}}}}
\fontsize{12}{0}\selectfont\put(40.8,65.5519){\makebox(0,0)[r]{\textcolor[rgb]{0.15,0.15,0.15}{{0.2}}}}
\fontsize{12}{0}\selectfont\put(40.8,99.1139){\makebox(0,0)[r]{\textcolor[rgb]{0.15,0.15,0.15}{{0.4}}}}
\fontsize{12}{0}\selectfont\put(40.8,132.676){\makebox(0,0)[r]{\textcolor[rgb]{0.15,0.15,0.15}{{0.6}}}}
\fontsize{12}{0}\selectfont\put(40.8,166.238){\makebox(0,0)[r]{\textcolor[rgb]{0.15,0.15,0.15}{{0.8}}}}
\fontsize{12}{0}\selectfont\put(40.8,199.8){\makebox(0,0)[r]{\textcolor[rgb]{0.15,0.15,0.15}{{1}}}}
\fontsize{13}{0}\selectfont\put(186.3,7.99999){\makebox(0,0)[t]{\textcolor[rgb]{0.15,0.15,0.15}{{x}}}}
\fontsize{13}{0}\selectfont\put(18.8,115.895){\rotatebox{90}{\makebox(0,0)[b]{\textcolor[rgb]{0.15,0.15,0.15}{{$P(A_2 \le x | \fq_T^* = 0)$}}}}}
\fontsize{10}{0}\selectfont\put(82.8474,182.284){\makebox(0,0)[l]{\textcolor[rgb]{0,0,0}{{$\ell=-\sigma$}}}}
\fontsize{10}{0}\selectfont\put(82.8474,167.279){\makebox(0,0)[l]{\textcolor[rgb]{0,0,0}{{$\ell=-1.4\sigma$}}}}
\end{picture}

%% file: Wvaryl2.tex
\setlength{\unitlength}{1pt}
\begin{picture}(0,0)
\includegraphics[scale=1]{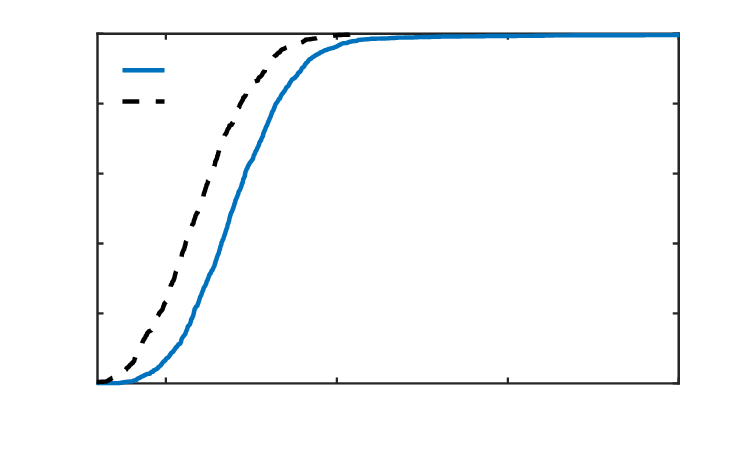}
\end{picture}%
\begin{picture}(360,216)(0,0)
\fontsize{12}{0}\selectfont\put(79.6235,23){\makebox(0,0)[t]{\textcolor[rgb]{0.15,0.15,0.15}{{-5}}}}
\fontsize{12}{0}\selectfont\put(161.682,23){\makebox(0,0)[t]{\textcolor[rgb]{0.15,0.15,0.15}{{0}}}}
\fontsize{12}{0}\selectfont\put(243.741,23){\makebox(0,0)[t]{\textcolor[rgb]{0.15,0.15,0.15}{{5}}}}
\fontsize{12}{0}\selectfont\put(325.8,23){\makebox(0,0)[t]{\textcolor[rgb]{0.15,0.15,0.15}{{10}}}}
\fontsize{12}{0}\selectfont\put(40.8,31.9898){\makebox(0,0)[r]{\textcolor[rgb]{0.15,0.15,0.15}{{0}}}}
\fontsize{12}{0}\selectfont\put(40.8,65.5519){\makebox(0,0)[r]{\textcolor[rgb]{0.15,0.15,0.15}{{0.2}}}}
\fontsize{12}{0}\selectfont\put(40.8,99.1139){\makebox(0,0)[r]{\textcolor[rgb]{0.15,0.15,0.15}{{0.4}}}}
\fontsize{12}{0}\selectfont\put(40.8,132.676){\makebox(0,0)[r]{\textcolor[rgb]{0.15,0.15,0.15}{{0.6}}}}
\fontsize{12}{0}\selectfont\put(40.8,166.238){\makebox(0,0)[r]{\textcolor[rgb]{0.15,0.15,0.15}{{0.8}}}}
\fontsize{12}{0}\selectfont\put(40.8,199.8){\makebox(0,0)[r]{\textcolor[rgb]{0.15,0.15,0.15}{{1}}}}
\fontsize{13}{0}\selectfont\put(186.3,7.99999){\makebox(0,0)[t]{\textcolor[rgb]{0.15,0.15,0.15}{{x}}}}
\fontsize{13}{0}\selectfont\put(18.8,115.895){\rotatebox{90}{\makebox(0,0)[b]{\textcolor[rgb]{0.15,0.15,0.15}{{$P(A_2 \le x | \fq_T > 0)$}}}}}
\fontsize{10}{0}\selectfont\put(82.8474,182.284){\makebox(0,0)[l]{\textcolor[rgb]{0,0,0}{{$\ell=-\sigma$}}}}
\fontsize{10}{0}\selectfont\put(82.8474,167.279){\makebox(0,0)[l]{\textcolor[rgb]{0,0,0}{{$\ell=-1.4\sigma$}}}}
\end{picture}

%% file: WvarykqT0I1.tex
\setlength{\unitlength}{1pt}
\begin{picture}(0,0)
\includegraphics[scale=1]{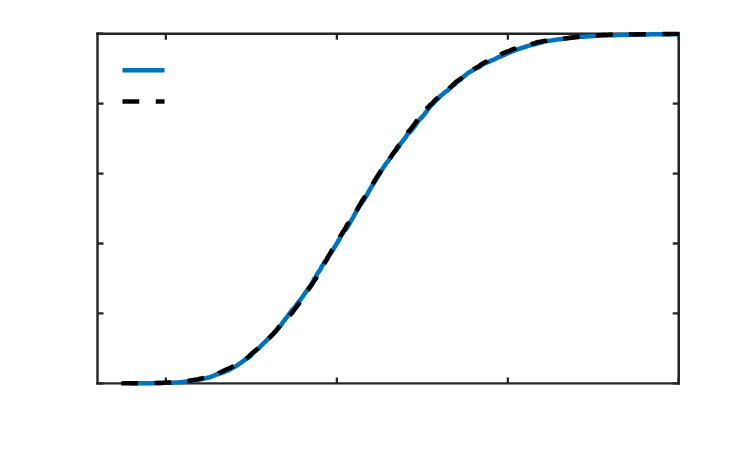}
\end{picture}%
\begin{picture}(360,216)(0,0)
\fontsize{12}{0}\selectfont\put(79.6235,23){\makebox(0,0)[t]{\textcolor[rgb]{0.15,0.15,0.15}{{-5}}}}
\fontsize{12}{0}\selectfont\put(161.682,23){\makebox(0,0)[t]{\textcolor[rgb]{0.15,0.15,0.15}{{0}}}}
\fontsize{12}{0}\selectfont\put(243.741,23){\makebox(0,0)[t]{\textcolor[rgb]{0.15,0.15,0.15}{{5}}}}
\fontsize{12}{0}\selectfont\put(325.8,23){\makebox(0,0)[t]{\textcolor[rgb]{0.15,0.15,0.15}{{10}}}}
\fontsize{12}{0}\selectfont\put(40.8,31.9898){\makebox(0,0)[r]{\textcolor[rgb]{0.15,0.15,0.15}{{0}}}}
\fontsize{12}{0}\selectfont\put(40.8,65.5519){\makebox(0,0)[r]{\textcolor[rgb]{0.15,0.15,0.15}{{0.2}}}}
\fontsize{12}{0}\selectfont\put(40.8,99.1139){\makebox(0,0)[r]{\textcolor[rgb]{0.15,0.15,0.15}{{0.4}}}}
\fontsize{12}{0}\selectfont\put(40.8,132.676){\makebox(0,0)[r]{\textcolor[rgb]{0.15,0.15,0.15}{{0.6}}}}
\fontsize{12}{0}\selectfont\put(40.8,166.238){\makebox(0,0)[r]{\textcolor[rgb]{0.15,0.15,0.15}{{0.8}}}}
\fontsize{12}{0}\selectfont\put(40.8,199.8){\makebox(0,0)[r]{\textcolor[rgb]{0.15,0.15,0.15}{{1}}}}
\fontsize{13}{0}\selectfont\put(186.3,7.99999){\makebox(0,0)[t]{\textcolor[rgb]{0.15,0.15,0.15}{{x}}}}
\fontsize{13}{0}\selectfont\put(18.8,115.895){\rotatebox{90}{\makebox(0,0)[b]{\textcolor[rgb]{0.15,0.15,0.15}{{$P(A_2 \le x | \fq_T = 0)$}}}}}
\fontsize{10}{0}\selectfont\put(82.8474,182.284){\makebox(0,0)[l]{\textcolor[rgb]{0,0,0}{{$kV/\eta=4/3$}}}}
\fontsize{10}{0}\selectfont\put(82.8474,167.279){\makebox(0,0)[l]{\textcolor[rgb]{0,0,0}{{$kV/\eta=4/30$}}}}
\end{picture}

%% file: WvarykqTg0I1.tex
\setlength{\unitlength}{1pt}
\begin{picture}(0,0)
\includegraphics[scale=1]{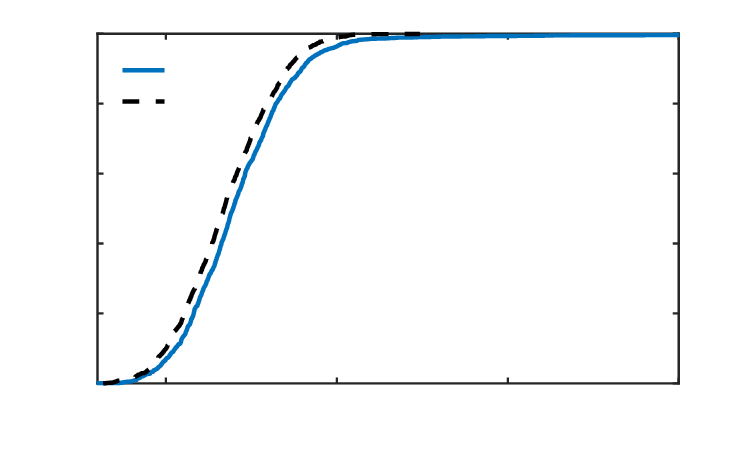}
\end{picture}%
\begin{picture}(360,216)(0,0)
\fontsize{12}{0}\selectfont\put(79.6235,23){\makebox(0,0)[t]{\textcolor[rgb]{0.15,0.15,0.15}{{-5}}}}
\fontsize{12}{0}\selectfont\put(161.682,23){\makebox(0,0)[t]{\textcolor[rgb]{0.15,0.15,0.15}{{0}}}}
\fontsize{12}{0}\selectfont\put(243.741,23){\makebox(0,0)[t]{\textcolor[rgb]{0.15,0.15,0.15}{{5}}}}
\fontsize{12}{0}\selectfont\put(325.8,23){\makebox(0,0)[t]{\textcolor[rgb]{0.15,0.15,0.15}{{10}}}}
\fontsize{12}{0}\selectfont\put(40.8,31.9898){\makebox(0,0)[r]{\textcolor[rgb]{0.15,0.15,0.15}{{0}}}}
\fontsize{12}{0}\selectfont\put(40.8,65.5519){\makebox(0,0)[r]{\textcolor[rgb]{0.15,0.15,0.15}{{0.2}}}}
\fontsize{12}{0}\selectfont\put(40.8,99.1139){\makebox(0,0)[r]{\textcolor[rgb]{0.15,0.15,0.15}{{0.4}}}}
\fontsize{12}{0}\selectfont\put(40.8,132.676){\makebox(0,0)[r]{\textcolor[rgb]{0.15,0.15,0.15}{{0.6}}}}
\fontsize{12}{0}\selectfont\put(40.8,166.238){\makebox(0,0)[r]{\textcolor[rgb]{0.15,0.15,0.15}{{0.8}}}}
\fontsize{12}{0}\selectfont\put(40.8,199.8){\makebox(0,0)[r]{\textcolor[rgb]{0.15,0.15,0.15}{{1}}}}
\fontsize{13}{0}\selectfont\put(186.3,7.99999){\makebox(0,0)[t]{\textcolor[rgb]{0.15,0.15,0.15}{{x}}}}
\fontsize{13}{0}\selectfont\put(18.8,115.895){\rotatebox{90}{\makebox(0,0)[b]{\textcolor[rgb]{0.15,0.15,0.15}{{$P(A_2 \le x | \fq_T > 0)$}}}}}
\fontsize{10}{0}\selectfont\put(82.8474,182.284){\makebox(0,0)[l]{\textcolor[rgb]{0,0,0}{{$k V/\eta =4/3$}}}}
\fontsize{10}{0}\selectfont\put(82.8474,167.279){\makebox(0,0)[l]{\textcolor[rgb]{0,0,0}{{$kV/\eta=4/30$}}}}
\end{picture}

%% file: A3qT0dist.tex
\setlength{\unitlength}{1pt}
\begin{picture}(0,0)
\includegraphics[scale=1]{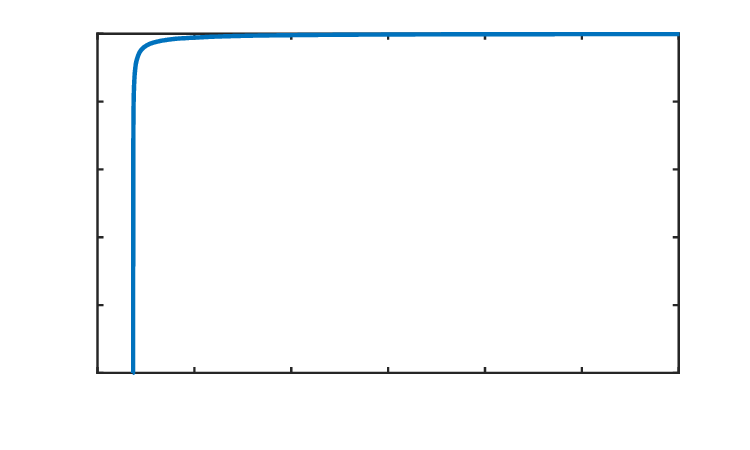}
\end{picture}%
\begin{picture}(360,216)(0,0)
\fontsize{12}{0}\selectfont\put(46.8,28){\makebox(0,0)[t]{\textcolor[rgb]{0.15,0.15,0.15}{{1}}}}
\fontsize{12}{0}\selectfont\put(93.3,28){\makebox(0,0)[t]{\textcolor[rgb]{0.15,0.15,0.15}{{1.05}}}}
\fontsize{12}{0}\selectfont\put(139.8,28){\makebox(0,0)[t]{\textcolor[rgb]{0.15,0.15,0.15}{{1.1}}}}
\fontsize{12}{0}\selectfont\put(186.3,28){\makebox(0,0)[t]{\textcolor[rgb]{0.15,0.15,0.15}{{1.15}}}}
\fontsize{12}{0}\selectfont\put(232.8,28){\makebox(0,0)[t]{\textcolor[rgb]{0.15,0.15,0.15}{{1.2}}}}
\fontsize{12}{0}\selectfont\put(279.3,28){\makebox(0,0)[t]{\textcolor[rgb]{0.15,0.15,0.15}{{1.25}}}}
\fontsize{12}{0}\selectfont\put(325.8,28){\makebox(0,0)[t]{\textcolor[rgb]{0.15,0.15,0.15}{{1.3}}}}
\fontsize{12}{0}\selectfont\put(40.8,36.9895){\makebox(0,0)[r]{\textcolor[rgb]{0.15,0.15,0.15}{{0}}}}
\fontsize{12}{0}\selectfont\put(40.8,69.5516){\makebox(0,0)[r]{\textcolor[rgb]{0.15,0.15,0.15}{{0.2}}}}
\fontsize{12}{0}\selectfont\put(40.8,102.114){\makebox(0,0)[r]{\textcolor[rgb]{0.15,0.15,0.15}{{0.4}}}}
\fontsize{12}{0}\selectfont\put(40.8,134.676){\makebox(0,0)[r]{\textcolor[rgb]{0.15,0.15,0.15}{{0.6}}}}
\fontsize{12}{0}\selectfont\put(40.8,167.238){\makebox(0,0)[r]{\textcolor[rgb]{0.15,0.15,0.15}{{0.8}}}}
\fontsize{12}{0}\selectfont\put(40.8,199.8){\makebox(0,0)[r]{\textcolor[rgb]{0.15,0.15,0.15}{{1}}}}
\fontsize{13}{0}\selectfont\put(186.3,13){\makebox(0,0)[t]{\textcolor[rgb]{0.15,0.15,0.15}{{$x$}}}}
\fontsize{13}{0}\selectfont\put(18.8,118.395){\rotatebox{90}{\makebox(0,0)[b]{\textcolor[rgb]{0.15,0.15,0.15}{{${\mathbb P}(A_3< \le x)$}}}}}
\end{picture}

%% file: A3vsqT.tex
\setlength{\unitlength}{1pt}
\begin{picture}(0,0)
\includegraphics[scale=1]{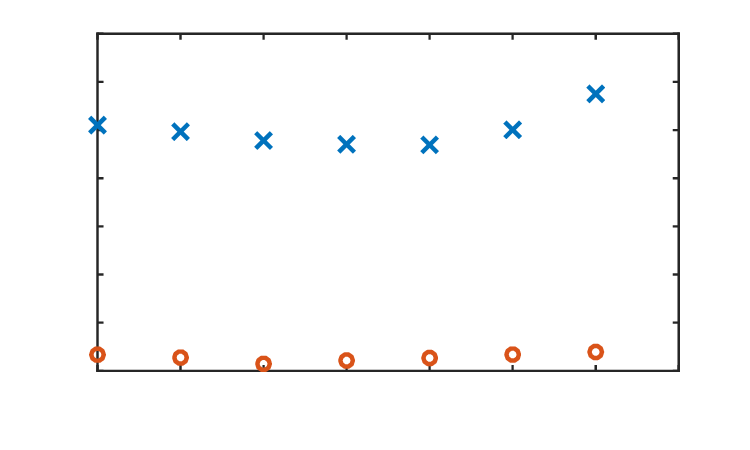}
\end{picture}%
\begin{picture}(360,216)(0,0)
\fontsize{12}{0}\selectfont\put(46.8,29){\makebox(0,0)[t]{\textcolor[rgb]{0.15,0.15,0.15}{{0}}}}
\fontsize{12}{0}\selectfont\put(86.6572,29){\makebox(0,0)[t]{\textcolor[rgb]{0.15,0.15,0.15}{{0.1}}}}
\fontsize{12}{0}\selectfont\put(126.514,29){\makebox(0,0)[t]{\textcolor[rgb]{0.15,0.15,0.15}{{0.2}}}}
\fontsize{12}{0}\selectfont\put(166.371,29){\makebox(0,0)[t]{\textcolor[rgb]{0.15,0.15,0.15}{{0.3}}}}
\fontsize{12}{0}\selectfont\put(206.229,29){\makebox(0,0)[t]{\textcolor[rgb]{0.15,0.15,0.15}{{0.4}}}}
\fontsize{12}{0}\selectfont\put(246.086,29){\makebox(0,0)[t]{\textcolor[rgb]{0.15,0.15,0.15}{{0.5}}}}
\fontsize{12}{0}\selectfont\put(285.943,29){\makebox(0,0)[t]{\textcolor[rgb]{0.15,0.15,0.15}{{0.6}}}}
\fontsize{12}{0}\selectfont\put(325.8,29){\makebox(0,0)[t]{\textcolor[rgb]{0.15,0.15,0.15}{{0.7}}}}
\fontsize{12}{0}\selectfont\put(40.8,37.9895){\makebox(0,0)[r]{\textcolor[rgb]{0.15,0.15,0.15}{{0}}}}
\fontsize{12}{0}\selectfont\put(40.8,61.1053){\makebox(0,0)[r]{\textcolor[rgb]{0.15,0.15,0.15}{{0.2}}}}
\fontsize{12}{0}\selectfont\put(40.8,84.221){\makebox(0,0)[r]{\textcolor[rgb]{0.15,0.15,0.15}{{0.4}}}}
\fontsize{12}{0}\selectfont\put(40.8,107.337){\makebox(0,0)[r]{\textcolor[rgb]{0.15,0.15,0.15}{{0.6}}}}
\fontsize{12}{0}\selectfont\put(40.8,130.453){\makebox(0,0)[r]{\textcolor[rgb]{0.15,0.15,0.15}{{0.8}}}}
\fontsize{12}{0}\selectfont\put(40.8,153.568){\makebox(0,0)[r]{\textcolor[rgb]{0.15,0.15,0.15}{{1}}}}
\fontsize{12}{0}\selectfont\put(40.8,176.684){\makebox(0,0)[r]{\textcolor[rgb]{0.15,0.15,0.15}{{1.2}}}}
\fontsize{12}{0}\selectfont\put(40.8,199.8){\makebox(0,0)[r]{\textcolor[rgb]{0.15,0.15,0.15}{{1.4}}}}
\fontsize{13}{0}\selectfont\put(18.8,118.895){\rotatebox{90}{\makebox(0,0)[b]{\textcolor[rgb]{0.15,0.15,0.15}{{${\mathbb E}[A_3| \fq_T ], \sqrt{\var(A_3|\fq_T)}$}}}}}
\fontsize{13}{0}\selectfont\put(186.3,14){\makebox(0,0)[t]{\textcolor[rgb]{0.15,0.15,0.15}{{$\fq_T$}}}}
\end{picture}

%% file: EtavarylqT0I1.tex
\setlength{\unitlength}{1pt}
\begin{picture}(0,0)
\includegraphics[scale=1]{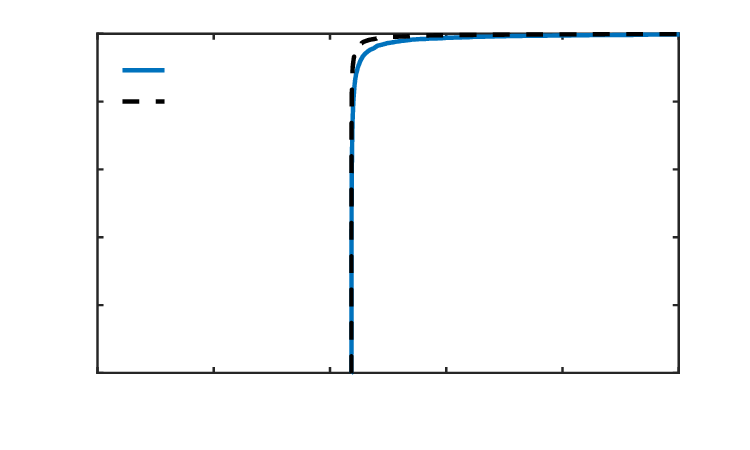}
\end{picture}%
\begin{picture}(360,216)(0,0)
\fontsize{12}{0}\selectfont\put(46.8,28){\makebox(0,0)[t]{\textcolor[rgb]{0.15,0.15,0.15}{{0.8}}}}
\fontsize{12}{0}\selectfont\put(102.6,28){\makebox(0,0)[t]{\textcolor[rgb]{0.15,0.15,0.15}{{0.9}}}}
\fontsize{12}{0}\selectfont\put(158.4,28){\makebox(0,0)[t]{\textcolor[rgb]{0.15,0.15,0.15}{{1}}}}
\fontsize{12}{0}\selectfont\put(214.2,28){\makebox(0,0)[t]{\textcolor[rgb]{0.15,0.15,0.15}{{1.1}}}}
\fontsize{12}{0}\selectfont\put(270,28){\makebox(0,0)[t]{\textcolor[rgb]{0.15,0.15,0.15}{{1.2}}}}
\fontsize{12}{0}\selectfont\put(325.8,28){\makebox(0,0)[t]{\textcolor[rgb]{0.15,0.15,0.15}{{1.3}}}}
\fontsize{12}{0}\selectfont\put(40.8,36.9895){\makebox(0,0)[r]{\textcolor[rgb]{0.15,0.15,0.15}{{0}}}}
\fontsize{12}{0}\selectfont\put(40.8,69.5516){\makebox(0,0)[r]{\textcolor[rgb]{0.15,0.15,0.15}{{0.2}}}}
\fontsize{12}{0}\selectfont\put(40.8,102.114){\makebox(0,0)[r]{\textcolor[rgb]{0.15,0.15,0.15}{{0.4}}}}
\fontsize{12}{0}\selectfont\put(40.8,134.676){\makebox(0,0)[r]{\textcolor[rgb]{0.15,0.15,0.15}{{0.6}}}}
\fontsize{12}{0}\selectfont\put(40.8,167.238){\makebox(0,0)[r]{\textcolor[rgb]{0.15,0.15,0.15}{{0.8}}}}
\fontsize{12}{0}\selectfont\put(40.8,199.8){\makebox(0,0)[r]{\textcolor[rgb]{0.15,0.15,0.15}{{1}}}}
\fontsize{13}{0}\selectfont\put(18.8,118.395){\rotatebox{90}{\makebox(0,0)[b]{\textcolor[rgb]{0.15,0.15,0.15}{{$P(A_3 \le x | \fq_T = 0)$}}}}}
\fontsize{13}{0}\selectfont\put(186.3,13){\makebox(0,0)[t]{\textcolor[rgb]{0.15,0.15,0.15}{{$x$}}}}
\fontsize{10}{0}\selectfont\put(82.8474,182.284){\makebox(0,0)[l]{\textcolor[rgb]{0,0,0}{{$\ell=-\sigma$}}}}
\fontsize{10}{0}\selectfont\put(82.8474,167.279){\makebox(0,0)[l]{\textcolor[rgb]{0,0,0}{{$\ell=-1.4\sigma$}}}}
\end{picture}

%% file: EtavarylqTg0I1.tex
\setlength{\unitlength}{1pt}
\begin{picture}(0,0)
\includegraphics[scale=1]{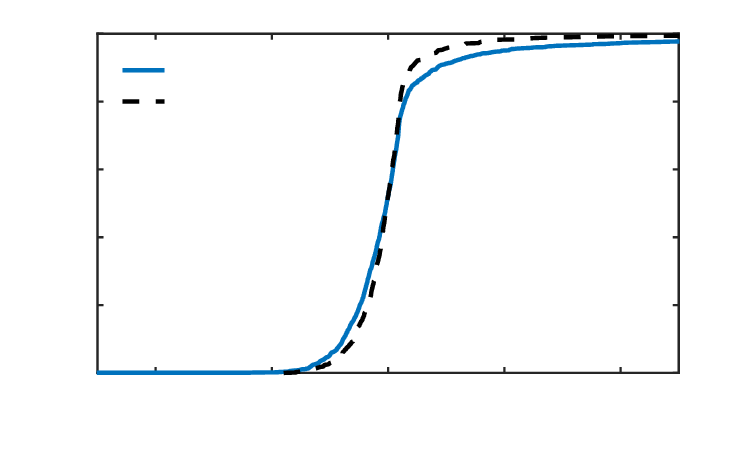}
\end{picture}%
\begin{picture}(360,216)(0,0)
\fontsize{12}{0}\selectfont\put(74.7,28){\makebox(0,0)[t]{\textcolor[rgb]{0.15,0.15,0.15}{{0.6}}}}
\fontsize{12}{0}\selectfont\put(130.5,28){\makebox(0,0)[t]{\textcolor[rgb]{0.15,0.15,0.15}{{0.8}}}}
\fontsize{12}{0}\selectfont\put(186.3,28){\makebox(0,0)[t]{\textcolor[rgb]{0.15,0.15,0.15}{{1}}}}
\fontsize{12}{0}\selectfont\put(242.1,28){\makebox(0,0)[t]{\textcolor[rgb]{0.15,0.15,0.15}{{1.2}}}}
\fontsize{12}{0}\selectfont\put(297.9,28){\makebox(0,0)[t]{\textcolor[rgb]{0.15,0.15,0.15}{{1.4}}}}
\fontsize{12}{0}\selectfont\put(40.8,36.9895){\makebox(0,0)[r]{\textcolor[rgb]{0.15,0.15,0.15}{{0}}}}
\fontsize{12}{0}\selectfont\put(40.8,69.5516){\makebox(0,0)[r]{\textcolor[rgb]{0.15,0.15,0.15}{{0.2}}}}
\fontsize{12}{0}\selectfont\put(40.8,102.114){\makebox(0,0)[r]{\textcolor[rgb]{0.15,0.15,0.15}{{0.4}}}}
\fontsize{12}{0}\selectfont\put(40.8,134.676){\makebox(0,0)[r]{\textcolor[rgb]{0.15,0.15,0.15}{{0.6}}}}
\fontsize{12}{0}\selectfont\put(40.8,167.238){\makebox(0,0)[r]{\textcolor[rgb]{0.15,0.15,0.15}{{0.8}}}}
\fontsize{12}{0}\selectfont\put(40.8,199.8){\makebox(0,0)[r]{\textcolor[rgb]{0.15,0.15,0.15}{{1}}}}
\fontsize{13}{0}\selectfont\put(18.8,118.395){\rotatebox{90}{\makebox(0,0)[b]{\textcolor[rgb]{0.15,0.15,0.15}{{$P(A_3 \le x | \fq_T > 0)$}}}}}
\fontsize{13}{0}\selectfont\put(186.3,13){\makebox(0,0)[t]{\textcolor[rgb]{0.15,0.15,0.15}{{$x$}}}}
\fontsize{10}{0}\selectfont\put(82.8474,182.284){\makebox(0,0)[l]{\textcolor[rgb]{0,0,0}{{$\ell=-\sigma$}}}}
\fontsize{10}{0}\selectfont\put(82.8474,167.279){\makebox(0,0)[l]{\textcolor[rgb]{0,0,0}{{$\ell=-1.4\sigma$}}}}
\end{picture}

%% file: EtavarykqT0I1.tex
\setlength{\unitlength}{1pt}
\begin{picture}(0,0)
\includegraphics[scale=1]{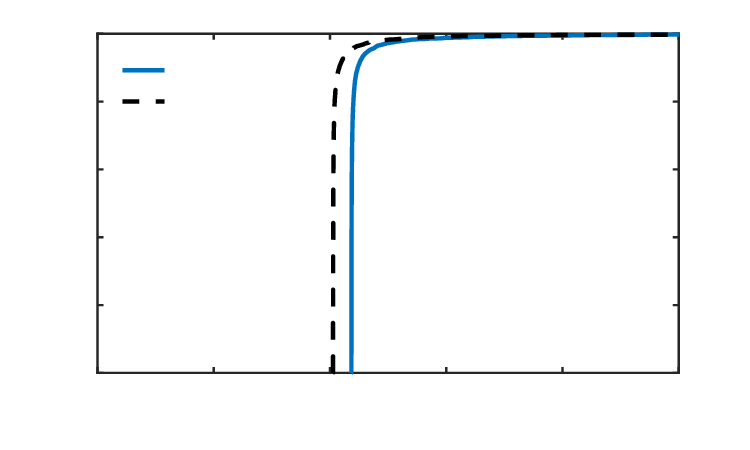}
\end{picture}%
\begin{picture}(360,216)(0,0)
\fontsize{12}{0}\selectfont\put(46.8,28){\makebox(0,0)[t]{\textcolor[rgb]{0.15,0.15,0.15}{{0.8}}}}
\fontsize{12}{0}\selectfont\put(102.6,28){\makebox(0,0)[t]{\textcolor[rgb]{0.15,0.15,0.15}{{0.9}}}}
\fontsize{12}{0}\selectfont\put(158.4,28){\makebox(0,0)[t]{\textcolor[rgb]{0.15,0.15,0.15}{{1}}}}
\fontsize{12}{0}\selectfont\put(214.2,28){\makebox(0,0)[t]{\textcolor[rgb]{0.15,0.15,0.15}{{1.1}}}}
\fontsize{12}{0}\selectfont\put(270,28){\makebox(0,0)[t]{\textcolor[rgb]{0.15,0.15,0.15}{{1.2}}}}
\fontsize{12}{0}\selectfont\put(325.8,28){\makebox(0,0)[t]{\textcolor[rgb]{0.15,0.15,0.15}{{1.3}}}}
\fontsize{12}{0}\selectfont\put(40.8,36.9895){\makebox(0,0)[r]{\textcolor[rgb]{0.15,0.15,0.15}{{0}}}}
\fontsize{12}{0}\selectfont\put(40.8,69.5516){\makebox(0,0)[r]{\textcolor[rgb]{0.15,0.15,0.15}{{0.2}}}}
\fontsize{12}{0}\selectfont\put(40.8,102.114){\makebox(0,0)[r]{\textcolor[rgb]{0.15,0.15,0.15}{{0.4}}}}
\fontsize{12}{0}\selectfont\put(40.8,134.676){\makebox(0,0)[r]{\textcolor[rgb]{0.15,0.15,0.15}{{0.6}}}}
\fontsize{12}{0}\selectfont\put(40.8,167.238){\makebox(0,0)[r]{\textcolor[rgb]{0.15,0.15,0.15}{{0.8}}}}
\fontsize{12}{0}\selectfont\put(40.8,199.8){\makebox(0,0)[r]{\textcolor[rgb]{0.15,0.15,0.15}{{1}}}}
\fontsize{13}{0}\selectfont\put(18.8,118.395){\rotatebox{90}{\makebox(0,0)[b]{\textcolor[rgb]{0.15,0.15,0.15}{{$P(A_3 \le x | \fq_T = 0)$}}}}}
\fontsize{13}{0}\selectfont\put(186.3,13){\makebox(0,0)[t]{\textcolor[rgb]{0.15,0.15,0.15}{{$x$}}}}
\fontsize{10}{0}\selectfont\put(82.8474,182.284){\makebox(0,0)[l]{\textcolor[rgb]{0,0,0}{{$kV/\eta=4/3$}}}}
\fontsize{10}{0}\selectfont\put(82.8474,167.279){\makebox(0,0)[l]{\textcolor[rgb]{0,0,0}{{$kV/\eta=4/30$}}}}
\end{picture}

%% file: EtavarykqTg0I1.tex
\setlength{\unitlength}{1pt}
\begin{picture}(0,0)
\includegraphics[scale=1]{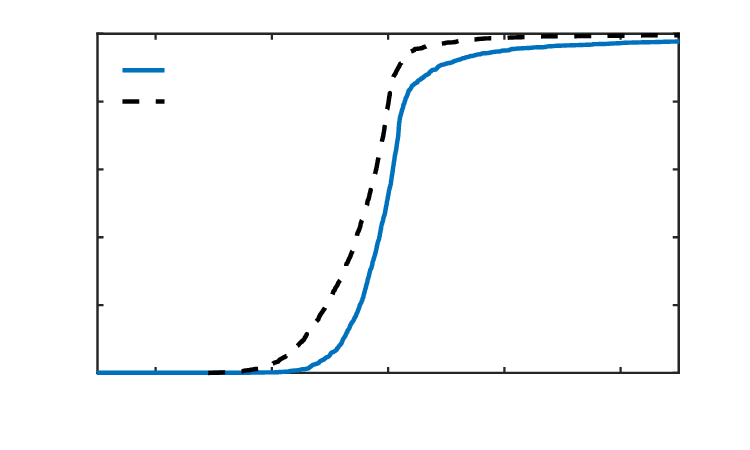}
\end{picture}%
\begin{picture}(360,216)(0,0)
\fontsize{12}{0}\selectfont\put(74.7,28){\makebox(0,0)[t]{\textcolor[rgb]{0.15,0.15,0.15}{{0.6}}}}
\fontsize{12}{0}\selectfont\put(130.5,28){\makebox(0,0)[t]{\textcolor[rgb]{0.15,0.15,0.15}{{0.8}}}}
\fontsize{12}{0}\selectfont\put(186.3,28){\makebox(0,0)[t]{\textcolor[rgb]{0.15,0.15,0.15}{{1}}}}
\fontsize{12}{0}\selectfont\put(242.1,28){\makebox(0,0)[t]{\textcolor[rgb]{0.15,0.15,0.15}{{1.2}}}}
\fontsize{12}{0}\selectfont\put(297.9,28){\makebox(0,0)[t]{\textcolor[rgb]{0.15,0.15,0.15}{{1.4}}}}
\fontsize{12}{0}\selectfont\put(40.8,36.9895){\makebox(0,0)[r]{\textcolor[rgb]{0.15,0.15,0.15}{{0}}}}
\fontsize{12}{0}\selectfont\put(40.8,69.5516){\makebox(0,0)[r]{\textcolor[rgb]{0.15,0.15,0.15}{{0.2}}}}
\fontsize{12}{0}\selectfont\put(40.8,102.114){\makebox(0,0)[r]{\textcolor[rgb]{0.15,0.15,0.15}{{0.4}}}}
\fontsize{12}{0}\selectfont\put(40.8,134.676){\makebox(0,0)[r]{\textcolor[rgb]{0.15,0.15,0.15}{{0.6}}}}
\fontsize{12}{0}\selectfont\put(40.8,167.238){\makebox(0,0)[r]{\textcolor[rgb]{0.15,0.15,0.15}{{0.8}}}}
\fontsize{12}{0}\selectfont\put(40.8,199.8){\makebox(0,0)[r]{\textcolor[rgb]{0.15,0.15,0.15}{{1}}}}
\fontsize{13}{0}\selectfont\put(18.8,118.395){\rotatebox{90}{\makebox(0,0)[b]{\textcolor[rgb]{0.15,0.15,0.15}{{$P(A_3 \le x | \fq_T > 0)$}}}}}
\fontsize{13}{0}\selectfont\put(186.3,13){\makebox(0,0)[t]{\textcolor[rgb]{0.15,0.15,0.15}{{$x$}}}}
\fontsize{10}{0}\selectfont\put(82.8474,182.284){\makebox(0,0)[l]{\textcolor[rgb]{0,0,0}{{$k V/\eta =4/3$}}}}
\fontsize{10}{0}\selectfont\put(82.8474,167.279){\makebox(0,0)[l]{\textcolor[rgb]{0,0,0}{{$kV/\eta=4/30$}}}}
\end{picture}

%% file: olv24.bbl
\def\cprime{$'$}
\begin{thebibliography}{10}

\bibitem{mahd:popi:seze:21}
M.~Ahmadi, A.~Popier, and A.~D. Sezer.
\newblock Backward stochastic differential equations with non-{M}arkovian
  singular terminal conditions for general driver and filtration.
\newblock {\em Electron. J. Probab.}, 26:Paper No. 64, 27, 2021.

\bibitem{almg:12}
R.~Almgren.
\newblock Optimal trading with stochastic liquidity and volatility.
\newblock {\em SIAM Journal on Financial Mathematics}, 3(1):163--181, 2012.

\bibitem{almg:chri:01}
R.~Almgren and N.~Chriss.
\newblock Optimal execution of portfolio transactions.
\newblock {\em Journal of Risk}, 3:5--40, 2001.

\bibitem{anki:jean:krus:13}
S.~Ankirchner, M.~Jeanblanc, and T.~Kruse.
\newblock B{SDE}s with {S}ingular {T}erminal {C}ondition and a {C}ontrol
  {P}roblem with {C}onstraints.
\newblock {\em SIAM J. Control Optim.}, 52(2):893--913, 2014.

\bibitem{barl:buck:pard:97}
G.~Barles, R.~Buckdahn, and E.~Pardoux.
\newblock Backward stochastic differential equations and integral-partial
  differential equations.
\newblock {\em Stochastics Stochastics Rep.}, 60(1-2):57--83, 1997.

\bibitem{buck:zuza:gunt:16}
Z.~Bu{\v{c}}kov{\'a}, M.~Ehrhardt, and M.~G{\"u}nther.
\newblock Fichera theory and its application in finance.
\newblock In G.~Russo, V.~Capasso, G.~Nicosia, and V.~Romano, editors, {\em
  Progress in Industrial Mathematics at ECMI 2014}, pages 103--111, Cham, 2016.
  Springer International Publishing.

\bibitem{cran:ishi:lion:92}
M.~G. Crandall, H.~Ishii, and P.-L. Lions.
\newblock User's guide to viscosity solutions of second order partial
  differential equations.
\newblock {\em Bull. Amer. Math. Soc. (N.S.)}, 27(1):1--67, 1992.

\bibitem{drap:heyn:kupp:13}
S.~Drapeau, G.~Heyne, and M.~Kupper.
\newblock Minimal supersolutions of convex {BSDE}s.
\newblock {\em Ann. Probab.}, 41(6):3973--4001, 2013.

\bibitem{elka:peng:quen:97}
N.~{El Karoui}, S.~Peng, and M.~Quenez.
\newblock Backward stochastic differential equations in finance.
\newblock {\em Math. Finance}, 7(1):1--71, 1997.

\bibitem{fors:kenn:12}
P.~A. Forsyth, J.~S. Kennedy, S.~Tse, and H.~Windcliff.
\newblock Optimal trade execution: a mean quadratic variation approach.
\newblock {\em Journal of Economic Dynamics and Control}, 36(12):1971--1991,
  2012.

\bibitem{friedman}
A.~Friedman.
\newblock {\em Partial differential equations of parabolic type}.
\newblock Dover Publications, 2008.

\bibitem{gath:10}
J.~Gatheral.
\newblock No-dynamic-arbitrage and market impact.
\newblock {\em Quant. Finance}, 10(7):749--759, 2010.

\bibitem{guea:15}
O.~Gu\'{e}ant.
\newblock Optimal execution and block trade pricing: a general framework.
\newblock {\em Appl. Math. Finance}, 22(4):336--365, 2015.

\bibitem{gueant2016financial}
O.~Gu\'{e}ant.
\newblock {\em The financial mathematics of market liquidity}.
\newblock Chapman \& Hall/CRC Financial Mathematics Series. CRC Press, Boca
  Raton, FL, 2016.
\newblock From optimal execution to market making.

\bibitem{krus:popi:14}
T.~Kruse and A.~Popier.
\newblock Bsdes with monotone generator driven by brownian and poisson noises
  in a general filtration.
\newblock {\em Stochastics}, 88(4):491--539, 2016.

\bibitem{krus:popi:15}
T.~Kruse and A.~Popier.
\newblock Minimal supersolutions for {BSDEs} with singular terminal condition
  and application to optimal position targeting.
\newblock {\em Stochastic Processes and their Applications}, 126(9):2554 --
  2592, 2016.

\bibitem{lady:solo:ural:68}
O.~A. Lady{\v{z}}enskaja, V.~A. Solonnikov, and N.~N. Ural{\cprime}ceva.
\newblock {\em Linear and quasilinear equations of parabolic type}.
\newblock Translated from the Russian by S. Smith. Translations of Mathematical
  Monographs, Vol. 23. American Mathematical Society, Providence, R.I., 1968.

\bibitem{nual:06}
D.~Nualart.
\newblock {\em The {M}alliavin calculus and related topics}.
\newblock Probability and its Applications (New York). Springer-Verlag, Berlin,
  second edition, 2006.

\bibitem{pard:rasc:14}
E.~Pardoux and A.~Rascanu.
\newblock {\em {Stochastic Differential Equations, Backward SDEs, Partial
  Differential Equations}}, volume~69 of {\em Stochastic Modelling and Applied
  Probability}.
\newblock Springer-Verlag, 2014.

\bibitem{popi:06}
A.~Popier.
\newblock Backward stochastic differential equations with singular terminal
  condition.
\newblock {\em Stochastic Process. Appl.}, 116(12):2014--2056, 2006.

\bibitem{popi:16}
A.~Popier.
\newblock Limit behaviour of bsde with jumps and with singular terminal
  condition.
\newblock {\em ESAIM: PS}, 20:480--509, 2016.

\bibitem{popi:17}
A.~Popier.
\newblock Integro-partial differential equations with singular terminal
  condition.
\newblock {\em Nonlinear Anal.}, 155:72--96, 2017.

\bibitem{popier2020continuity}
S.~A. Samuel, A.~Popier, and A.~D. Sezer.
\newblock Continuity problem for singular {BSDE} with random terminal time.
\newblock {\em ALEA Lat. Am. J. Probab. Math. Stat.}, 19(2):1185--1220, 2022.

\bibitem{schi:scho:09}
A.~Schied and T.~Sch\"{o}neborn.
\newblock Risk aversion and the dynamics of optimal liquidation strategies in
  illiquid markets.
\newblock {\em Finance Stoch.}, 13(2):181--204, 2009.

\bibitem{schi:scho:tehr:10}
A.~Schied, T.~Sch\"{o}neborn, and M.~Tehranchi.
\newblock Optimal basket liquidation for {CARA} investors is deterministic.
\newblock {\em Appl. Math. Finance}, 17(6):471--489, 2010.

\bibitem{krus:popi:seze:18}
A.~D. Sezer, T.~Kruse, and A.~Popier.
\newblock Backward stochastic differential equations with non-{M}arkovian
  singular terminal values.
\newblock {\em Stoch. Dyn.}, 19(2):1950006, 34, 2019.

\end{thebibliography}
